\documentclass[12pt]{amsart}
\usepackage{latexsym,amscd, amsmath, amssymb,epsfig,xypic,tikz,graphicx} 
\usepackage{bbm}

\theoremstyle{plain}
  \newtheorem{theorem}{Theorem}[section]
  \newtheorem{proposition}[theorem]{Proposition}
  \newtheorem{lemma}[theorem]{Lemma}
  \newtheorem{corollary}[theorem]{Corollary}
  
\theoremstyle{definition}
  \newtheorem{definition}[theorem]{Definition}
  \newtheorem{example}[theorem]{Example}

\theoremstyle{remark}

\numberwithin{equation}{section}

\def\umapright#1{\smash{
   \mathop{\longrightarrow}\limits^{#1}}}

\def\umapleft#1{\smash{
   \mathop{\longleftarrow}\limits^{#1}}}

\def\lmapdown#1{\llap{$\vcenter{\hbox{$\scriptstyle#1$}}$}
    \Big\downarrow}
\def\rmapdown#1{\Big\downarrow\rlap
   {$\vcenter{\hbox{$\scriptstyle#1$}}$}}

\def\rmapup#1{\Big\uparrow\rlap
   {$\vcenter{\hbox{$\scriptstyle#1$}}$}}

\def\tempbaselines
{\baselineskip22pt\lineskip3pt
   \lineskiplimit3pt}
\def\diagram#1{\null\,\vcenter{\tempbaselines
\mathsurround=0pt
    \ialign{\hfil$##$\hfil&&\quad\hfil$##$\hfil\crcr
      \mathstrut\crcr\noalign{\kern-\baselineskip}
  #1\crcr\mathstrut\crcr\noalign{\kern-\baselineskip}}}\,}

\def\pullback#1&#2&#3&#4&#5&#6&#7&#8&{
\diagram{#1&\umapright{#2}&#3\cr
\rmapdown{#4}&&\rmapdown{#5}\cr
#6&\umapright{#7}&#8\cr}}

\def\calA{{\mathcal A}}
\def\calB{{\mathcal B}}
\def\calC{{\mathcal C}}
\def\calD{{\mathcal D}}
\def\calE{{\mathcal E}}
\def\calF{{\mathcal F}}
\def\calG{{\mathcal G}}
\def\calH{{\mathcal H}}

\def\calJ{{\mathcal J}}
\def\calK{{\mathcal K}}

\def\calP{{\mathcal P}}
\def\calQ{{\mathcal Q}}
\def\calS{{\mathcal S}}
\def\calU{{\mathcal U}}
\def\calV{{\mathcal V}}
\def\calX{{\mathcal X}}
\def\calY{{\mathcal Y}}

\def \all{\mathop{\rm all}\nolimits} 
\def \Aut{\mathop{\rm Aut}\nolimits} 
\def\Cat{{\mathop{\rm Cat}\nolimits}}
\def\CSet{{\mathop{\calC{\rm-Set}}\nolimits}}

\def\Cores{{\mathop{\rm Cores}\nolimits}}
\def\Corresp{{\mathop{\rm Corresp}\nolimits}}
\def\Cograph{{\mathop{\rm Cograph}\nolimits}}

\def \Dist{\mathop{\rm Dist}\nolimits}

\def \End{\mathop{\rm End}\nolimits} 
\def \Ess{\mathop{\rm Ess}\nolimits} 
\def \Ext{\mathop{\rm Ext}\nolimits} 

\def \Fun{{\mathop{\rm Fun}\nolimits} }

\def \Hom{\mathop{\rm Hom}\nolimits} 
\def \Im{\mathop{\rm Im}\nolimits} 
\def \Iness{\mathop{\rm Iness}\nolimits}

\def \Ker{\mathop{\rm Ker}\nolimits}
\def \LK{{\mathop{\rm LK}\nolimits}}
\def \Mat{{\mathop{\rm Mat}\nolimits}}
\def \mod{\hbox{\rm -mod}}
\def \Ob{{\mathop{\rm Ob}\nolimits}}
\def\one{{{\rm 1} \kern -.32em {\rm 1} }}
\def\one{{\mathbbm{1} }}
\def\one{{\mathbf{1} }}
\def \op{{\mathop{\rm op}\nolimits}}
\def \Out{\mathop{\rm Out}\nolimits} 
\def \Rad{\mathop{\rm Rad}\nolimits} 
\def \Rank{\mathop{\rm Rank}\nolimits} 
\def \Res{\mathop{\rm Res}\nolimits} 
\def\Set{\mathop{\rm Set}\nolimits} 
 
\def \Tor{\mathop{\rm Tor}\nolimits} 
\def \ydp{\mathsf{p}}
\def \ydP{\mathsf{P}}

\def\BB{{\Bbb B}}

\def\NN{{\Bbb N}}

\def\SS{{\Bbb S}}

\def\ZZ{{\Bbb Z}}

\begin{document}

\title[Biset functors for categories]
{Biset functors for categories}

\author{Peter Webb}
\email{webb@umn.edu}
\address{School of Mathematics\\
University of Minnesota\\
Minneapolis, MN 55455, USA}



\begin{abstract}
We introduce the theory of biset functors defined on finite categories. Previously, biset functors have been defined on groups, and in that context they are closely related to Mackey functors. Standard examples on groups include representation rings, the Burnside ring and group cohomology. The new theory allows these same examples, but defined for arbitrary finite categories, thus including, for instance, the representation rings of finite EI categories, among which are posets, and the free categories associated to quivers without oriented cycles. Group homology with trivial coefficients is replaced by the homology of the category with constant coefficients and this is a biset functor when bisets that are representable on one side are used. We give a definition of the Burnside ring of an arbitrary finite category. It is a biset functor that plays a key role throughout the theory. We discuss properties of the simple biset functors on categories, including their parametrization and calculation. We describe the symmetric monoidal structure on biset functors, Green biset functors and an approach to fibered biset functors for categories, together with the technicalities these entail.  Various examples are given, the most elaborate showing a connection between the correspondence functors of Bouc and Th\'evenaz and biset functors on Boolean lattices.
\end{abstract}

\maketitle
\tableofcontents

\section{Introduction}

We introduce the study of biset functors defined on finite categories. Biset functors defined on groups have been studied for a long time and include examples such as Grothendieck groups of group rings, the Burnside ring of a group, and group cohomology (with a restriction on the bisets used), among others. They are closely related to Mackey functors. These structures have been used as an aid in calculating values of these Grothendieck groups, group cohomology, and so on. The fundamental nature of the constructions also provides a framework in which questions about group representation theory may be phrased and studied.  For two overviews of such functors on groups see \cite{Bou3} and \cite{Webf} (where biset functors are called `globally-defined Mackey functors').

The goals of introducing biset functors on finite categories are similar, but the context is broader. The new theory applies to Grothendieck rings of representations of categories, thus including representations of quivers and posets, as well as representations of groups. We make a new definition of the Burnside ring of every finite category, generalizing the definition for groups. 

The theory applies also to the homology and cohomology of categories. In the context of finite groups, a restriction must be made to use only bisets that are free on one side if we wish to regard cohomology and homology as biset functors. We identify the generalization of this condition that works for categories, and it is that the bisets be representable on one side. With this restriction, homology and cohomology of categories are biset functors, as is Hochschild homology. One use of this is that the biset functor formalism provides a way to construct transfer maps.

As a special case, any simplicial complex determines a category, namely, the poset of simplices of the simplicial complex. The  homology of the poset in the sense of category homology is the same as the topological homology of the simplicial complex. By this means the usual topological homology of finite simplicial complexes is included within the theory of biset functors. 

Before presenting these examples we must first establish what we mean by a biset functor in the context of categories, and we begin by setting up the foundations of the theory. The development is parallel to that of biset functors for finite groups, and we will see that many things that work for finite groups carry straight over to finite categories, while other things do not work. 

In Section 2 we start with the notion of a set with an action of a category, or `category set'. By this we mean a functor from the category to the category of sets, a construction that appears throughout mathematics in various forms. For example, regarding groups as categories with a single object in which all morphisms are invertible, a category set for a group $G$  is the same thing as a $G$-set, or permutation representation.  There are many other examples of category sets that standardly arise.
For instance, simplicial sets are defined as functors from a certain category to sets; in combinatorics, directed graphs and rooted trees may both be regarded as sets for the action of certain categories; sometimes (contravariant) functors to the category of sets are termed presheaves; and so on. In developing the theory of category sets along the lines of permutation representations we find that we have to be careful about certain things. For example, there is more than one candidate for the definition of a transitive set, and the different candidates have different properties. It is in this section that we define the Burnside ring of a finite category to be the Grothendieck ring of its finite category sets.

Section 3 is about the first properties of bisets and biset functors. A biset for a pair of categories $\calC,\calD$ is a set for $\calC\times\calD^\op$. Such bisets have been studied previously, in a different context, under the name of distributors, or profunctors. We import results established for bisets in that other context whenever possible.  This includes the associativity of the product that we define on bisets, as well as other properties. We define the biset category: the objects are finite categories and the morphisms are linear combinations of bisets. Composition of morphisms is determined by the biset product and distributivity. 

There are features of bisets for categories that do not appear when we only consider bisets for groups. For instance, non-isomorphic categories may become isomorphic in the biset category, a possibility that does not happen with groups. There is, in general, no factorization of bisets as composites of special bisets that factor through smaller categories, as happens with groups. On the other hand, we can factor every biset in a different way as a composite of two bisets factoring through a larger category. This is a possibility not allowed when the objects in the biset category are just groups, as it requires categories to have more than one object. These and other properties are studied in Section 4.

In Section 5 we introduce representable sets and bisets, and in Section 6 we show how to define the homology and cohomology of a category as biset functors, using bisets that are representable on one side. 
In the case of biset functors defined on groups, the usual way to see that group cohomology is a biset functor is to observe that there are generating bisets corresponding to restriction, transfer and conjugation that satisfy defining relations also satisfied by cohomology. In the case of bisets and cohomology for categories we no longer have distinguished generating bisets, let alone relations they might satisfy. Our approach instead is to construct all the cohomology operations at once. We use the fact that this has already been done for Hochschild homology, which also has the structure of a biset functor on categories. We deduce the result for ordinary category homology by splitting it off from Hochschild homology in the same way as was done by Xu for ordinary cohomology and Hochschild cohomology. In the special case of group homology and cohomology, this process constructs restriction, transfer and conjugation maps by a single formula. For finite categories in general it constructs transfer maps in situations where induction need not be both the left and right adjoint of restriction (as is usually assumed), as we show by example

Treating biset functors in the spirit of representation theory, it is crucial to identify the simple biset functors, and we provide a way to do this in Section 7. There is a standard well-known relationship between simple functors on a category and simple functors on a full subcategory. It means that simple functors are either simple or zero on each full subcategory, and also that each simple functor on a subcategory extends uniquely to a simple functor on the larger category. As a consequence, simple functors are determined by each of their non-zero evaluations at categories, as modules for the endomorphism rings (in the biset category) of those categories. Parametrizing a simple functor is then a question of choosing a preferred category on which it is non-zero, and finding an appropriate way to describe its simple evaluation there. This is the basis of our approach to the parametrization of simple biset functors. Finding the appropriate description of simple modules for endomorphism rings of categories turns out to be not so straightforward as for groups.

As a more extended example, we make a connection in Section 8 between the correspondence functors studied by Bouc and Th\'evenaz and biset functors defined on Boolean lattices of subsets of some given set, using representable bisets. We show that there is a canonical bijection between the simple functors in the two cases. Thus each simple correspondence functor determines a simple biset functor and, conversely, each simple biset functor that is non-zero on a Boolean lattice determines a simple correspondence functor.

In Section 9 we describe further structures on the biset category and on biset functors. We show that the biset category is a rigid symmetric monoidal linear category with an internal Hom, implying that biset functors also have a symmetric monoidal structure, again with an internal Hom. The rigidity of the biset category translates to certain isomorphisms between biset functors. From all this, Green biset functors and modules are defined, extending the definition for groups. We conclude with a generalization of fibered biset functors to categories.

Some notation: we will write $\Set$ for the category of \textit{finite} sets, and $R$ will be a commutative ring. All arbitrary categories we consider, denoted by letters $\calC, \calD$ etc, will be \textit{finite}, and we will assume this without explicitly mentioning it. We write $\calC^\op$ for the opposite category to $\calC$ and $\one$ for the category with a single object and a single morphism (the identity group). We let $\Cat$ denote the category whose objects are finite categories, and whose morphisms are functors.

I express my thanks to the many people who have helped me with their knowledge and insight, and especially I thank Robert Boltje, Serge Bouc, Andrew Snowden, Ben Spitz and Jacques Th\'evenaz  for very valuable discussions.

\section{Sets with an action of a category}

Permutation actions of groups, or $G$-sets, are fundamental in the study of groups and in any situation where group actions are considered. We generalize this notion to categories. A preliminary version of some of this material was given in \cite{Webe}.

\begin{definition}
Given a category $\calC$, a  \textit{left $\calC$-set} is a functor $\Omega:\calC\to\Set$. By a \textit{right $\calC$-set} we mean a functor $\Psi:\calC^\op\to\Set$, where $\calC^\op$ is the opposite category. Without the specification of left or right, a $\calC$-set is taken to be a left $\calC$-set. The  $\calC$-sets form a category $\calC\hbox{-Set}$ whose objects are the $\calC$-sets and in which the morphisms are the natural transformations between the functors. 
\end{definition}

Thus a $\calC$-set means that for each object $x$ of $\calC$ we have a set $\Omega(x)$, and for each morphism $\alpha:x\to y$ in $\calC$ we have a map of sets $\Omega(\alpha):\Omega(x)\to\Omega(y)$ satisfying the usual functorial conditions. Rather than writing $\Omega(\alpha)(u)$ for the image of an element $u\in\Omega(x)$, we often write simply $\alpha u$. The \textit{left} $\calC$-set terminology is consistent with the convention of applying mappings from the left, so that $\beta(\alpha u) = (\beta\alpha)u$ when $\beta: y\to z$. Given a contravariant functor on $\calC$, meaning a functor $\Psi:\calC^\op\to\Set$, we may write $v\beta$ for the element $\Psi(\beta)(v)$ where $v\in\Psi(z)$. Now $v(\beta\alpha)=(v\beta)\alpha$.

When $\calC$ is a group (a category with one object in which all morphisms are invertible) then a $\calC$-set is the same thing as a $G$-set, where $G$ is the set of morphisms of $\calC$.   In general a $\calC$-set $\Omega$ is just a diagram of sets having the shape of $\calC$. Many examples of $\calC$-sets for other categories $\calC$ arise in a fundamental way throughout mathematics.
We continue with basic definitions.

\begin{definition}
A \textit{sub-$\calC$-set} of a $\calC$-set $\Omega$ is a functor $\Psi:\calC\to\Set$ for which each set $\Psi(x)$ is a subset of $\Omega(x)$ and the morphisms $\Psi(\alpha)$ are restrictions of the morphisms $\Omega(\alpha)$. If $\Psi_1,\Psi_2$ are both sub-$\calC$-sets of $\Omega$ then $\Psi_1\cap\Psi_2$ is the sub-$\calC$-set with $(\Psi_1\cap\Psi_2)(x)=\Psi_1(x)\cap\Psi_2(x)$ for all objects $x$.
Given two $\calC$-sets $\Omega_1$ and $\Omega_2$ we define their \textit{disjoint union} $\Omega_1\sqcup\Omega_2$ to be the $\calC$-set defined at each object $x$ of $\calC$ by
$$
(\Omega_1\sqcup\Omega_2)(x):=\Omega_1(x)\sqcup\Omega_2(x)
$$
 and on morphisms $\alpha:x\to y$ in $\calC$ 
 $$
 (\Omega_1\sqcup\Omega_2)(\alpha)(u):=
\begin{cases}\Omega_1(\alpha)(u)&\hbox{if }u \in \Omega_1(x)\\
\Omega_2(\alpha)(u)&\hbox{if }u \in \Omega_2(x).\\
\end{cases}
 $$
We write $\Omega=\Omega_1\sqcup\Omega_2$ to mean that $\Omega(x)=\Omega_1(x)\sqcup\Omega_2(x)$  for all objects $x$ of $\calC$. We will call a non-empty $\calC$-set \textit{indecomposable}  if it is not isomorphic to a proper disjoint union of $\calC$-sets (one where both $\calC$-sets are non-empty). 

A $\calC$-set $\Omega$ is \textit{finite} if there are only finitely many objects on which $\Omega$ is non-empty, and on each of those objects $\Omega$ takes value a finite set. 

We define the \textit{support} of a $\calC$-set $\Omega$ to be the set of objects $x$ in $\calC$ for which $\Omega(x)$ is non-empty. We will say that a category $\calC$ is \textit{connected} if it cannot be written $\calC=\calC_1\sqcup\calC_2$ where $\calC_1$ and $\calC_2$ are non-empty categories and there are no morphisms going from one of $\calC_1,\calC_2$ to the other. 
\end{definition}

An equivalent condition for a category $\calC$ to be connected is that for every pair of objects $a,b$ in $\calC$ there is a chain of morphisms in $\calC$ of the form $a \to x_1 \gets x_2 \to x_3 \cdots \gets x_n\to b$. It is also equivalent to say that the nerve of $\calC$ is connected.

When $G$ is a group, transitive $G$-sets or orbits can be defined in more than one way: as sets generated by a single element; or as indecomposable $G$-sets. In the generality of categories it turns out that the definition as an indecomposable set works better than the other one. Thus the decomposition of any permutation representation of a group uniquely as a disjoint union of transitive sets works for categories as a disjoint union of indecomposable sets.

\begin{proposition}
Let $\calC$ be a category.
\begin{enumerate}
\item If $\Psi$ is a sub-$\calC$-set of $\Omega$ and $\Omega=\Omega_1\sqcup\Omega_2$ then 
$$
\Psi=(\Psi\cap\Omega_1)\sqcup(\Psi\cap\Omega_2).
$$
\item Every finite $\calC$-set $\Omega$ has a unique decomposition
$$\Omega=\Omega_1\sqcup \Omega_2\sqcup\cdots\sqcup\Omega_n$$
where each $\Omega_i$ is indecomposable. 
\item If $\Omega$ is an indecomposable $\calC$-set then the full subcategory of $\calC$ whose objects are the support of $\Omega$ is connected.
\end{enumerate}
\end{proposition}

\begin{proof}
(1) This is because for any set of the form $A=A_1\sqcup A_2$ and any subset $B\subseteq A$ we have $B=(B\cap A_1)\sqcup (B\cap A_2)$. 

(2) The finite $\calC$-set $\Omega$ may happen to be the disjoint union of two $\calC$-sets, or not; if it can be broken up as a disjoint union we can ask if either of the factors is a disjoint union, and by repeating this we end up with a disjoint union of $\calC$-sets each of which is indecomposable. Uniqueness of the decomposition follows from (1).

(3) Let $\calD$ be the full subcategory of $\calC$ whose objects are the support of $\Omega$. If $\calD$ is disconnected then $\calD=\calD_1\sqcup\calD_2$ with no morphisms between $\calD_1$ and $\calD_2$. In such a situation we may write $\Omega=\Omega_1\sqcup\Omega_2$, where for $i=1,2$, $\Omega_i(x)=\Omega(x)$ if $x\in\calD_i$, and is empty otherwise and thus one of  $\Omega_1,\Omega_2$ must be empty. If, say, $\Omega_1$ is empty then $\calD_1$ is not in the support of $\Omega$ which is a contradiction. Thus $\calD$ is connected.
\end{proof}

\begin{definition}
Given a $\calC$-set $\Omega$ and a subset
$$
X\subseteq \bigsqcup_{x\in\Ob\calC}\Omega(x),
$$
we say that the sub-$\calC$-set  of $\Omega$ \textit{generated} by $X$ is the intersection of the sub-$\calC$-sets of $\Omega$ that contain $X$. We say that $\Omega$ itself is generated by $X$ if it equals the sub-$\calC$-set of $\Omega$ generated by $X$. 

Furthermore, if $u\in \bigsqcup_{x\in\Ob\calC}\Omega(x)$ we define the \textit{orbit} of $u$ to be the unique indecomposable summand $\Omega_i$ of $\Omega$ for which $u\in\Omega_i(x)$
\end{definition}

In the case of sets for a group, the indecomposable subsets are precisely the sets generated by one element, and these are also the orbits of a $G$-set. The same is not true in general for $\calC$-sets, as the next example shows. It illustrates the care we must take in defining these concepts for categories.

\begin{example}
\label{a2-sets}
Let $\calC$ be the category 
$$\calA_2 =\mathop{\bullet}\limits_x\xrightarrow{\alpha}\mathop{\bullet}\limits_y$$
that has two objects $x$ and $y$, a single morphism $\alpha$ from $x$ to $y$, and the identity morphisms at $x$ and $y$. Let $\Omega$ be a $\calC$-set. By considering the points of $\Omega(x)$ that map to each point of $\Omega(y)$, we see that the indecomposable $\calA_2$-sets have the form
$$
\Omega_n:= \underline n \to \underline 1,\quad n\ge 0
$$
where $\underline n = \{1,\ldots,n\}$ is a set with $n$ elements, the mapping between the two sets sending every element onto a single element. We see from this example that a finite category may have infinitely many non-isomorphic indecomposable sets, and also that indecomposable sets need not be generated by any single element, because $\Omega_n$ needs $n$ elements to generate it. The orbits of a $\calC$-set such as $\Omega_5\sqcup\Omega_7$ are  precisely the sets $\Omega_5,\Omega_7$, even though there are other indecomposable sets such as $\Omega_0,\Omega_1$ etc. as subsets. There is also an indecomposable infinite set $\Omega_\infty$ defined with $\Omega_\infty (x)$ of infinite cardinality.

More generally, let $\calA_n=\mathop{\bullet}\limits_{x_n}\to \mathop{\bullet}\limits_{x_{n-1}}\to \cdots\to\mathop{\bullet}\limits_{x_1}$ be the category that is the poset with $n$ elements in a chain. The indecomposable $\calA_n$-sets are the rooted trees whose nodes are all at distance $\le n-1$ from the root. 
\end{example}

We do have the following relation between indecomposability and generation by a single element.

\begin{lemma}
\label{cyclic-implies-indecomposable}
Let $\Omega$ be a $\calC$-set. If there exists an object $x$ and an element $u\in\Omega(x)$ that generates $\Omega$ then $\Omega$ is indecomposable.
\end{lemma}

\begin{proof}
If $\Omega=\Omega_1\sqcup \Omega_2$ then without loss of generality $u\in\Omega_1(x)$ and we deduce that the sub-$\calC$-set of $\Omega$ generated by $u$ must be contained in $\Omega_1$. It follows that $\Omega=\Omega_1$ so that $\Omega_2$ is empty, and $\Omega$ must be indecomposable.
\end{proof}

\begin{example}
\label{Kronecker-example}
Let $\calC$  be the \textit{Kronecker category} with two objects $x$ and $y$, and two morphisms $\alpha,\beta:x\to y$. This finite category has infinitely many indecomposable finite sets, all of whose morphisms are injective maps. To see this, consider the $\calC$-sets $\Omega_n$ and $\Psi_n$, where $n\ge 1$, described diagrammatically as follows.
$$
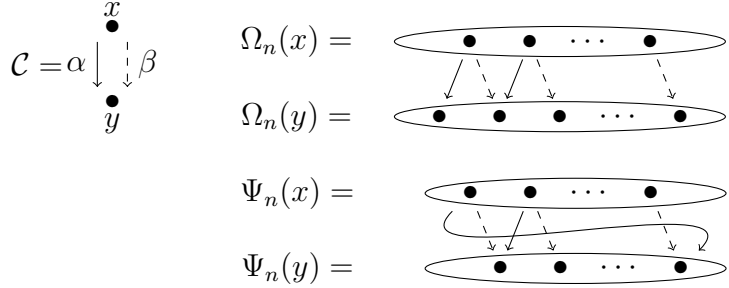
\begin{aligned}
\begin{tikzpicture}[xscale=1,yscale=1]
\node at (-1,0.5) {$\calC =$};
\node (A) at (0, 1) {$\bullet$};
\node (B) at (0,0) {$\bullet$};
\draw [->] (-0.2,0.8)--(-0.2,0.2) node[midway,left]{$\alpha$};
\draw [->, dashed] (0.2,0.8)--(0.2,0.2)  node[midway,right]{$\beta$};
\node [ above] at  (A)   {$x$};
\node [ below] at  (B)   {$y$};
\end{tikzpicture}
\qquad
&\begin{tikzpicture}[xscale=1,yscale=1]
\node at (-2.3,1) {$\Omega_n(x) =$};
\node at (-2.3,0) {$\Omega_n(y) =$};
\node (A) at (0, 1) {$\bullet$};
\node (B) at (0.8,1) {$\bullet$};
\node (C) at (1.6,1) {$\cdots$};
\node (E) at (0.4,0) {$\bullet$};
\node (F) at (1.2,0) {$\bullet$};
\node (G) at (2.0,0) {$\cdots$};
\node (I) at (2.8,0) {$\bullet$};
\node (H) at (2.4,1) {$\bullet$};
\node (J) at (-0.4,0) {$\bullet$};
\draw [->] (B)--(E);
\draw [->, dashed] (A)--(E);
\draw [->, dashed] (B)--(F);
\draw [->, dashed] (H)--(I);
\draw [->] (A)--(J);
\draw (1.2,1) ellipse (2.2cm and 0.2cm);
\draw (1.2,0) ellipse (2.2cm and 0.2cm);
\end{tikzpicture}
\cr
&\begin{tikzpicture}[xscale=1,yscale=1]
\node at (-2.3,1) {$\Psi_n(x) =$};
\node at (-2.3,0) {$\Psi_n(y) =$};
\node (A) at (0, 1) {$\bullet$};
\node (B) at (0.8,1) {$\bullet$};
\node (C) at (1.6,1) {$\cdots$};
\node (E) at (0.4,0) {$\bullet$};
\node (F) at (1.2,0) {$\bullet$};
\node (G) at (2.0,0) {$\cdots$};
\node (I) at (2.8,0) {$\bullet$};
\node (H) at (2.4,1) {$\bullet$};
\draw [->] (B)--(E);
\draw [->, dashed] (A)--(E);
\draw [->, dashed] (B)--(F);
\draw [->, dashed] (H)--(I);
\draw [->] (A) to [out=225,in=45] (I);
\draw (1.4,1) ellipse (2cm and 0.2cm);
\draw (1.4,0) ellipse (2cm and 0.2cm);
\end{tikzpicture}
\cr
\end{aligned}$$
These $\calC$-sets have $|\Omega_n(x)|=|\Psi_n(x)|=|\Psi_n(y)|=n$ and $|\Omega_n(y)|=n+1$. There are also  sets of the form $\Omega_\infty$ that continue indefinitely in one or both directions.
The effects of $\alpha$ and $\beta$ are indicated by solid and dashed arrows. The $\calC$-sets of this kind are all of the indecomposable sets with morphisms acting as injective maps, the other $\calC$-sets having maps are not injective.

The $\calC$-sets for this category are in bijection with directed graphs: given a $\calC$-set $\Omega$ we may construct a graph whose vertex set is $\Omega(y)$, and where there is an edge from vertex $a$ to vertex $b$ if there is an element $u\in\Omega(x)$ so that  $\alpha(u)=a$ and $\beta(u)=b$. Every directed graph arises in this way. The connected directed graphs coming from $\calC$-sets $\Theta$ where both $\Theta(\alpha)$ and $\Theta(\beta)$ are injective are the oriented chains and circuits.
\end{example}

We now give a definition of the Burnside ring of a finite category, generalizing the definition of the Burnside ring of a finite group.  Our definition is different to other definitions we may find, such as those given in \cite{May, Yos, Yos2}.

\begin{definition}
\label{burnside-ring-definition}
We define the \textit{Burnside ring} of a finite category $\calC$ to be
$$
B(\calC):=\hbox{Grothendieck group of finite $\calC$-sets with respect to $\sqcup$}.
$$
Thus $B(\calC)$ is the free abelian group with the (isomorphism classes of) indecomposable $\calC$-sets as a basis.

There is a product operation on $\calC\hbox{-Set}$ making this a symmetric monoidal category: given two $\calC$-sets $\Omega$ and $\Psi$ we define $(\Omega\times\Psi)(x)=\Omega(x)\times\Psi(x)$ for every object $x$, and $(\Omega\times\Psi)(\alpha)=\Omega(\alpha)\times\Psi(\alpha)$ on morphisms $\alpha$ of $\calC$.  This allows us to define a multiplication on $B(\calC)$, given by $\times$ on the basis elements.

We let $*$ denote the constant $\calC$-set whose value is a single point. This is a functor $*:\calC\to\Set$ whose value at each object is a single point, and whose value on every morphism is the identity morphism. Note that $*$ is indecomposable if and only if $\calC$ is connected. 
\end{definition}

The following is evident:

\begin{proposition}
Let $\calC$ be a finite category. Then $B(\calC)$ is a commutative ring with identity $*$, and with basis the isomorphism classes of indecomposable $\calC$-sets.
\end{proposition}

The literature on the Burnside rings of groups is well established, so we give an example of a Burnside ring of a category that is not a group.

\begin{example}
As in Example~\ref{a2-sets} consider the category $$\calA_2 =\mathop{\bullet}\limits_x\xrightarrow{\alpha}\mathop{\bullet}\limits_y.$$
We have seen that its indecomposable sets have the form $\Omega_n:= \underline n \to \underline 1,\quad n\ge 0
$ and we see that $\Omega_m\times\Omega_n \cong \Omega_{mn}$. Thus $B(\calA_2)\cong \ZZ\NN^\times$ is the monoid algebra over $\ZZ$ of the monoid $\NN^\times$ of natural numbers under multiplication.
\end{example}

When $G$ is a group we can recover $G$ from its category of $G$-sets. More generally, we can only recover a category $\calC$ from its category of $\calC$-sets up to idempotent completion, as we shall see.  When $\calC$ is a finite category we take its \textit{idempotent completion} to be the category $\bar\calC$ whose objects are the idempotent endomorphisms in $\calC$, and where if $e\in\End_\calC(x)$ and $f\in\End_\calC(y)$ are idempotent then $\Hom_{\bar\calC}(e,f):=f\cdot \Hom_\calC(x,y)\cdot e $. This is called the Cauchy completion in \cite{Bor} and it is also known as the Karoubian envelope.

The following result appears as Theorem 6.5.11 in \cite{Bor}, where a proof is given.

\begin{theorem}
Categories $\calC$ and $\calD$ have equivalent idempotent completions if and only if  the categories $\calC\hbox{-Set}$ and $\calD\hbox{-Set}$ are naturally equivalent.
\end{theorem}

In particular, equivalent categories $\calC\simeq\calD$ have $\calC\hbox{-Set}\simeq\calD\hbox{-Set}$, and composition with an equivalence $F:\calC\to\calD$ provides an equivalence $\calD\hbox{-Set}\to\calC\hbox{-Set}$. Also, composition with the natural inclusion $\calC\to\bar\calC$ provides an equivalence $\bar\calC\hbox{-Set}\to\calC\hbox{-Set}$.  We may deduce from this the following. 

\begin{theorem}
\label{idempotent-completion-Burnside-theorem}
If $\calC$ and $\calD$ have equivalent idempotent completions then the Burnside rings $B(\calC)$ and $B(\calD)$ are isomorphic as rings.
\end{theorem}

\begin{proof}
Suppose that $\calC$ and $\calD$ have equivalent idempotent completions, so that $\calC\hbox{-Set}\simeq\calD\hbox{-Set}$. Disjoint union decompositions of $\calC$-sets are determined as coproducts, so that  indecomposable $\calC$-sets correspond to indecomposable $\calD$-sets under the equivalence. This means that any such equivalence induces an isomorphism of abelian groups $B(\calC)\cong B(\calD)$.

To show that we have a ring isomorphism, let us first observe that if we have an inclusion of categories $I:\calU\to\calV$ that induces an equivalence $\calV\hbox{-Set}\to\calU\hbox{-Set}$ then we have a ring isomorphism $B(\calU)\cong B(\calV)$. This is because the map on Burnside rings induced by an inclusion is a ring homomorphism. 

To apply this, the equivalence of $\bar\calC$ and $\bar\calD$ means there is a common skeletal subcategory $\calS$ of $\bar\calC$ and of $\bar\calD$ so that the inclusions $\calS\to\bar\calC$ and $\calS\to\bar\calD$ are equivalences, and now all of the inclusions $\calC\to\bar\calC$ and $\calS\to\bar\calC$, $\calS\to\bar\calD$ and $\calD\to\bar\calD$ induce equivalences of the categories of sets. The result follows by composing the ring isomorphisms obtained from this.
\end{proof}

\begin{example}
The Burnside ring of the monoid $\{1, \alpha\}$ where $\alpha^2=\alpha$ is isomorphic to the Burnside ring of the category shown in Example~\ref{monoid-idempotent-completion}, which is its idempotent completion.
\end{example}

\section{Bisets for categories}
Bisets have been studied for a long time under several different names. In the context of category theory they are known as \textit{distributors}, or \textit{profunctors}. They appear in the 1966 thesis of Bunge \cite{Bun} as well as in the work of B\'enabou \cite{Ben} from 1973, and there is a description of the theory of distributors in \cite[Sec. 7.8]{Bor}. When the categories on which they are defined are groups, bisets and a biset category using bisets as morphisms were introduced in a topological context in \cite{AGM}. The bisets there had a free action on one side, The similar category without this restriction was introduced in \cite{Bou1}.  We review these definitions and introduce our notation, using the term biset, rather distributor or profunctor, noting that the theory of bisets differs from that of distributors once we impose the relation that disjoint union corresponds to sum.

\begin{definition}
Given finite categories $\calC$ and $\calD$ we define a \textit{$(\calC,\calD)$-biset} to be a $\calC\times\calD^\op$-set; that is, a functor $\Omega:\calC\times\calD^\op\to\Set$. To record the presence of $\calC$ and $\calD$ we will also write ${}_\calC\Omega_\calD$, instead of just  $\Omega$, for this biset. 
\end{definition}

Given $(\calC,\calD)$-biset $\Omega$, a morphism $\alpha:x\to x_1$ in $\calC$, a morphism $\beta:y_1\to y$ in $\calD$, and an element $u\in\Omega(x,y)$ we get elements
$$
\alpha u:=\Omega(\alpha\times 1_y)(u)\in\Omega(x_1,y)
$$
and 
$$
u\beta:=\Omega(1_x\times\beta)(u)\in\Omega(x,y_1).
$$
Thus we have commuting actions of $\calC$ from the left and $\calD$ from the right on $\Omega$. Often it is easier to use the simpler notation on the left side of these equations rather than the full functorial notation on the right.

As with any sets for a category, given $(\calC,\calD)$-bisets ${}_\calC\Omega_\calD$ and ${}_\calC\Psi_\calD$ we may form their disjoint union ${}_\calC\Omega_\calD \sqcup {}_\calC\Psi_\calD$, and the properties of this construction established in the last section hold.

There is a \textit{product}, or \textit{composition}, of bisets that appears in \cite{Ben}, and is also described in  \cite[Prop. 7.8.2]{Bor}. It is as follows.

\begin{definition}
\label{biset-composition-definition}
Given a $(\calC,\calD)$-biset ${}_\calC\Omega_\calD$ and a $(\calD,\calE)$-biset ${}_\calD\Psi_\calE$ we construct a $(\calC,\calE)$-biset $\Omega\circ\Psi$ by the formula, for $x\in\Ob\calC$ and $z\in\Ob\calE$,
$$
\Omega\circ\Psi(x,z)=\left(\bigsqcup_{y\in\Ob\calD} \Omega(x,y)\times\Psi(y,z)\right)/\sim
$$
where $\sim$ is the equivalence relation generated by  $(u\beta,v)\sim(u,\beta v)$ whenever $u\in\Omega(x,y_1)$, $v\in\Psi(y_2,z)$ and $\beta:y_2\to y_1$ in $\calD$. The left functorial action of $\calC$ comes from the action of $\calC$ on $\Omega$ and the right action of $\calE$ comes from the right action of $\calE$ on $\Omega$.
\end{definition}{}

The expression on the right in the above definition is the coend $\int^y \Omega(x,y)\times\Psi(y,z)$ as described in \cite{MacLane}.
 When the categories $\calC$, $\calD$ and $\calE$ are groups, this composition of bisets coincides with the usual definition, given in \cite{AGM, Bou1}. The composition is associative up to isomorphism of bisets, as indicated in the next result. Of special importance is the \textit{identity $(\calC,\calC)$-biset} denoted ${}_\calC\calC_\calC$, for each category $\calC$, whose value at objects $x,y$ in $\calC$ is $\Hom_\calC(y,x)$.

\begin{proposition}[7.8.2 of \cite{Bor}]
\label{Dist-definition}
Small categories, with $(\calC,\calD)$-bisets as morphisms from $\calD$ to $\calC$, and with natural transformations of bisets as 2-morphisms, form a bicategory $\Dist$. In particular, the composition operation $\circ$ is associative up to isomorphism of bisets. For each category $\calC$ the biset ${}_\calC\calC_\calC$ is the identity endomorphism.
\end{proposition}

\begin{proof} We only sketch the proof because a full account appears in \cite{Bor}.
Equivalence classes in $\bigsqcup_{y\in\Ob\calD} \Omega(x,y)\times\Psi(y,z)$ are preserved by the functorial action of $\calC$ from the left and $\calE$ from the right. For, if $\alpha:x\to x_1$ in $\calC$ and $\gamma:z_1\to z$ in $\calE$ then $(\alpha u\beta,v)\sim (\alpha u,\beta v)$ and $(u\beta,v\gamma)\sim (u,\beta v\gamma)$. Hence the definition of $\Omega\circ\Psi$ does indeed give a biset. The operation is associative because if $\Theta$ is a $(\calB,\calC)$-biset then
$(\Theta\circ(\Omega\circ\Psi))(w,z)$ may be identified with the set of equivalence classes of triples $(t,u,v)\in \Theta(w,x)\times\Omega(x,y)\times\Psi(y,z)$ under the equivalence relation defined by requiring that $(t\alpha,u,v)\sim(t,\alpha u,v)$ and $(t,u\beta, v)\sim (t,u,\beta v)$ for all morphisms $\alpha$ in $\calC$ and $\beta$ in $\calD$. The same is true of $((\Theta\circ\Omega)\circ\Psi)(w,z)$, and the functorial actions of $\calB$ and $\calE$ are the same under these identifications.

To see that ${}_\calC\calC_\calC$ acts as the identity under $\circ$, we have a natural transformation ${}_\calC\calC_\calC\circ \Omega\to\Omega$ given at the object $(x,z)\in \calC\times\calD^\op$ by the map specified at the level of elements as $(\alpha,r)\mapsto \alpha r$, where $\alpha: y\to x$ and $r\in\Omega(y,z)$. This map is evidently surjective to $\Omega(x,z)$. If $(\alpha,r)$ and $(\beta,s)$ have the same image then $\alpha r=\beta s$, so that $(\alpha,r)\sim (1_x,\alpha r)=(1_x,\beta s)\sim (\beta,s)$. This shows that the map is bijective and hence we have a natural isomorphism. The argument that $\Omega\circ {}_\calD\calD_\calD\cong\Omega$ as bisets is similar.
\end{proof}

The composition of bisets is a precursor of the tensor product of bimodules, as indicated in the next proposition. Given a commutative ring $R$, for each finite category $\calC$ we form its category algebra $R\calC$ (see \cite{Hai} or \cite{Webb}),  and for each $(\calC,\calD)$-biset $\Omega$ we obtain by linearization an $(R\calC,R\calD)$-bimodule $R\Omega$.

\begin{proposition}
Let  ${}_\calC\Omega_\calD$ be a $(\calC,\calD)$-biset and ${}_\calD\Psi_\calE$ a$(\calD,\calE)$-biset. Passing to bimodules $R\Omega$ and $R\Psi$ we have 
$$
R(\Omega\circ\Psi) \cong R\Omega\otimes_{R\calD} R\Psi
$$
as $(R\calC,R\calE)$-bimodules. Furthermore $R[{}_\calC\calC_\calC]$ is the category algebra $R\calC$ regarded as a bimodule over itself.
\end{proposition}

\begin{proof}
This is evident from the definitions.
\end{proof}

The technicalities of the bicategory $\Dist$ can be avoided by constructing instead the 1-category with the same objects, and with the isomorphism classes of distributors (or bisets) as morphisms. We are about to define a quotient category of this, that we call the biset category, by requiring that disjoint unions of bisets are identifed as sums. The reason for imposing this relation is so that the Mackey formula for biset functors defined on groups is satisfied.

\begin{definition}
Let $\calC$ and $\calD$ be finite categories. We define $A(\calD,\calC)$ to be the Grothendieck group of finite $(\calD,\calC)$-bisets with respect to disjoint union $\sqcup$. This means that $A(\calD,\calC)$ is the quotient of the free abelian group with the $(\calD,\calC)$-bisets as basis, by the subgroup generated by elements $\Omega - (\Omega_1 + \Omega_2)$ whenever $\Omega\cong \Omega_1\sqcup \Omega_2$ as bisets. Thus $A(\calD,\calC)$ is the free abelian group with symbols in bijection with the isomorphism types of $\sqcup$-indecomposable $(\calD,\calC)$-bisets as basis. If $\Omega$ is a $(\calD,\calC)$-biset we use the same symbol $\Omega$ for the coset in $A(\calD,\calC)$ that it represents. If $R$ is a commutative ring we put
$$
A_R(\calD,\calC):=R\otimes_\ZZ A(\calD,\calC).
$$
Given a further finite category $\calE$ we define a mapping
$$
A_R(\calE,\calD)\times A_R(\calD,\calE)\to{A_R(\calE,\calC)}
$$
given on bisets by the composition $\circ$. This composition sends disjoint unions to disjoint unions, so extends to an $R$-bilinear map.
\end{definition}

 In the context of groups and with the restriction that bisets be free on one side, the ring $A(\calC,\calC)$ is the \textit{double Burnside ring}, introduced in \cite{AGM}. 
Our next definition (ignoring the free condition on bisets) is an analogue for categories of the \textit{Burnside category} for groups in the terminology of \cite{AGM}, see also \cite{Bou1} and \cite{Weba}.

\begin{definition} Letting $R$ be a commutative ring, the \textit{biset category} $\BB_R$ is the category whose objects are finite categories, and where $\Hom_{\BB_R}(\calC,\calD)= A_R(\calD,\calC)$. The reversal of the symbols $\calC,\calD$ arises because of the convention that we apply morphisms from the left. Composition is determined by the composition of bisets, extended by $R$-linearity.
A \textit{biset functor} over $R$ is an $R$-linear functor $\BB_R\to R\hbox{-mod}$. Biset functors are the objects in a category $\calF_R$ in which the morphisms are natural transformations.  If the commutative ring $R$ is clear from context, we will omit it from the notation, writing simply $\BB$ instead of $\BB_R$ and $\calF$ instead of $\calF_R$.
\end{definition}

As with all functor categories to the category of $R$-modules, the category of biset functors is abelian.
We now present examples of parts of the biset category, and of biset functors. Some of them are the same as examples familiar from the theory of biset functors on groups, and the difference is they are now defined on all finite categories. 

\begin{example}
\label{discrete-category-bisets}
For each natural number $n$, consider the discrete category $[n]$ whose objects are the elements of the set $\underline n=\{1,\ldots,n\}$, and where the only morphisms are the identity morphisms. An $([m],[n])$-biset $\Omega$ is a set for the category $[m]\times [n]^\op$, which is a discrete category with  $mn$ objects that can be written as pairs $(i,j)$ where $i\in[m]$ and $j\in [n]$. A set for this category is the specification of a list of $mn$ sets $S_{ij}$ that we can list in an array
$$
\Omega=(S_{ij})=\begin{bmatrix}
S_{11}&S_{12}&\cdots&S_{1m}\cr
\vdots&&&\vdots\cr
S_{n1}&S_{n2}&\cdots&S_{nm}\cr
\end{bmatrix}.
$$
If we have a second $([\ell],[m])$-biset $\Psi=(T_{jk})$ the composite $\Omega\Psi$ has $(i,k)$ entry
$$(S_{i1}\times T_{1k})\sqcup(S_{i2}\times T_{2k})\sqcup\cdots
\sqcup (S_{im}\times T_{mk}),
$$
there being no relation to put on the disjoint unions because there are no non-identity morphisms in the discrete categories. Identifying these bisets by the sizes of the sets that are their entries, we see that isomorphism classes of bisets biject with matrices with non-negative integer entries, with composition given by matrix multiplication. Under the operation of disjoint union the bisets are the free commutative monoid spanned by the bisets $E_{ij}$ that consist of a single point set in position $(i,j)$ and are empty elsewhere. Thus the Grothendieck group  $A_\ZZ([m],[n])\cong \Mat_{m,n}(\ZZ)$ and, more generally, $A_R([m],[n])\cong \Mat_{m,n}(R)$ for any commutative ring $R$. We summarize this in the following result.

\begin{proposition}
\label{discrete-category-proposition}
Let $R$ be a commutative ring. The full subcategory of the biset category $\BB_R$ that has as its objects the discrete categories $[m]$ is equivalent to the category whose objects are free modules of finite rank over $R$, with $R$-module homomorphisms as morphisms. Thus $\End_{\BB_R}([m])\cong \Mat_{m,m}(R)$.
\end{proposition}

We see in this example the effect of imposing the relation $\Omega\sqcup\Psi = \Omega+\Psi$ on bisets, as morphisms in the biset category. Without considering the operation of disjoint union, we would obtain a theory in which $([m],[m])$-bisets form a monoid under composition that may be identified as the monoid $\Mat_{m,m}(\NN)$ under matrix multiplication. Representations of this monoid are linear representations of the monoid algebra  $R\Mat_{m,m}(\NN)$ and these include more than the linear representations of $\Mat_{m,m}(R)$ that appear with biset functors.
\end{example}

The last example suggests a useful notation for storing information about a biset.

\begin{definition}
Let $\Omega$ be a $(\calC,\calD)$-biset. Place the objects of $\calC$ and also the objects of $\calD$ in some total order. The \textit{size matrix} $|\Omega|$ of $\Omega$ is the matrix with rows indexed by the objects of $\calC$, columns indexed by the objects of $\calD$, and where the $(x,y)$ entry is the size $|\Omega(x,y)|$.
\end{definition}

\begin{example}
\label{size-matrix-example}
Let $\calC=\calE= \calA_2 =1\to 2$ and $\calD = \calA_3 = 1\xrightarrow{\alpha} 2\xrightarrow{\beta} 3$ with the objects placed in the order indicated. Consider a $(\calC,\calD)$-biset $\Omega$ and a $(\calD,\calE)$-biset  $\Psi$ with the size matrices indicated:
$$
|\Omega|=\begin{bmatrix} 1&1&0\\ 1&1&1\\ \end{bmatrix}, \quad |\Psi| = 
\begin{bmatrix}1&0\\1&0\\ 1&1\\ \end{bmatrix},\quad
|\Omega\circ\Psi| = \begin{bmatrix}1&0\\ 1&1\\ \end{bmatrix}.
$$
In this particular case, $\Omega,\Psi$ and $\Omega\circ\Psi$ are uniquely specified by their size matrices. In more general examples this would not happen. For instance, $\Omega(1,2)$ is a single point, and there is only one possibility for how it can map to the single point $\Omega(1,1)$ by the right action of $\calD$, and also to the single point $\Omega(2,2)$ by the left action of $\calC$. The size matrix for $\Omega\circ\Psi$ is computed by matrix multiplication as in Example \ref{discrete-category-bisets}, followed this time by imposing a non-trivial equivalence relation. For example, in computing $(\Omega\circ\Psi)(2,1)$ the set
$$
\bigsqcup_{y\in\Ob\calD} \Omega(2,y)\times\Psi(y,1)
$$
that appears in the definition of the product is the disjoint union of three one-point sets $\{(a,r)\},\{(b,s)\},\{(c,t)\}$ where
$$
\begin{aligned}
&\Omega(2,1)=\{a\},\quad\Omega(2,2)=\{b\},\quad\Omega(2,3)=\{c\}\\
&\Psi(1,1)=\{r\},\quad\Psi(2,1)=\{s\},\quad\Psi(3,1)=\{t\}.\\
\end{aligned}
$$
The three elements are equivalent because $(a,r) = (b\alpha,r)\sim (b,\alpha r) = (b,s)$, and similarly $(b,s)\sim (c,t)$.

The size matrix for $\Omega\circ\Psi$ is not the product $|\Omega||\Psi| = \begin{bmatrix} 2&0\\ 3&1\\ \end{bmatrix}$ and the discrepancy between the entries of these two matrices is accounted for by the fact that we imposed a non-trivial equivalence relation. Without the equivalence relation we would get straight matrix multiplication, as happens with discrete categories in Example~\ref{discrete-category-bisets}. In general the size after the equivalence relation has been imposed is, if anything, smaller than the entry in the matrix product; it must also be at least 1 if any of the sets $ \Omega(x,y)\times\Psi(y,z)$ is non-empty.
\end{example}

We have proved the following useful aid in doing calculations with bisets:

\begin{proposition}
Let $\Omega$ be a $(\calC,\calD)$-biset and $\Psi$ a $(\calD,\calE)$-biset. The entries of the size matrix $|\Omega\circ\Psi|$ are bounded by the entries of $|\Omega||\Psi|$. An entry of $|\Omega\circ\Psi|$ is zero if and only if the corresponding entry of  $|\Omega||\Psi|$ is zero.
\end{proposition}

We continue with examples of biset functors.

\begin{example}
Let $R$ be a commutative ring that is either a field or a complete local domain, so that whenever $\calC$ is a finite category, the finitely generated modules for the category algebra $R\calC$ satisfy the Krull-Schmidt theorem.
Let $K_0(R\calC,\oplus)$ be the Grothendieck group of finitely generated $R\calC$-modules with relations given only by direct sum decompositions. The assignment $F(\calC)=K_0(R\calC,\oplus)$ is a biset functor. The functorial effect on bisets is a follows. Given a $(\calC,\calD)$-biset $\Omega$ and an $R\calD$-module $M$ we get an $R\calC$-module $R\Omega\otimes_{R\calD} M$. This defines a homomorphism $K_0(R\calD,\oplus)\to K_0(R\calC,\oplus)$ that is functorial on $\BB$ because tensor product preserves direct sums.
\end{example}

\begin{example}
\label{Burnside-functor-example}
For each finite category $\calC$ the representable functor $\Hom_\BB(\calC,-)$ is a biset functor. When $\calC=\one$ and the underlying commutative ring $R$ is $\ZZ$, this is the Burnside ring functor $B$ defined in \ref{burnside-ring-definition}. This is because $\Hom_\BB(\one,\calD)$ is the Grothendieck group of $(\calD,\one)$-bisets with respect to $\sqcup$, and sets with an action of $\one$ from the right have no more structure than sets. By Yoneda's lemma, representable functors are projective in the category of biset functors (see \cite{Webb}) and the endomorphism ring of $\Hom_\BB(\calC,-)$ (as a biset functor) is isomorphic to $\End_\BB(\calC)^\op$. In particular, the Burnside ring functor is projective, and it is indecomposable over any connected ring $R$ because its endomorphism ring is $R$.
\end{example}

\section{Bisets obtained from functors}
\label{functors-to-bisets-section}

The following construction appears in \cite[Ex. 7.8.3]{Bor} in different notation.   

\begin{definition}
\label{basic-bisets-definition} Given categories $\calC$, $\calD$ and $\calE$, and functors $F:\calC\to\calE$ and $G:\calD\to\calE$  we obtain a $(\calC,\calD)$-biset that we denote ${}_{\calC^F}\calE_{{}^G\calD}$. On objects $x$ of $\calC$ and $y$ of $\calD$ this biset is defined by
$$
{}_{\calC^F}\calE_{{}^G\calD}(x,y):=\Hom_\calE(G(y),F(x)).
$$
The effect of this functor on morphisms of $\calC$ and $\calD$ is given by composition in $\calE$, after first applying  $F$ and $G$. Thus if $\alpha:x_1\to x_2$ in $\calC$, $\beta:y_2\to y_1$ in $\calD$ and $\phi:G(y_2)\to F(x_1)$ in $\calE$, then
${}_{\calC^F}\calE_{{}^G\calD}(\alpha,\beta)$ sends $\phi$ to $F(\alpha)\phi G(\beta)$.
If either of the functors $F$ or $G$ is the inclusion of a subcategory in a bigger category, we might omit it from the notation. This is consistent with the notation for the identity biset ${}_\calC\calC_\calC$. 
\end{definition}

\begin{example}
\label{group-operation-example}
In the context of biset functors defined on groups, when $H$ is a subgroup of a group $G$ the biset ${}_HG_G$ encodes the restriction operation on biset functors, while ${}_GG_H$ encodes transfer or corestriction, and if $Q$ is a factor group of $G$ the biset ${}_GQ_Q$ encodes inflation and ${}_QQ_G$ encodes deflation.
\end{example}

The bisets constructed from functors in this way are the basic connection between the category of finite categories and the biset category. We will see that some of their properties familiar from the situation of groups still hold for categories, and some need to be modified.
The following result, implied by \cite[Prop. 7.8.5]{Bor}, says that these constructions of bisets are functorial, and determines when two functors between categories yield isomorphic bisets. We use the notation $\Cat$ to denote the category whose objects are finite categories, and whose morphisms are functors.

\begin{proposition}
\label{cat-embedding}
\begin{enumerate}
\item There is a functor $\phi:\Cat\to\BB$ defined to be the identity on objects, and defined on functors $F:\calC\to\calD$ as $\phi(F)={}_\calD\calD_{{}^F\calC}:\calC\to\calD$. There is also a contravariant functor $\hat\phi:\Cat^\op\to\BB$ that is again the identity on objects, and with $\hat\phi(F)={}_{\calC^F}\calD_\calD$.
\item Under these functors $\phi$ and $\hat\phi$, two functors $F,G:\calC\to\calD$ are sent to the same morphism in $\BB$ if and only if $F$ and $G$ are naturally isomorphic.
\end{enumerate}
\end{proposition}

\begin{proof}
We describe the case of $\phi$, the situation with $\hat\phi$ being similar. In \cite[Prop. 7.8.5]{Bor} it is described that there is a is a pseudo-functor from  $\Cat$ to the distributor bicategory $\Dist$ (see Proposition~\ref{Dist-definition}). The pseudo-functor is the identity on objects, and on a functor $F:\calC\to\calD$ we get a biset (or distributor) ${}_\calD\calD_{{}^F\calC}$. This implies that if $H:\calD\to\calE$ is another functor then ${}_\calE\calE_{{}^H\calD}\circ {}_\calD\calD_{{}^F\calC}\cong {}_\calE\calE_{{}^{HF}\calC}$ as bisets.  Because bisets are taken up to isomorphism in $\BB$, the definition of $\phi:\Cat\to\BB$ preserves composition of morphisms, and we obtain a functor. This proves (1).

It is also shown in \cite[Prop. 7.8.5]{Bor} that the sets of natural transformations are preserved under the pseudo-functor: if $F, G:\calC\to\calD$ are both functors then $\mathrm{Nat}_\Cat(F,G)\leftrightarrow \mathrm{Nat}_{\Dist}( {}_\calD\calD_{{}^G\calC}, {}_\calD\calD_{{}^F\calC})$. Thus a natural transformation between the bisets is an isomorphism  if and only if the functors are naturally isomorphic,  and this is equivalent to the bisets being equal in $\BB$. This proves (2).
\end{proof}

We continue with an application of $\phi$ and $\hat\phi$. In the theory of  bisets defined when the categories are groups it is significant that the outer automorphism group $\Out(G)$ embeds in $\End_\BB(G)$,  this being one of the ingredients in the parametrization of simple biset functors by pairs $(G,V)$ where $V$ is a simple representation of $\Out G$ (see  \cite{Bou1}). There is an analogue for categories of this embedding, although we will see later that it fails to have the same significance in general that it has for groups. We describe this now. The first issue is to decide what the outer automorphism group of a category should be.

\begin{definition}
\label{out-definition}
Given a category $\calC$, we write $\Out(\calC)$ for the group of natural isomorphism classes of self-equivalences  $F:\calC\to\calC$.
\end{definition}

\begin{example}
An easy exercise shows that, when $\calC$ is a group, $\Out(\calC)$ is the usual outer automorphism group of $\calC$. Note that, for an arbitrary category $\calC$, the set of self-equivalences  $F:\calC\to\calC$ is not a group under composition.
\end{example}

In the next result we write $\Out_\Cat(\calC)$ and $ \Aut_\BB(\calC)$ to distinguish when $\calC$ is regarded as an object of $\Cat$, or of $\BB$.

\begin{corollary}
\label{out-embeds}
Let $\calC$ and $\calD$ be a finite categories.
\begin{enumerate}
\item If $\calC$ and $\calD$ are equivalent categories then they are isomorphic in the biset category.
\item The monoid homomorphism $\End_\Cat(\calC)\to \End_\BB(\calC)$ determined by the functor $\phi$ induces an injective group homomorphism
$$
\Out_\Cat(\calC)\to \Aut_\BB(\calC).
$$
Similarly, $\hat\phi$ induces an injective group homomorphism
$$
\Out_\Cat(\calC)^\op\to \Aut_\BB(\calC).
$$
\item The injection $\Out_\Cat(\calC)\to \Aut_\BB(\calC)$ realizes the group algebra $R\Out_\Cat(\calC)$ as a subalgebra of $\End_\BB(\calC)=A_R(\calC,\calC)$.
\item  For every biset functor $F$, the evaluation $F(\calC)$ has the structure of an $R\Out\calC$-module.
\end{enumerate}
\end{corollary}

\begin{proof}
(1) If $\calC$ and $\calD$ are equivalent, there are functors $F:\calC\to \calD$ and $G:\calD\to \calC$ so that $FG\simeq 1_\calD$ and $GF\simeq 1_\calC$. Now the bisets ${}_\calD\calD_{{}^F\calC}$ and ${}_\calC\calC_{{}^G\calD}$ are inverse isomorphisms in $\BB$ between $\calC$ and $\calD$ by Proposition~\ref{cat-embedding}. This is because $$
{}_\calD\calD_{{}^F\calC}\circ {}_\calC\calC_{{}^G\calD}\cong {}_\calD\calD_{{}^{FG}\calD}\cong {}_\calD\calD_{\calD}
$$
and similarly 
$$
{}_\calC\calC_{{}^G\calD}\circ {}_\calD\calD_{{}^F\calC}\cong {}_\calC\calC_{{}^{GF}\calC}\cong {}_\calC\calC_{\calC}.
$$

To prove (2), we observe from Proposition~\ref{cat-embedding} part (2) that $\phi$ is well-defined on natural isomorphism classes of functors, so induces a group homomorphism as claimed. Under $\phi$, a self-equivalence of $\calC$ is sent to the identity automorphism of $\calC$ in $\BB$ if and only if it is naturally equivalent to the identity functor on $\calC$, and this shows that the induced map on $\Out_\Cat(\calC)$ is injective.

(3) Each element of $\Out_\Cat(\calC)$ is mapped to a distinct isomorphism class of bisets by Proposition~\ref{cat-embedding}, and so the span of these elements in $A_R(\calC,\calC)$ is a copy of $R\Out_\Cat(\calC)$.

(4) is immediate from (3).
\end{proof}

We conclude this section with some results that show behavior not predicted by the part of the biset category with groups as its objects. We have already seen in Theorem~\ref{idempotent-completion-Burnside-theorem} that categories with equivalent idempotent completions have isomorphic Burnside rings, and the following fills out this picture.

\begin{theorem}
\label{idempotent-completion-isomorphism-theorem}
Suppose that $\calC$ and $\calD$ are categories whose idempotent completions are equivalent. Then $\calC$ and $\calD$ are isomorphic in the biset category $\BB$.
\end{theorem}

\begin{proof}
It is part of Theorem 7.9.4 of \cite{Bor} that if the idempotent completions of $\calC$ and $\calD$ are equivalent then  $\calC$ and $\calD$ are equivalent in the bicategory of small categories with distributors as morphisms. The biset category has natural isomorphism classes of distributors as morphisms with the relation imposed that a disjoint union of bisets is equal to the sum of the terms in the disjoint union. This means that an equivalence in the distributor category implies an isomorphism in the biset category.
\end{proof}

When we come to the parametrization of simple biset functors, Theorem~\ref{idempotent-completion-isomorphism-theorem} is an obstacle to extending the approach used for biset functors on groups. An ingredient in the parametrization of a simple biset functor $S$ is to take a group $H$ of least size such that $S(H)$ is non-zero. Because an isomorphism type of categories in $\BB$ may have representatives of different sizes, an appropriate choice of category on which $S$ is non-zero is less clear than it is with groups.

Our final observation in this section has to do with factorization of bisets.
It is shown in \cite{Bou1} that every indecomposable $(G_1,G_2)$-biset, where $G_1,G_2$ are groups, can be written as a composite of four bisets of the form ${}_HJ_K$ obtained from group homomorphisms $H\to J\gets K$ where $H,J,K$ are all no bigger than the largest of $G_1$ and $G_2$, and this turns out to be important in the development of the theory for groups in many ways.  The analogous statement is not true for bisets for categories in general, as can be seen by taking a situation where there are infinitely many $(\calC,\calD)$-bisets (as in Example~\ref{size-matrix-example}), but only finitely many functors between  categories smaller than $\calC$ or $\calD$. We do, however, show that every $(\calC,\calD)$-biset ${}_\calC\Omega_\calD$ can be written as
$$
{}_\calC\Omega_\calD={}_{\calC^F}\calE_{{}^G\calD}={}_{\calC}\calE_\calE\circ{}_\calE\calE_{\calD}
$$
where $\calE$ is a larger category that has $\calC$ and $\calD$ as full subcategories and $F:\calC\to\calE$ and $G:\calD\to\calE$ are the inclusion functors. This is a decomposition that is not available with groups, and it is possible because categories may have more than one object. It implies that bisets of the form ${}_{\calC}\calE_\calE$ and ${}_\calE\calE_{\calD}$, with $\calC$ and $\calD$ full subcategories of a category $\calE$, generate the biset category.

Given a $(\calC,\calD)$-biset ${}_\calC\Omega_\calD$ we will construct the category $\calE$ that is variously known as the \textit{cograph}, \textit{bridge} or \textit{collage} of the biset. We will denote it $\Cograph(\Omega)$. See  \cite{Pec} for an account of this theory. The objects of $\Cograph(\Omega)$ are $\Ob\calC\sqcup\Ob\calD$ and the morphisms are defined as follows:
$$
\Hom_{\Cograph(\Omega)}(x,y)=
\begin{cases}
\Hom_{\calC}(x,y)&\hbox{if }x,y\in\calC,\cr
\Hom_{\calD}(x,y)&\hbox{if }x,y\in\calD,\cr
\Omega(y,x)&\hbox{if }x\in\calD\hbox{ and }y\in\calC,\cr
\emptyset&\hbox{if }x\in\calC\hbox{ and }y\in\calD.\cr
\end{cases}
$$
The composition of two morphisms that are both in $\calC$, or both in $\calD$, is the same as it was before in $\calC$ or $\calD$. If $\alpha:x_1\to x_2$ in $\calD$ and $\phi\in\Omega(y,x_2)$ then $\phi\alpha:=\Omega(1_y,\alpha)(\phi)$. If $\beta:y_1\to y_2$ in $\calC$ and $\phi\in\Omega(y_1,x)$ then $\beta\phi:=\Omega(\beta,1_x)(\phi)$.

Part of the following appears as Exercise 7.10.11 in \cite{Bor}.

\begin{proposition}
\begin{enumerate}
    \item Let $\Omega$ be a $(\calC,\calD)$-biset and let $\calE = \Cograph(\Omega)$. Then $\calE$ has $\calC$ and $\calD$ as full subcategories, and $\Omega$ can be written $\Omega ={}_{\calC^F}\calE_{{}^G\calD}$ as a $(\calC,\calD)$-biset, where $F$ and $G$ are the inclusion functors.
    \item Let $\calE$ be a category with full subcategories $\calC, \calD$ whose objects partition the objects of $\calE$ as $\Ob\calE=\Ob\calC\sqcup\Ob\calD$. Suppose that there are no morphisms $x\to y$ with $x\in\calC$ and $y\in\calD$. Then $\calE=\Cograph({}_\calC\calE_\calD)$ as categories.
\end{enumerate}
\end{proposition}

\begin{proof}
(1) From the definitions we have 
$$
\Omega(y,x)=\Hom_{\Cograph(\Omega)}(x,y)
$$
whenever $x\in\calD$ and $y\in\calC$.

(2) From the definitions again, $\calE$ and $\Cat({}_\calC\calE_\calD)$ have the same objects, and the same morphisms.
\end{proof}

\section{Representable sets and bisets}
In the theory of biset functors defined on groups it may be necessary to restrict the kinds of bisets that appear, because some important biset functors are not defined on all bisets. A natural condition when the categories are groups is to consider only bisets where the stabilizers on each side lie in specified classes of subgroups. In \cite{Weba} bisets free on restriction to both sides, and also on restriction to only one side, were considered, because group cohomology is not defined on arbitrary bisets. By developing the theory using such bisets an application to the computation of group cohomology was made in that paper. More general stabilizer restrictions on bisets were considered in \cite{Bou1}.

In the case of bisets for categories it no longer makes sense to talk about stabilizers of elements, and so we replace a condition on stabilizers by a condition that does make sense. 

\begin{definition}
We will say that a $\calC$-set is \textit{representable} if it is isomorphic to a disjoint union of sets of the form $\Hom_\calC(x,-)$ for various objects $x$ of $\calC$. Furthermore, a $\calC$-set $\Omega$ is \textit{representably generated by $\omega\in \Omega(x)$} if and only if the map $\Hom_\calC(x,-)\to\Omega$, given at $y$ on a morphism $\alpha:x\to y$ by $\alpha\mapsto (\Omega(\alpha))(\omega)$, is an isomorphism. 
\end{definition}

Note in this definition that we are making a departure from standard terminology, in that a representable functor is one isomorphic to $\Hom(x,-)$, whereas we are defining a representable $\calC$-set to be one that is a \textit{disjoint union} of representable functors. The convenience of this abuse of terminology justifies its use. If a $\calC$-set is representably generated by some $\omega$ then it is  representable in the standard sense.

\begin{example}
When the category $\calC$ is a group, a $\calC$-set is representable if and only if the action of the group on the set is free, meaning that the stabilizer of each element is 1. This is because the group has only one object, and the functor represented by that object is the regular representation.
\end{example}

\begin{proposition}
\label{representable-generation}
\begin{enumerate}
\item The representable functor $\Hom_\calC(x,-)$ is representably generated at $z$ by $f:x\to z$ if and only if $f$ is an isomorphism in $\calC$. In particular, $\Hom_\calC(x,-)$ is representably generated at $x$ by $1_x$.
\item The representable functors $\Hom_\calC(x,-)$ and $\Hom_\calC(y,-)$ are isomorphic if and only if $x\cong y$ in $\calC$.
\item The isomorphism classes of indecomposable representable $\calC$-sets are precisely the $\calC$-sets $\Hom_\calC(x,-)$, where $x\in\calC$ is taken up to isomorphism.
\end{enumerate}
\end{proposition}

\begin{proof}
These are well known properties of the Yoneda embedding. 

(1) The condition that $\Hom_\calC(x,-)$ be representably generated at $z\in\Ob\calC$ by $f:x\to z$ is that the map $\Hom_\calC(z,-)\to \Hom_\calC(x,-)$ given at $y$ on a morphism $\alpha:z\to y$ by $\alpha\mapsto \alpha f$ is an isomorphism of functors. If $f$ is an isomorphism in $\calC$ then $\beta\mapsto \beta f^{-1}$ is inverse to the previous map, so the condition on the right implies the condition on the left. Conversely, if  $\alpha\mapsto \alpha f$ is an isomorphism of functors the fact that the Yoneda embedding is faithful means that $f$ is an isomorphism.

(2) It is immediate that the condition on the right implies the condition on the left. Conversely, inverse isomorphisms
$$
\Hom_\calC(x,-)\to \Hom_\calC(y,-)\to \Hom_\calC(x,-)
$$
are determined by mappings $x\xleftarrow{\phi} y\xleftarrow{\theta} x$ that are the images of $1_y$ and $1_x$, that compose in both directions to give the identities on $x$ and $y$, so isomorphism of the functors implies isomorphism of $x$ and $y$. 

For (3), each $\Hom_\calC(x,-)$ is indecomposable as a $\calC$-set by Lemma~\ref{cyclic-implies-indecomposable}, because it is generated by $1_x$, and the result now follows from the definition and (2).
\end{proof}

\begin{example}
\label{monoid-idempotent-completion}
It is possible that a representable $\calC$-set $\Hom_\calC(y,-)$ can be generated at some other object $x\in\calC$ with $x\not\cong y$. It follows from Proposition~\ref{representable-generation} part (2) that the $\calC$-set  cannot be representably generated at $x$. Consider the category
\begin{center}
\begin{tikzpicture}[xscale=1,yscale=1]
\node at (-1.8,0) {$\calC =$};
\node (A) at (0, 0) {$\bullet$};
\node (B) at (1.2,0) {$\bullet$};
\node [ left] at  (A)   {$x$};
\node [ right] at  (B)   {$y$};
\draw [->] (0.2,0.1)--(1,0.1) node[midway,above]{$u$};
\draw [->] (1,-0.1)--(0.2,-0.1) node[midway,below]{$v$};
\draw [->] (-0.1,0.2) to [out=110,in=0] (-0.4,0.4) to  [out=180,in=90] (-0.8,0) to [out=270,in=180] (-0.4,-0.4) to [out=0,in=260] (-0.1,-0.2);
\node[left] at (-0.8,0) {$\alpha$};
\end{tikzpicture}
\end{center}
with compositions $uv=1_y$, $vu=\alpha$, so that $\alpha^2=\alpha\ne 1_x$. The representable functor $\Hom_\calC(y,-)$ is representably generated at $y$, and is also generated at $x$ by $v\in\Hom_\calC(y,x)$, but not representably so. The mapping $\Hom_\calC(x,-)\to\Hom_\calC(y,-)$, given by $f\mapsto fv$, is onto but not one-to-one, because $\alpha v=v=1_x v$.  
\end{example}

\begin{proposition}
\label{representably-generated-criterion}
Let $\Omega:\calC\to\Set$ be a $\calC$-set, let $\omega\in\Omega(x)$ for some object $x$, and let $\mathrm{ev}_x:\calC\hbox{-Set}\to\Set$ be the functor that is evaluation at $x$. The following are equivalent.
\begin{enumerate}
\item $\Omega$ is representably generated at $x$ by $\omega$.
\item for all objects $y\in\calC$, for all $u\in \Omega(y)$, there exists a unique morphism $\alpha:x\to y$ so that $u=\Omega(\alpha)(\omega)$.
\item As functors $\calC\hbox{-Set}\to\hbox{Set}$, there is a natural isomorphism
$$
\Hom_\CSet(\Omega,-)\cong\mathrm{ev}_x
$$
whose value  $\Hom_\CSet(\Omega,\Omega)\cong \Omega(x)$ at $\Omega$ sends  $1_\Omega\mapsto \omega$.
\end{enumerate}
\end{proposition}

\begin{proof}
(1) implies (2) because condition (2) holds for the representable functor $\Hom_\calC(x,-)$ with regard to the element $1_x$ that representably generates it. 

(2) implies (1). Consider the map $\Hom_\calC(x,-)\to\Omega$ given by $\alpha\mapsto \Omega(\alpha)(\omega)$. The hypothesis is that this map is bijective, so it is an isomorphism. From this, (1) holds.

(1) implies (3) because if $\Omega$ is representably generated at $x$ then, by definition, there is a natural isomorphism $\Hom_\calC(x,-)\to\Omega$ given at $y$ on a morphism $\alpha:x\to y$ by $\alpha\mapsto (\Omega(\alpha))(\omega)$.
Yoneda's lemma says that, for all $\calC$-sets $P$, we have an isomorphism
$$
\Hom_\CSet(\Hom_\calC(x,-),P)\cong P(x),
$$
so that $\Hom_\CSet(\Hom_\calC(x,-),-)\cong\mathrm{ev}_x$. Combining this with the previous natural isomorphism we get  $\Hom_\CSet(\Omega,-)\cong\mathrm{ev}_x$. 

Consider the following two isomorphisms, the first of which is pre-composition with the earlier isomorphism $\Hom_\calC(x,-)\to\Omega$ and the second is the isomorphism in Yoneda's lemma:
$$
\Hom_\CSet(\Omega,\Omega)\to\Hom_\CSet(\Hom_\calC(x,-),\Omega) \to \Omega(x).
$$
The first of these maps $1_\Omega$ to the natural map $\alpha\mapsto (\Omega(\alpha))(\omega)$, and the second sends this to its value on taking $\alpha=1_x$.  Thus $1_\Omega$ is sent to $\omega$.

(3) implies (1). The isomorphism given in (3) is a bijection of natural transformations $\Hom_\CSet(\Omega,-)\cong\Hom_\CSet(\Hom_\calC(x,-),-)$ and the fact that the Yoneda embedding is faithful means that $\Omega\cong \Hom_\calC(x,-)$ as $\calC$-sets. The condition that $1_\Omega\mapsto \omega$ in the isomorphism $\Hom_\calC(x,-)\to \Omega$ at $x$ implies that $\Omega$ is representably generated by $\omega\in\Omega(x)$, because if $\alpha: x\to y$ in $\calC$ then $\alpha=\alpha1_x\mapsto \Omega(\alpha)\Omega(1_x) = \Omega(\alpha) \omega$ by naturality of the isomorhism, and this verifies the definition of being representably generated at $x$ by $\omega$.
\end{proof}

\begin{corollary}
If the $\calC$-set $\Omega$ is representably generated by both $\omega\in\Omega(x)$ and $\psi\in\Omega(y)$ then there is an isomorphism $\alpha:x\to y$ in $\calC$ with $\Omega(\alpha)(\omega)=\psi$.
\end{corollary}

\begin{proof}
Suppose that $\Omega$ is representably generated by both $\omega$ and $\psi$. By Proposition~\ref{representably-generated-criterion} there are unique morphisms $\alpha:x\to y$ and $\beta:y\to x$ for which $\Omega(\alpha)(\omega)=\psi$ and $\Omega(\beta)(\psi)=\omega$. Now $\beta\alpha$ is the unique morphism with $\Omega(\beta\alpha)(\omega)=\omega$, so that $\beta\alpha=1$ because 1 also has this property. Similarly $\alpha\beta=1$ so that $\alpha$ is an isomorphism.
\end{proof}.

We now consider representable bisets. Given a $(\calC,\calD)$-biset ${}_\calC\Omega_\calD$ we can regard it as a left $\calC$-set ${}_\calC\Omega$ and also as a right $\calC$-set $\Omega_\calD$ by restricting the actions: as a left $\calC$-set its value on an object $x$ of $\calC$ is defined to be $\bigsqcup_{y\in \calD}\Omega(x,y)$ and as a right $\calD$ set the value on an object $y$ of $\calD$ is defined to be $\bigsqcup_{x\in \calC}\Omega(x,y)$. 

\begin{definition}
We say that a $(\calC,\calD)$-biset $\Omega$ is \textit{representable on the left} if the left $\calC$-set ${}_\calC\Omega$ is representable, and \textit{representable on the right} if the right $\calD$-set $\Omega_\calD$ is representable (that is, representable as a $\calD^\op$-set).  It is \textit{birepresentable} if it is both representable on the left and representable on the right.
\end{definition}

We are about to show that the composite of bisets that are representable on one side is also representable on that side, and to do this it is convenient to consider the intermediate set constructed from two bisets before the relation is imposed that defines the composition of the bisets. Given bisets ${}_\calC\Omega_\calD$ and ${}_\calD\Psi_\calE$ we define
$$ \Omega \odot\Psi(x,z): =\bigsqcup_{y\in\calD} \Omega(x,y)\times\Psi(y,z).
$$

Thus $\Omega\circ\Psi(x,z)= \Omega \odot\Psi(x,z)/\sim$ where $\sim$ is the relation generated by $(u\delta,v)\sim(u,\delta v)$, with $\delta\in\calD$. The relation $\sim$ will appear throughout the next lemmas. We will denote the equivalence class of $(u,v)$ under $\sim$ in $\Omega \odot\Psi(x,z)$ by $u\boxtimes v$.

As notation, we recall that if $\Omega= \bigsqcup\Omega_i$ is a disjoint union of indecomposable $\calC$-subsets and $u\in \bigsqcup_{x\in\Ob\calC}\Omega(x)$ we define the \textit{orbit} of $u$ to be the unique indecomposable summand $\Omega_i$ of $\Omega$ for which $u\in \bigsqcup_{x\in\Ob\calC}\Omega_i(x)$. We will denote this orbit by $\calC u$, with the understanding that $\calC u$ might not be generated as a $\calC$-set by $u$, or indeed by any single element. Note that $\calC u = \calC \alpha u$ for all morphisms $\alpha$ in $\calC$ for which $\alpha u$ is defined.

The next two results are stated for orbits that are representable on the left. Similar statements hold for orbits that are representable on the right.

\begin{lemma}
\label{equivalence-class-lemma}
Let ${}_\calC\Omega_\calD$ and ${}_\calD\Psi_\calE$ be bisets, suppose that $\hat u\in\Omega(x,y)$ representably generates the $\calC$-orbit $\calC\hat u$ and that $\hat v\in\Psi(y,z)$ representably generates $\calD\hat v$. Then, if $\alpha:x\to x_1$ in $\calC$, the elements of $\Omega\odot\Psi(x_1,z)$ equivalent to $(\alpha\hat u,\hat v)$ are precisely the set
$$
\begin{aligned}
E_\alpha^{(\hat u,\hat v)}:&=\{(u,\delta\hat v)\bigm| u\delta=\alpha\hat u,\;u\in\Omega(x_1,y_1),\;y_1\in\calD,\;\delta:y\to y_1\}\\&\subseteq\bigsqcup_{w\in\calD} \Omega(x_1,w)\times\Psi(w,z).\\
\end{aligned}
$$
Thus $\alpha\hat u\boxtimes \hat v=E_\alpha^{(\hat u,\hat v)}$.
\end{lemma}

\begin{proof}
We use the notation that appears in the definition of $E_\alpha^{(\hat u,\hat v)}$.
Because $(u,\delta\hat v)\sim (u\delta,\hat v)=(\alpha\hat u,\hat v)$, all elements of $E_\alpha^{(\hat u,\hat v)}$ are equivalent to $(\alpha\hat u,\hat v)$. On the other hand, if $(u,\delta\hat v)\in E_\alpha$ is equivalent to some further element, then it arises through relations of two kinds:
\begin{itemize}
\item $(u,\delta\hat v)=(u_1\delta_1,\delta\hat v)\sim (u_1,\delta_1\delta\hat v)$. In this case $u_1\delta_1\delta=\alpha\hat u$, so this $(u_1,\delta_1\delta\hat v)$ lies in $E_\alpha^{(\hat u,\hat v)}$ already (taking $u$ to be $u_1$).
\item $(u,\delta\hat v)=(u,\gamma v_1)\sim (u\gamma,v_1)$ for some morphism $\gamma\in\calD$ and some $v_1$. Now $\calD\hat v=\calD v_1$ so $v_1=\epsilon\hat v$ for some unique morphism $\epsilon\in\calD$, because $\hat v$ generates this orbit representably. Then $\gamma v_1=\gamma\epsilon\hat v=\delta\hat v$ so $\gamma\epsilon=\delta$ for the same reason.
Now $(u\gamma,v_1)=(u\gamma,\epsilon\hat v)$  where $u\gamma\epsilon=u\delta=\alpha\hat u$. Thus $(u\gamma,v_1)\in E_\alpha^{(\hat u,\hat v)}$ already, so $E_\alpha^{(\hat u,\hat v)}$ is an equivalence class under $\sim$.
\end{itemize}
This completes the proof.
\end{proof}

\begin{corollary}
\label{hats-generate}
As in Lemma~\ref{equivalence-class-lemma}, let $\hat u\in\Omega(x,y)$ and $\hat v\in\Psi(y,z)$ representably generate their orbits under $\calC$ and $\calD$, respectively. Then in $\Omega\circ\Psi$ the equivalence class $\hat u\boxtimes \hat v$ representably generates its $\calC$-orbit.
\end{corollary}

\begin{proof}
From the definition of the composition of two bisets, the left action of $\calC$ on $\Omega$ determines a left action of $\calC$ on the equivalence classes in $\Omega\odot\Psi$, and it is given by $\alpha (u\boxtimes v) = \alpha u\boxtimes v$ when $\alpha$ is a morphism in $\calC$ (and when these expressions are defined). The issue is that if $(u,v)\sim(u',v')$ then  $(\alpha u,v)\sim(\alpha u',v')$, which is apparent.

We first show that any element in the $\calC$-orbit of  $\hat u\boxtimes \hat v$ can be written  $\alpha\hat u\boxtimes \hat v$ for some morphism $\alpha$ in $\calC$. The set of such elements is closed under applying morphisms $\beta$ of $\calC$, meaning that $\beta(\alpha \hat u\boxtimes \hat v) =( \beta\alpha) \hat u\boxtimes \hat v$ remains in this set. We must also show that if $\alpha \hat u\boxtimes \hat v=\beta r\boxtimes s$ for some morphism $\beta$ in $\calC$, so that $r\boxtimes s$ lies in the same $\calC$-orbit as $u\boxtimes v$, then $r\boxtimes s$ can also be written in the form $\gamma\hat u\boxtimes \hat v$ for some $\gamma$. This equation says that $(\beta r,s)\in E_\alpha^{(\hat u,\hat v)}$, so $(\beta r,s) = (u,\delta\hat v)$ for some $u,\delta$ with $u\delta = \alpha\hat u$, by  Lemma~\ref{equivalence-class-lemma} . Thus $u=\beta r$, $s=\delta\hat v$, so $(r,s) = (r,\delta\hat v)\sim (r\delta,\hat v)$ and $u\delta = \beta r\delta =\alpha\hat u$ so $r\delta$ lies in the $\calC$-orbit of $\hat u$ in $\Omega$. Thus $r\delta = \gamma\hat u$ for some $\gamma$ in $\calC$ because $\hat u$ generates its $\calC$-orbit representably. It follows that $r\boxtimes s = r\boxtimes \delta \hat v = r\delta\boxtimes \hat v = \gamma\hat u\boxtimes\hat v$, as required.

To show that $\hat u\boxtimes\hat v$ generates its orbit representably it only remains to show that if $\alpha_1\hat u\boxtimes\hat v= \alpha_2\hat u\boxtimes\hat v$ then $\alpha_1=\alpha_2$, by Proposition~\ref{representably-generated-criterion}(2). For this we again use the description of these equivalence classes in Lemma~\ref{equivalence-class-lemma} as $\alpha\hat u\boxtimes\hat v= E_\alpha^{(\hat u,\hat v)}$.
If $E_{\alpha_1}^{(\hat u,\hat v)}=E_{\alpha_2}^{(\hat u,\hat v)}$ then  $(\alpha_1\hat u,\hat v)\in E_{\alpha_2}^{(\hat u,\hat v)}$ so $(\alpha_1\hat u,\hat v)=(u,\delta\hat v)$ where $u\delta=\alpha_2\hat u$, for some morphism $\delta$ in $\calD$. Since $\hat v$ generates its $\calD$-orbit representably, $\delta=1$, so $u=\alpha_2\hat u=\alpha_1\hat u$. Since $\hat u$ generates its $\calC$-orbit representably, $\alpha_1=\alpha_2$, as required.
\end{proof}

The following result, that composites of representable bisets are representable, is one of the main goals of this section.

\begin{theorem}
Let $\calC,\calD$ and $\calE$ be categories. 
If the bisets ${}_\calC\Omega_\calD$ and ${}_\calD\Psi_\xi$ are both representable on the left, then so is $\Omega\circ\Psi$.
If ${}_\calC\Omega_\calD$ and ${}_\calD\Psi_\xi$ are both representable on the right, then so is $\Omega\circ\Psi$.
\end{theorem}

\begin{proof}
We only prove the statement about left representability. The statement for right representability is proved in the same way as for the left by reversing the roles of left and right. 

We have seen in Corollary~\ref{hats-generate}  that the $\calC$-orbits generated by elements $\hat u\boxtimes\hat v$ are representable (if $\hat u$ and $\hat v$ representably generate their orbits), but we have not seen that all $\calC$-orbits have this form.
Let $(u,v)\in\Omega\odot\Psi$. To show that $u\boxtimes v$ lies in a representable $\calC$-orbit, we find its representable generator. Now $v=\delta\hat v$ for some $\delta\in \calD$, for some representable generator $\hat v$ of the $\calD$-orbit of $v$, and $(u,v)=(u,\delta\hat v)\sim (u\delta,\hat v)$. Next, $u\delta=\alpha\hat u$ for some $\alpha\in\calC$ and some $\hat u$ that representably generates the $\calC$-orbit of $u\delta$. Now $u\boxtimes v$ is in the $\calC$-orbit of $\hat u\boxtimes \hat v$, which is representably generated, by Corollary~\ref{hats-generate}.
\end{proof}

The next result examines when bisets obtained from functors are representable, and has the implication that the identity bisets are birepresentable. We continue with the notation of Definition~\ref{basic-bisets-definition}.

\begin{proposition}
\label{regular-is-representable}
Let $F:\calD\to\calC$ be a functor. Then the $(\calC,\calD)$-biset ${}_\calC\calC_{{}^F\calD}$ is representable on the left and the $(\calD,\calC)$-biset ${}_{\calD^F}\calC_\calC$ is representable on the right. Furthermore, if $F$ has a right adjoint then  ${}_\calC\calC_{{}^F\calD}$ is also representable on the right, and if  $F$ has a left adjoint then  ${}_{\calD^F}\calC_\calC$  is also representable on the left.
\end{proposition}

\begin{proof} We prove only the assertions about ${}_\calC\calC_{{}^F\calD}$, the case of ${}_{\calD^F}\calC_\calC$ being similar
As a left $\calC$-set ${}_\calC\calC_{{}^F\calD}=\bigsqcup_{y\in\calD}\Hom_\calC(Fy,-)$, which is a disjoint union of representable sets, representable by the objects $Fy$, without any condition on $F$. 

If $F$ has a right adjoint $G:\calC\to\calD$ then, as a right $\calD$-set,
$$
{}_\calC\calC_{{}^F\calD}=\bigsqcup_{x\in\calC}\Hom_\calC(F(-),x).
$$
Each right $\calD$-set $\Hom_\calC(F(-),x)\cong\Hom_\calD(-,G(x))$ is representable, represented by $G(x)$. 
\end{proof}

We point out that if a functor $F:\calD\to\calC$ is an equivalence then it has both a left and a right adjoint, so that the bisets ${}_\calC\calC_{{}^F\calD}$ and ${}_{\calD^F}\calC_\calC$ are both birepresentable. In particular, the identity biset ${}_\calC\calC_\calC$ is birepresentable, as also follows from the next criterion.

\begin{corollary}
A left $\calC$-set ${}_\calC\Omega$ is representable if and only if ${}_\calC\Omega$ is a summand (with respect to $\sqcup$) of the disjoint union of $n$ copies of ${}_\calC\calC$, for some $n$. Similarly, a right $\calD$-set $\Omega_\calD$ is representable if and only if $\Omega$ is a summand of the disjoint union of $n$ copies of  $\calD_\calD$, for some $n$. In particular, the identity biset ${}_\calC\calC_\calC$ is birepresentable.
\end{corollary}

\begin{proof}
We see from Proposition~\ref{regular-is-representable} that ${}_\calC\calC$ is representable, and so every summand of ${}_\calC\calC^n$ has representable orbits and so is representable. Conversely, every isomorphism type of indecomposable representable left $\calC$-sets appears as a summand of ${}_\calC\calC$. Thus if $n$ is the largest multiplicity of an indecomposable summand in ${}_\calC\Omega$ then ${}_\calC\Omega$ is a summand of ${}_\calC\calC^n$. The proof for right sets is similar.
\end{proof}

\begin{example}
This example shows the use of Proposition~\ref{regular-is-representable} in establishing representability, and also shows that representability may hold in more general circumstances, not predicted by  Proposition~\ref{regular-is-representable}. For future reference in Section~\ref{cohomology-section} it also provides an example of a subcategory of another category where the left and right adjoints of the restriction functor on representations do not coincide, but where we will see that we have transfer maps in homology and cohomology. While the example is specific, it is clear that there are many other similar examples with these properties.

\label{birepresentable-example}
Consider the category
\begin{center}
\begin{tikzpicture}[xscale=1,yscale=1]
\node at (-2.5,0) {$\calC =$};
\node (A) at (0, 0) {$\bullet$};
\node (B) at (1.2,0) {$\bullet$};
\node [ left] at  (A)   {$x$};
\node [ right] at  (B)   {$y$};
\draw [->] (0.2,0.1)--(1,0.1) node[midway,above]{$u$};
\draw [->] (0.2,-0.1)--(1,-0.1) node[midway,below]{$v$};
\draw [->] (-0.1,0.2) to [out=110,in=0] (-0.4,0.4) to  [out=180,in=90] (-0.8,0) to [out=270,in=180] (-0.4,-0.4) to [out=0,in=260] (-0.1,-0.2);
\node[left] at (-0.8,0) {$C_4$};
\draw [->] (1.3,0.2) to [out=80,in=180] (1.6,0.4) to  [out=0,in=90] (2,0) to [out=270,in=0] (1.6,-0.4) to [out=180,in=280] (1.3,-0.2);
\node[right] at (2,0) {$C_2$};
\end{tikzpicture}
\end{center}
with objects $x$ and $y$, where $\End_\calC(x) = C_4 = \langle g\rangle$ is a cyclic group of order 4 generated by $g$ and $\End_\calC(y) = C_2 = \langle h\rangle$ is a cyclic group of order 2 generated by $h$. There are two morphisms $u,v:x\to y$ with compositions $ug=v$, $vg=u$, $hu=v$, $hv=u$. Let $\calD$ be the full subcategory of $\calC$ whose only object is $y$.

\begin{proposition}
\label{C4C2-category-proposition}
Let $\calC$ and $\calD$ be the categories just described and let $R$ be a commutative ring. 
\begin{enumerate}
\item The inclusion functor $\calD\to\calC$ has a left adjoint, but no right adjoint.
\item The bisets ${}_\calC\calC_\calD$ and ${}_\calD\calC_\calC$ are both birepresentable.
\item The left and right adjoints of the restriction functor from $R\calC$-modules to $R\calD$-modules do not coincide.
\end{enumerate}
\end{proposition}

\begin{proof}
(1) The left adjoint of inclusion is $G:\calC\to\calD$ determined by $G(x)=G(y)=y$ and $G(g)=G(h)=h$.

If the inclusion were to have a right adjoint $H:\calC\to\calD$ then necessarily $H(x)=H(y)=y$. We would have to have a bijection $\Hom_\calD(y,H(x)) \leftrightarrow \Hom_\calC(y,x)$, but the term on the right is empty and the term on the left is not, so such $H$ does not exist.

(2) As a consequence of  Proposition~\ref{regular-is-representable}, ${}_\calD\calC_\calC$ is representable on both sides, and ${}_\calC\calC_\calD$ is representable on the left.

It is also the case that ${}_\calC\calC_\calD$ is representable on the right. 
As a $\calD$-set it is $\bigsqcup_{u\in\calC}\Hom_\calC(-,u)$. Its value at the only object $y\in\calD$ is $\Hom_\calC(y,y)=\Hom_\calD(y,y)$ so this $\calD$-set is the representable set at $y$.

(3) It suffices to prove that the left and right adjoints do not coincide when $R$ is a field: for a general ring, if the left and right adjoints were to coincide, they would on factoring out a maximal ideal. Now the left adjoint of restriction, applied to the regular representation of  $\End_\calD(y) = C_2$, returns the projective cover over $R\calC$ of the simple module with support at $y$, whereas the right adjoint returns the injective hull of the simple with support at $y$. These are non-isomorphic, by standard calculations with EI categories, see~\cite{Webb}.
\end{proof}

\end{example}

\begin{definition}
We define $A^{1,\all}(\calC,\calD)$ to be the Grothendieck group of finite $(\calC,\calD)$-bisets that are representable on the left, with respect to disjoint union $\sqcup$.  Thus $A^{1,\all}(\calC,\calD)$ is the free abelian group with basis in bijection with the isomorphism types of $\sqcup$-indecomposable $(\calC,\calD)$-bisets, representable on the left. If $R$ is a commutative ring with 1 we put
$$
A_R^{1,\all}(\calC,\calD):=R\otimes_\ZZ A^{1,\all}(\calC,\calD)
$$
but often omit $R$ from the notation when it is understood. Similarly $A_R^{\all,1}(\calC,\calD)$ and $A_R^{1,1}(\calC,\calD)$ are the same using only bisets that are representable on the right, or on both sides, in the two cases.

We define $\BB_R^{1,\all}, \BB_R^{\all,1}$ and $\BB_R^{1,1}$ to be the subcategories of $\BB_R$ with the same objects as $\BB_R$ (that is, finite categories) where
$$
\Hom_{\BB_R^{\calX,\calY}}(\calD,\calC)=A_R^{\calX,\calY}(\calC,\calD),
$$
with $\calX,\calY\in\{1,\all\}$. To be consistent, we define $\BB_R^{\all,\all}=\BB_R$ to be two different names for the same biset category. Again, we often omit $R$ from the notation. 

We use the term \textit{biset functor over $R$} to include the case of $R$-linear functors $\BB_R^{\calX,\calY}\to R\hbox{-mod}$ for any choice of $\calX,\calY\in\{1,\all\}$. The category of such biset functors may be denoted $\calF_R^{\calX,\calY}$, but usually we will omit the subscript and superscripts.
This generalizes biset functors of a kind that appear in \cite{AGM, Bou1, Weba} for groups.
\end{definition}

We summarize some of the development we have made.

\begin{theorem}
The constructions $\BB_R^{1,\all}, \BB_R^{\all,1}$ and $\BB_R^{1,1}$ are all categories. The functor $\phi:\Cat\to\BB_R$ of Proposition~\ref{cat-embedding}, that is the identity on objects and sends a functor $F:\calC\to\calD$ to the biset ${}_\calD\calD_{{}^F\calC}$, has image in the subcategory $\BB_R^{1,\all}$. The functor $\hat\phi:\Cat^\op\to\BB_R$ that is the identity on objects and sends a functor $F:\calC\to\calD$ to the biset ${}_{\calC^F}\calD_\calD$ has image in
$\BB_R^{\all,1}$. Both $\phi$ and $\hat\phi$ send functors that are equivalences to morphisms in $\BB_R^{1,1}$. Thus $\phi$ induces a group homomorphism
$
\Out_\Cat(\calC)\to \Aut_{\BB_R^{1,1}}(\calC)
$
and $\hat\phi$ induces a group homomorphism
$
\Out_\Cat(\calC)^\op\to \Aut_{\BB_R^{1,1}}(\calC).
$
\end{theorem}

Note that $ \Aut_{\BB_R^{1,1}}(\calC)\subseteq  \Aut_{\BB_R^{\calX,\calY}}(\calC)$ for all choices of $\calX$ and $\calY$.

\begin{example}
\label{Kronecker-biset-example}
Although a finite category may have infinitely many indecomposable sets, it will only have finitely many indecomposable representable sets, because these are in bijection with the isomorphism classes of objects in the category. In spite of this, there may be infinitely many indecomposable birepresentable $(\calC,\calD)$-bisets for a pair of finite categories $\calC, \calD$. 
As an example, take $\calC=\calD$  to be the \textit{Kronecker category} we considered in Example~\ref{Kronecker-example}. This has two objects $x$ and $y$, and two morphisms $\alpha,\beta:x\to y$.
 To distinguish the two copies of this category, will write $\bar\alpha, \bar\beta$ instead of $\alpha,\beta$ when these morphisms lie in the second category $\calD$. For each $n\ge 1$ consider the biset $\Omega_n:\calC\times\calD^\op\to\Set$ with 
$$
\begin{aligned}
\Omega_n((x,x)) &= \{u_1,\ldots ,u_n\},\\
\Omega_n((y,y))&= \{v_1,\ldots ,v_n\},\\
\Omega_n((y,x))&= \{w_1,\ldots ,w_{2n+1}\},\\
\Omega_n((x,y))&=\emptyset.\\
\end{aligned}
$$
We define
$$
\begin{aligned}
\Omega_n(\alpha,1_x)(u_i)&=w_{2i-1}\\
\Omega_n(\beta,1_x)(u_i)&=w_{2i}\\
\Omega_n(1_y,\bar\alpha)(v_i)&=w_{2i}\\
\Omega_n(1_y,\bar\beta)(v_i)&=w_{2i+1}\\
\end{aligned}
$$
This is depicted in the following diagram for $n=3$, where the effects of $\Omega_3(\alpha,1_x), \Omega_3(\beta,1_x)$ are written simply $\alpha,\beta$, and the effects of $\Omega_3(1_y,\bar\alpha),\Omega_3(1_y,\bar\beta)$ are written $\bar\alpha, \bar\beta$.
\begin{center}
\begin{tikzpicture}[xscale=1,yscale=1]
\node at (-1,0) {$\Omega_3(y,x) =$};
\node at (-1,2) {$\Omega_3(x,x) =$};
\node at (-1,-2) {$\Omega_3(y,y) =$};
\node (A) at (1,0) {$w_1$};
\node (B) at (2,0) {$w_2$};
\node (C) at (3,0) {$w_3$};
\node (D) at (4,0) {$w_4$};
\node (E) at (5,0) {$w_5$};
\node (F) at (6,0) {$w_6$};
\node (G) at (7,0) {$w_7$};
\node (H) at (1.5,2) {$u_1$};
\node (I) at (3.5,2) {$u_2$};
\node (J) at (5.5,2) {$u_3$};
\node (K) at (2.5,-2) {$v_1$};
\node (L) at (4.5,-2) {$v_2$};
\node (M) at (6.5,-2) {$v_3$};

\draw [->] (1.5,1.75)--(1,0.25) node[midway,left]{$\alpha$};
\draw [->] (3.5,1.75)--(3,0.25) node[midway,left]{$\alpha$};
\draw [->] (5.5,1.75)--(5,0.25) node[midway,left]{$\alpha$};
\draw [->] (1.5,1.75)--(2,0.25) node[midway,right]{$\beta$};
\draw [->] (3.5,1.75)--(4,0.25) node[midway,right]{$\beta$};
\draw [->] (5.5,1.75)--(6,0.25) node[midway,right]{$\beta$};

\draw [->] (2.5,-1.75)--(2,-0.25) node[midway,left]{$\bar\alpha$};
\draw [->] (4.5,-1.75)--(4,-0.25) node[midway,left]{$\bar\alpha$};
\draw [->] (6.5,-1.75)--(6,-0.25) node[midway,left]{$\bar\alpha$};
\draw [->] (2.5,-1.75)--(3,-0.25) node[midway,right]{$\bar\beta$};
\draw [->] (4.5,-1.75)--(5,-0.25) node[midway,right]{$\bar\beta$};
\draw [->] (6.5,-1.75)--(7,-0.25) node[midway,right]{$\bar\beta$};

\end{tikzpicture}
\end{center}
The connectivity of the underlying graph of arrows that describe this biset shows that it is indecomposable. It is representable on both sides because on considering either the left action of $\calC$ or the right action it is a disjoint union of the copies of the two types of representable set for each side. There are infinitely many of these indecomposable birepresentable bisets. As a consequence, $\End_{\BB_R^{1,1}}(\calC)$ has infinite rank over $R$.
\end{example}

We next mention some examples of biset functors defined on $\BB_R^{\calX,\calY}$ with $\calX,\calY\in\{1,\all\}$.

\begin{example}
All examples of biset functors defined on $\BB=\BB_R^{\all,\all}$ provide examples of biset functors defined on $\BB_R^{\calX,\calY}$ with $\calX,\calY\in\{1,\all\}$, by restricting the functor to the subcategory. The representable functors defined on $\BB_R^{\calX,\calY}$ provide more examples.
\end{example}

\begin{example}
Let $k$ be a field and let $G_0(k\calC)$ be the Grothendieck group of finitely generated $k\calC$-modules, using relations $[B]=[A]+[C]$ whenever $0\to A\to B\to C\to 0$ is a short exact sequence of $k\calC$-modules. This defines a biset functor $G_0(k\calC)$ on $\BB_\ZZ^{\all,1}$ and on $\BB_\ZZ^{1,1}$. The reason this works is that when ${}_\calC\Omega_\calD$ is a biset that is representable for $\calD$ then $k\Omega$ is projective as a right $k\calD$-module (see \cite[Prop. 4.4]{Webb}), so the functor $k\Omega\otimes_{k\calD}-: k\calD\hbox{-mod}\to k\calC\hbox{-mod}$ is exact and defines a homomorphism $G_0(k\calD)\to G_0(k\calC)$.
\end{example}

\begin{example}
Considering only the discrete categories $[n]$ of Example~\ref{discrete-category-bisets}, all bisets for these categories are birepresentable, and so  Proposition~\ref{discrete-category-proposition} holds with $\BB_R$ replaced by $\BB_R^{1,1}$.
\end{example}

We show in the next section that homology and cohomology of the category are biset functors on these modified biset categories, just as happens with biset functors defined on groups.

When we study the internal tensor product and Hom of biset functors we will need the following straightforward lemma about products of representable sets.

\begin{lemma}
\label{representable-product-lemma}
Let $\calC_i,\calD_i$ be categories, $i=1,2$. 
\begin{enumerate}
\item If $\Omega_i$ is a representable $\calC_i$-set, $i=1,2$, then $\Omega_1\times\Omega_2$ is a representable $\calC_1\times\calC_2$-set.
\item If $\Omega_i$ is a $(\calC_i,\calD_i)$-biset that is representable on the left, $i=1,2$, then $\Omega_1\times\Omega_2$ is a  $(\calC_1\times\calC_2,\calD_1\times\calD_2)$-biset representable on the left. The similar result holds for bisets that are representable on the right.
\end{enumerate}
\end{lemma}

\begin{proof}
(1) It suffices to consider the case where each $\Omega_i\cong\Hom_{\calC_i}(x_i,-)$ is representably generated at $x_i$. We readily check that
$$\Omega_1\times\Omega_2\cong\Hom_{\calC_1\times\calC_2}((x_1,x_2),-)$$
is representably generated at $(x_1,x_2)$.
Part (2) follows immediately.
\end{proof}

\section{Realizing homology and cohomology as biset functors}
\label{cohomology-section}
In the classical case of biset functors for groups, group cohomology $H^*(-,R)$ is not a biset functor when all bisets are allowed, because cohomology does not admit an operation of deflation, although it does admit operations of restriction, corestriction and inflation. If we allow only bisets that have a free action on one side then $H^*(-,R)$ is indeed a biset functor. This point of view was exploited in \cite{Weba} to provide a method of computing group cohomology at a prime $p$, in terms of the cohomology of $p$-groups. 

This tells us that the cohomology of a category, defined as $H^*(\calC,R):=\Ext_{R\calC}^*(\underline R,\underline R)$ where $\underline R$ is the constant functor on $\calC$, will not be a biset functor on categories if we allow all bisets as morphisms, because this is not the case when the categories are groups. We show in this section that it is indeed a biset functor defined on $\BB_R^{\all,1}$ (at least when $R$ is a field) and that homology $H_*(\calC,R)$ is a biset functor defined on $\BB_R^{1,\all}$ (without restriction on $R$). This extends the known situation for groups to arbitrary finite categories.

The most accessible way to show that cohomology of groups is a biset functor on groups is to observe that the biset category $\BB^{\all,1}$ with groups as objects is generated by standard bisets that correspond to restriction, corestriction, conjugation and inflation, with certain defining relations between them, and that these relations are also satisfied by the corresponding operations on cohomology. The problems with this approach for biset functors on categories are that we do not have an appropriate set of standard bisets that generate the biset category and, even if we did, we do not know what defining relations hold between different candidate bisets. We take a different approach, looking for a functorial definition of the cohomology operations that is not specific to structures like groups. Such an approach is to observe first that Hochschild homology is a biset functor; and then to show that ordinary homology of a category is naturally a direct summand of its Hochschild homology. This follows the work of Xu~\cite{Xu2}, who did the same thing for cohomology. He did not write out the argument for homology, although the ideas are exactly the same. Because there are some technicalities, we write out the argument here.

See \cite{Webb} for background on the homology and cohomology of a category. As well as having an algebraic definition, these groups are also the homology and cohomology in the topological sense of the nerve of the category. The generality also includes the homology and cohomology of simplicial complexes. Given a simplicial complex, we may form the poset of its simplices, whose nerve is the barycentric subdivision of the original complex. Maps of simplicial complexes give order preserving maps of these posets, which are functors when the posets are regarded as categories. The functor $\phi:\Cat\to\BB^{1,\all}$ that we have considered now realizes maps of simplicial complexes as morphisms in the biset category, and applying the homology biset functor we get the usual homology of spaces.

\subsection{Hochschild homology as a biset functor}
Because the category algebra $R\calC$ is free as an $R$-module, we may take the Hochschild homology of the finite category $\calC$ over $R$ to be
$$
HH_*(R\calC)=\Tor_*^{R\calC^e}(R\calC, R\calC)
$$
where $R\calC^e:= R\calC\otimes_R R\calC^\op$ is the \textit{enveloping algebra} of the category algebra $R\calC$. The next result is immediate from work of Bouc and Keller.

\begin{theorem}[Bouc \cite{Bou4}, Keller \cite{Kel}]
\label{hochschild-homology-as-biset-functor}
 Hochschild homology $HH_*(R\calC)$ has the structure of a biset functor on $\BB_R^{1,\all}$.
\end{theorem}

Keller's approach is explained more recently in \cite[Theorem 2.1]{AK} and Bouc's approach is explained more recently in \cite[5.8.1]{Zim}.

\begin{proof}
Bouc and Keller both show that Hochschild homology is a functor $\calP\to R\hbox{-mod}$ where $\calP$ is the category whose objects are $R$-algebras and where the morphisms are finite complexes of bimodules, perfect on one side. Specifically, if $A$ and $B$ are $R$-algebras, we will take a morphism $B\to A$ in $\calP$ to be an object in the derived category of $(A,B)$-bimodules (that is, a complex of $A\otimes_R B^\op$-modules) that is perfect on the left: as $A$-modules the modules are all finitely generated projective. There is a functor $\BB_R^{1,\all}\to\calP$ that sends a category $\calC$ to its category algebra $R\calC$, and a biset $\Omega$ to the bimodule $R\Omega$, regarded as a complex concentrated in degree 0. This is projective on the left if $\Omega$ is representable on the left. Composing these two functors establishes the result. 
\end{proof}

\subsection{Ordinary homology and cohomology as biset functors}
We now set about showing that ordinary homology of a category is naturally a direct summand of its Hochschild homology.
For this we review the functorial dependence of $\Tor_*^{R\calC}$ on $\calC$. It is more complicated with categories than with rings, because a functor $F:\calC\to\calD$ gives rise to an algebra homomorphism $R\calC \to R\calD$ only when $F$ is one-to-one on objects \cite[2.2.3]{Xu1}, and we have to take account of this. Given a functor $F:\calC\to\calD$ there is a restriction functor $\Res_F: R\calD\mod\to R\calC\mod$ given by composition with $F$ when representations are viewed as functors. It has a left adjoint, which is the left Kan extension $\LK_F:R\calC\mod\to R\calD\mod$; see \cite[2.3]{Xu2} for details in this context. The restriction functor is exact and so the left Kan extension sends projective objects to projective objects. 

We will regard $R\calC$-modules $N$ not just as functors defined on $\calC$, but also as $R$-modules $\bigoplus_{x\in\Ob\calC} N(x)$ on which the ring $R\calC$ acts. We distinguish between the $R$-module $\bigoplus_{x\in\Ob\calC} (\Res_F\LK_F N)(x)$ and the $R$-module $\bigoplus_{y\in\Ob\calD} (\LK_F N)(y)$.

The unit of the adjunction for the left Kan extension is a natural transformation whose value at $N$ is $\eta^N: N\to \Res_F\LK_F N$. We now construct a 
map of $R$-modules
 $\bar\eta^N:N\to\LK_F N$, natural with respect to $R\calC$-module homomorphisms,
in the following way: for each object $x\in\calC$ and $u\in N(x)$ we have $(\Res_F\LK_FN)(x)=(\LK_F(N))F(x)$, so $\eta_x^N(u)\in \LK_FN(F(x))$. We take these maps $\eta_x^N$ as components of a map
$$
\bar\eta^N: \bigoplus_{x\in\Ob\calC} N(x)\to \bigoplus_{y\in\Ob\calD} (\LK_FN)(y)
$$
which is the desired map of $R$-modules.  The difference between $\eta^N$ and $\bar\eta^N$ is that in $\bar\eta^N$ we sum the values of $\eta^N$ on the objects in the preimage of each object $y$ of $\calD$. 

\begin{example}
Let $\calC$ be a discrete category with two objects $\{x,y\}$ and only identity morphisms, and let $\calD=\one$ be the category with one object and one morphism. Let $F:\calC\to\calD$ be the unique possible functor, and let $N$ be the representation of $\calC$ that is $R$ on $x$ and $0$ on $y$. Then $\LK_F N$ is the constant functor on $\calD$, and $\Res_F\LK_F N=\underline R_\calC$ is the constant functor on $\calC$. Now $\eta:N\to \Res_F\LK_F N$ identifies as a map $R^2\to R^2$ with matrix $\begin{bmatrix}1&0\cr 0&0\cr  \end{bmatrix}$ and $\bar\eta:R^2\to R$ has matrix $\begin{bmatrix}1&0\cr  \end{bmatrix}$
\end{example}

\begin{proposition}
\label{map-on-tensor-product}
Let $F:\calC\to\calD$ and $G:\calD\to \calE$ be functors and let $\eta:1\to \Res_F\LK_F$ and $\theta:1\to \Res_G\LK_G$ be the units of the corresponding adjunctions. Let $N$ be a left $R\calC$-module and $M$ a right $R\calC$-module.
\begin{enumerate}
    \item  The map $\bar\eta^N$ satisfies  $\bar\eta^N(\alpha(u)) = F(\alpha)\bar\eta^N(u)$ for all morphisms $\alpha$ in $\calC$ and elements $u$ of $N$.
    \item We have $\LK_G\LK_F \simeq \LK_{GF}$ and $\overline{\theta\eta}=\bar\theta\bar\eta$.
    \item  There is an $R$-module homomorphism $$\bar \eta^M\otimes \bar\eta^N: M\otimes_{R\calC}N\to \LK_F M\otimes_{R\calD} \LK_F N$$
    natural in both $M$ and $N$, given by $u\otimes_{R\calC} v\mapsto \bar\eta^M (u)\otimes_{R\calD} \bar\eta^N(v)$.
\end{enumerate}
\end{proposition}

\begin{proof}
(1) This follows from naturality of $\eta^N$.

(2) comes from the universal property of the left Kan extension and the fact that, for an object $z$ of $\calE$, the preimage under $GF$ of $z$ in $\calC$ is the disjoint union of the preimages under $F$ of objects $y\in\Ob\calD$ of the preimage of $z$ under $G$.

For (3), we observe that the mapping
$$
M\times N\to \LK_F M\otimes_{R\calD} \LK_F N
$$
given by $(u,v)\mapsto \bar\eta^M (u)\otimes_{R\calD} \bar\eta^N(v)$ is balanced for $R\calC$.
\end{proof}

\begin{corollary}
\label{Tor-functor}
Let $F:\calC\to\calD$ be a functor and let $M, N$ be right and left $R\calC$-modules. There is a morphism
$$
F_*:\Tor_*^{R\calC}(M,N)\to \Tor_*^{R\calD}(\LK_FM,\LK_FN)
$$ that is functorial in F and natural in both $M$ and $N$.
\end{corollary}

\begin{proof}
Take a resolution $\calP_M\to M$ of $M$ by projective right $R\calC$-modules. Applying $\bar\eta$ to all the terms, we get a map of complexes $\calP_M\to \LK_F\calP_M$ that realizes $\bar\eta^M$ on zero homology. We denote this map of complexes also by $\bar\eta^M$. Because $\LK_F$ is the left adjoint of an exact functor, $\LK_F\calP_M$ is a complex of modules that are projective for $R\calD$. Taking a projective resolution $\calQ\to\LK_F M$ over $R\calD$, the identity map on $H_0(\LK_F\calP_M)=H_0(\calQ)=\LK_FM$ lifts to a morphism of complexes $\lambda:\LK_F\calP_M\to\calQ$, because the first has projective terms and the second is acyclic. We obtain a morphism $(\lambda\circ\bar\eta^M)\otimes \bar\eta^N: \calP_M\otimes_{R\calC}N\to \calQ\otimes_{R\calD}LK_FN$ as in Proposition~\ref{map-on-tensor-product}. This gives rise to the desired morphism $\Tor_i^{R\calC}(M,N)\to \Tor_i^{R\calD}(\LK_F M,\LK_F N)$ for each $i\ge 0$ on taking homology.
\end{proof}

We apply this functorial action on $\Tor$ to the homology version of the situation considered by Xu in \cite{Xu2} where he relates the ordinary cohomology of a category $\calC$ to the Hochschild cohomology of $R\calC$. We explain how Xu's approach applies to homology. Xu defines the \textit{enveloping category} of $\calC$ to be the category $\calC^e=\calC\times\calC^\op$ and he observes that the category algebra $R[\calC^e]$ may be identified with the enveloping algebra of $R\calC$ thus: $R[\calC^e]\cong (R\calC)^e= R\calC\otimes_R (R\calC)^\op$. He makes use of Quillen's \textit{factorization category} $F(\calC)$ that has as its objects the morphisms of $\calC$. To avoid confusion, if $\alpha:x\to y$ is a morphism in $\calC$ we write $[\alpha]$ to denote the corresponding object of $F(\calC)$. If $\alpha':x'\to y'$ is another morphism in $\calC$, a morphism $[\alpha]\to[\alpha']$ in $F(\calC)$ is a pair of morphisms $(u,v)$ in $\calC$ making the following diagram commute:
$$
\diagram{y&\umapleft{\alpha}&x\cr
\lmapdown{u}&&\rmapup{v}\cr
y'&\umapleft{\alpha'}&x'\cr
}
$$
In other words, there is a morphism from $[\alpha]$ to $[\alpha']$ if and only if $\alpha'=u\alpha v$ for some morphisms $u,v$ in $\calC$ or, equivalently, $\alpha$ is a factor of $\alpha'$ in $\calC$. Composition of morphisms is achieved by vertical juxtaposition of diagrams and composition of morphisms vertically. The category $F(\calC)$ admits two natural covariant functors to $\calC$ and $\calC^\op$
$$
\calC \xleftarrow{t} F(\calC) \xrightarrow{s} \calC^\op,
$$
where $t$ and $s$ send an object $[\alpha]$ to its target and source, respectively.

We now have a commutative triangle of functors:
$$
\diagram{F(\calC)&\umapright{\tau}&\calC^e\cr
&\searrow^t&\rmapdown{\textrm{pr}}\cr
&&\calC\cr}
$$
where $\tau=(t,s)$ and $\textrm{pr}$ is projection onto the first factor. On page 1880 of \cite{Xu2} Xu establishes that $\LK_\tau\underline R_{F(\calC)}\cong R\calC$ as $R\calC^e$-modules and $\LK_{\textrm{pr}}R\calC\cong \underline R_\calC\cong \LK_t\underline R_{F(\calC)}$ as $R\calC$-modules where, for any category $\calD$, $\underline R_\calD$ (or just $\underline R$) denotes the constant functor on $\calD$.

Recall that, for any category $\calC$, we have an isomorphism
$$
\Tor_*^{R\calC}(\underline R,\underline R)\cong H_*(|\calC|;R)
$$
between the algebraically defined $\Tor$ groups of the constant functor and the homology over $R$ of the nerve $|\calC|$ of $\calC$ (see \cite{Webb} for an exposition). Quillen showed in \cite{Qui} that the functor $t:F(\calC)\to\calC$ is a homotopy equivalence between the nerves, and so it induces an isomorphism on homology. 

Putting all this together we get the analogue for homology of a theorem of Xu~\cite{Xu2}.

\begin{theorem}
\label{Hochschild-splitting-theorem}
Let $R$ be a commutative ring.
\begin{enumerate}
    \item The functors $t$, $\tau$ and $\textit{pr}$ induce homomorphisms 
    $$
    H_*(\calC;R)\xrightarrow{\tau_*(t_*)^{-1}} \Tor_*^{R\calC^e}(R\calC, R\calC) \xrightarrow{\textit{pr}_*}H_*(\calC;R)
    $$
    whose composite is the identity.
    \item $\Tor_*^{R\calC^e}(R\calC, R\calC)\cong HH_*(R\calC)$ is the Hochschild homology of $R\calC$.
\end{enumerate}
\end{theorem}

\begin{proof}
(1) We have constructed the functorial maps $t_*$, $\tau_*$ and $\textit{pr}_*$ in Corollary~\ref{Tor-functor} and the surrounding description. In applying Corollary~\ref{Tor-functor} we take $M=\underline R$ and $N=\underline R$ to be the right and left constant functors on $\calC$ and now the $\Tor$ groups in that corollary are $H_*(\calC;R)$. Xu's identification of the left Kan extensions gives the remaining two terms in the sequence. The fact that $\textit{pr}\circ \tau = t$ is a homotopy equivalence means that $\textit{pr}_*\circ\tau_*\circ (t_*)^{-1}$ is the identity. 

(2) holds because the category algebra $R\calC$ is projective as an $R$-module and so Hochschild homology can be computed as the $\Tor$ group.
\end{proof}

\begin{corollary}
\label{ordinary-summand-of-Hochschild}
$H_*(\calC;R)$ is a direct summand of $HH_*(R\calC)$ in a way that commutes with the functorial action of bisets in $\BB^{1,\all}$.
\end{corollary}

The statement about the action of bisets means that if ${}_\calC\Omega_\calD$ is a biset  in $\BB^{1,\all}$ then the morphism $HH_*(R\calD) \to HH_*(R\calC)$ determined by $\Omega$ maps the summand $H_*(\calD;R)$ to the summand $H_*(\calC;R)$.

\begin{proof}
The fact that $H_*(\calC;R)$ is a direct summand is immediate from Theorem~\ref{Hochschild-splitting-theorem}. The fact that the summand is preserved by the biset action comes from an examination of the functorial effect of the bimodule $R\Omega$ on Hochschild homology and its interaction with the functors $t,\tau$ and $\textit{pr}$.
\end{proof}

We are now ready to prove that homology and cohomology are biset functors. For any ring $A$ and right $A$-module $M$ we define the left $A$-module $M^\vee:=\Hom_A(M,A_A)$. As is well known, this provides a duality between the full subcategory of right $A$-modules whose objects are finitely generated projective and  the full subcategory of left $A$-modules whose objects are finitely generated projective, as can easily be seen by reducing to the case $M=A_A$, in which case $M^\vee\cong {}_AA$. 

\begin{theorem}
\label{homology-as-biset-functor}
Let $R$ be a commutative ring. The ordinary homology of a category $H_*(\calC;R)$ has the structure of a biset functor on $\BB_R^{1,\all}$.
When $R$ is a field, $H^*(\calC;R)$ is a biset functor on $\BB_R^{\all,1}$.
\end{theorem}

\begin{proof}
The homology statement follows from  Theorem~\ref{hochschild-homology-as-biset-functor} and Corollary~\ref{ordinary-summand-of-Hochschild}. If ${}_\calC\Omega_\calD$ is a $(\calC,\calD)$-biset, representable for $\calC$, we define its effect on $H_i(\calD;R)$ to be the composite
$$
    H_*(\calD;R)\xrightarrow{\tau_*(t_*)^{-1}} \Tor_*^{R\calD^e}(R\calD, R\calD)
    \xrightarrow{R\Omega_*}
    \Tor_*^{R\calC^e}(R\calC, R\calC)
    \xrightarrow{\textit{pr}_*}H_*(\calC;R)
$$
where $R\Omega_*$ denotes the operation on $\Tor$ groups defined by Bouc and Keller for the $(R\calC,R\calD)$-bimodule $R\Omega$, which is projective for $R\calC$.

For cohomology we exploit the fact that the homology and cohomology of $\calC$ are the homology and cohomology of a space (the nerve $|\calC|$). When $R$ is a field it follows that cohomology $H^*(\calC,R)\cong (H_*(\calC;R))^*$ is the vector space dual of homology, by the universal coefficient theorem (see \cite[Theorem 53.5]{Mun}).  To define cohomology on $\BB_R^{\all,1}$ let ${}_\calD\Omega_\calC$ be a biset that is representable on the right, as a $\calC$-set. The $(R\calD,R\calC)$-bimodule $R\Omega$ is finitely generated projective as a right $R\calC$-module, so that the $(R\calC,R\calD)$-bimodule $R\Omega^\vee$ is finitely generated projective as a left $\calC$-module and defines a mapping $H_*(\calD,R)\to H_*(\calC,R)$, by the result for homology. Taking the vector space dual we have a mapping $(H_*(\calC,R))^*\to (H_*(\calD,R))^*$. This provides the desired operation on cohomology
\end{proof}

The construction provides analogues of restriction, corestriction, and deflation (on homology) or inflation (on cohomology) for category cohomology, in the manner of group cohomology. In that context, when $H$ is a subgroup of a group $G$, the $(H,G)$-biset ${}_HG_G$ determines the restriction maps $H_*(G,R)\to H_*(H,R)$ and $H^*(G,R)\to H^*(H,R)$ (also called transfer in the case of homology), whereas the biset ${}_GG_H$ determines the corestriction maps $H_*(H,R)\to H_*(G,R)$ and $H^*(H,R)\to H^*(G,R)$ (also called transfer in the case of cohomology). This is explained in Example 2.6 of  \cite{Lin} and the subsequent discussion there.

We present examples to illustrate the kind of circumstances where we may construct restriction and corestriction maps.

\begin{example}
When $F:\calD\to\calC$ is a functor we obtain a restriction map in cohomology $\Res = F^*:H^*(\calC,R)\to H^*(\calD,R)$ and a corestriction map $\Cores=F_*: H_*(\calD,R)\to H_*(\calC,R)$ simply by functoriality of cohomology and homology. We may also derive this from our theory because ${}_\calC\calC_{{}^F\calD}$ is representable on the left and ${}_{\calD^F}\calC_\calC$ is representable on the right, by Proposition~\ref {regular-is-representable}. The same proposition implies that when $F$ has a left adjoint then there is a corestriction map in cohomology, and if $F$ has a right adjoint there is a restriction map in homology.
\end{example}

\begin{example}
Consider the inclusion of categories $\calD\to\calC$ of Example~\ref{birepresentable-example}. It was shown there in Proposition~\ref{C4C2-category-proposition} that both ${}_\calC\calC_\calD$ and ${}_\calD\calC_\calC$ are birepresentable, so that $H^*$ and $H_*$ are defined on these bisets giving a transfer map $H^*(\calD,R)\to H^*(\calC,R)$ (when $R$ is a field) and a restriction map $H_*(\calC)\to H_*(\calD)$. In this example, note that the category algebra $R\calC$ is not symmetric, and also that the left and right adjoints of the restriction map $R\calC\mod\to R\calD\mod$ are distinct. This is in contrast to the situation in group cohomology and some other constructions of transfer maps, such as  in   \cite{Lin}, where coincidence of the left and right adjoints of restriction is used.
\end{example}

\begin{example}
We find examples of transfer maps in the context of the \textit{Grothendieck construction} of a finite group $G$ acting on a finite category $\calC$. This includes the case of transporter categories described in \cite{Xu3}, which are the special case where $\calC$ is a poset. We refer to  \cite{Xu3} for further details on the Grothendieck construction and its connection with equivariant cohomology and group graded algebras

Let the finite group $G$ act on a finite category $\calC$. Such an action can be described by regarding $G$ as a category $\calG$ with a single object $*$, where the endomorphisms of $*$ are the set $G$, and by specifying a functor $F:\calG\to \hbox{Finite Categories}$ with $F(*)=\calC$. When $u$ is an object of $\calC$ we shorten the notation for the action of $G$ on $\calC$ to $gu:=F(g)(u)$. In this situation the Grothendieck construction $\calC\rtimes G$ may be taken to be the category with the same objects as $\calC$ and, if $u,v$ are objects of $\calC$, then $\Hom_{\calC\rtimes G}(u,v)$ is the set of pairs $(\theta,g)$ where $g\in G$ and $\theta:F(g)(u)\to v$ is a morphism in $\calC$. If $(\psi,h):v\to w$ is a second morphism then composition is defined to be $(\psi,h)(\theta,g):= (\psi\cdot F(h)(\theta),hg)$.

\begin{proposition}
Let $G$ be a finite group acting on a finite category $\calC$. Let $H$ be a subgroup of $G$. Then the bisets
$$
{}_{\calC\rtimes H}( \calC\rtimes G) _{\calC\rtimes G}
\quad\hbox{and}\quad
{}_{\calC\rtimes G}( \calC\rtimes G) _{\calC\rtimes H}
$$
are birepresentable. As a consequence we obtain restriction and transfer maps between $H_*( \calC\rtimes G, R)$ and $H_*(\calC\rtimes H, R)$, as well as between $H^*( \calC\rtimes G, R)$ and $H^*(\calC\rtimes H, R)$.
\end{proposition}

\begin{proof}
To show that ${}_{\calC\rtimes H} \calC\rtimes G$ is representable as a left $\calC\rtimes H$-set we show that $\Hom_{\calC\rtimes G}(u,-)$ is representable for each object $u\in\calC$. As $\calC$-sets we have
$$
\Hom_{\calC\rtimes G}(u,-)\cong \bigsqcup_{g\in G}\Hom_\calC(gu,-)
$$
via a natural isomorphism that at $v$ sends a morphism $(\theta,g)$, where $\theta:gu\to v$, to $\theta$ in the component indexed by $g$. From this we see that as a $\calC\rtimes H$-set it is $\bigsqcup_{g\in [H\backslash G]}\Hom_{\calC\rtimes H}(gu,-)$ where $[H\backslash G]$ denotes a set of representatives for the right cosets of $H$ in $G$. This is representable. 

We now show that $ \calC\rtimes G_{\calC\rtimes H}$ is representable on the right. As $\calC$-sets, we have for each object $v\in\calC$ 
$$
\Hom_{\calC\rtimes G}(-,v)\cong \bigsqcup_{g\in G}\Hom_\calC(-,g^{-1}v)
$$
via a natural isomorphism that at $u$ sends a morphism $(\theta,g)$, where $\theta:gu\to v$, to $F(g^{-1})(\theta):u\to g^{-1}v$ in the component indexed by $g$. This is isomorphic to  $\bigsqcup_{g\in [G/H]}\Hom_{\calC\rtimes H}(-,g^{-1}v)$ as a right $\calC\rtimes H$-set, and is representable.
\end{proof}
\end{example}

\section{Simple biset functors}
We first review the general properties of simple functors on linear categories when they are restricted to full subcategories. We then characterize simple functors in the manner of \cite{TW}, we show how to parametrize simple functors, and we provide a way to calculate their values. This approach copies what happens with biset functors on groups, but there are some differences. For example, simple biset functors on groups are parametrized by pairs $(H,V)$ where $H$ is a group and $V$ is a simple representation of $\Out(H)$. With categories, we will see that it is no longer the case that simple biset funtors are parametrized by all the pairs $(\calC,V)$, where $\calC$ ranges through the isomorphism classes of objects of $\BB$ and $V$ ranges through the isomorphism types of simple representations of $\Out(\calC)$. We point out other differences as well.

\begin{definition}
A biset functor is said to be simple if it has no subfunctors apart from the whole functor and the zero functor. 
\end{definition}

Our general results about simple functors apply to all of the categories $\BB_R^{\calX,\calY}$ regardless of the choice of $\calX,\calY\in\{1,\all\}$. We will write $\BB$ to denote one of these categories, suppressing the notation for $\calX,\calY$ and $R$.  A simple biset functor is always annihilated by a maximal ideal of $R$, so without loss of generality we will assume that $R$ is a field when considering them. 

\subsection{Restricting and extending simple biset functors between full subcategories}

We recall the basic relationship between simple representations of categories and of their full subcategories. This is described in the context of representations of categories in \cite[Props. 3.2 and 4.1]{Webb}, rather than $R$-linear representations of $R$-linear categories, which is what we need here. Furthermore, in \cite{Webb} the restriction was made to full subcategories with {finitely many} objects. In the proof of Proposition~\ref{simple-restriction-extension} we explain the small modifications needed to apply to our present situation.

\begin{proposition}
\label{simple-restriction-extension}
Let $S$ be a simple $R$-linear representation of an $R$-linear category $\calC$, and let $\calD$ be a full subcategory of $\calC$.
\begin{enumerate}
    \item The restriction $S\downarrow_\calD^\calC$ of $S$ to $\calD$ is either a simple representation of $\calD$ or zero.
    \item Every simple representation of $\calD$ arises uniquely as the restriction of a simple representation of $\calC$.
    \item For every object $x$ of $\calC$ the evaluation $S(x)$ is either a simple $\End_\calC(x)$-module or zero.
    \item If $T$ is any simple representation of $\calC$ over $R$ and $x$ is an object of $\calC$ for which $T(x)\cong S(x)\ne 0$ as $R\End_\calC(x)$-modules, then $S\cong T$.
\end{enumerate}
\end{proposition}

\begin{proof}
Parts (3) and (4) are immediate corollaries of parts (1) and (2), on taking $\calD$ to be the full subcategory whose only object is $x$. This is because saying that a representation of such $\calD$ is simple is the same as saying that its value at $x$ is a simple module for $\End_\calD(x)$.

Proofs of (1) and (2) are given in the proof of Proposition 3.2 (2) of \cite{Webb}. That proof made the assumption that we were considering all representations of categories, rather than $R$-linear representations, but the exact words used there apply to the linearized version as well. The other difference with \cite{Webb} is that it was supposed that $\calD$ has finitely many objects, this allowing notation using idempotents in the category algebra. 

If we wish we can avoid the notation used in \cite{Webb} to make the proof apply generally, but we can also deduce the general case from the finite case as follows. Assume the result when $\calD$ has finitely many objects, so that parts (3) and (4) may be assumed proved. We deduce (1): if $S\downarrow_\calD^\calC$ were to have a proper subfunctor $M$, there would be an object $x$ of $\calD$ with $M(x)\ne S(x)$, from which $M(x)=0<S(x)$ by simplicity of $S(x)$, and another object $y$ in $\calD$ where $M(y)\ne 0$ so that $M(y)=S(y)$ by simplicity of $S(y)$. This is not possible, because taking $\calE$ to be the full subcategory with the two objects $\{x,y\}$, the restriction of $S$ to $\calE$ is simple (by the case with finitely many objects), and it is non-zero at both $x$ and $y$, but it has a subfunctor that is zero at $x$ but not at $y$. This proves (1) in generality.

To deduce (2) from the finite case, let $T$ be a simple representation of $\calD$ and let $x$ be an object of $\calD$ with $T(x)\ne 0$. There is a unique simple representation $S$ of $\calC$ with $S(x)\cong T(x)$, by (4). Now $S$ restricts to $\calD$ as $T$ because, for each two-object full subcategory $\calE$ of $\calD$ with objects $\{x,y\}$ the restrictions of $S$ and $T$ to $\calE$ are the unique simple extension from  $\{x\}$ to $\calE$. 
\end{proof}

We summarize this situation in the context of biset categories.

\begin{corollary}
\label{simple-evaluations}
Let $S$ be a simple biset functor defined on one of the biset categories $\BB_R^{\calX,\calY}$ with $\calX,\calY\in\{ 1,\all\}$, and let $\BB'$ be a full subcategory.  
\begin{enumerate}
\item The restriction of $S$ to $\BB'$ is either zero or a simple functor. 
\item This establishes a bijection between isomorphism types of simple biset functors on $\BB$ that are non-zero on $\BB'$, and simple functors on $\BB'$. 
\item The list of pairs  $\{(\calC,S(\calC))\bigm| \calC\hbox{ is a category, } S(\calC)\ne 0\}$ determined by $S$ has the properties:
\begin{itemize} 
\item each $S(\calC)$ is a simple module for $\End_{ \BB_R^{\calX,\calY}}(\calC)$, and 
\item $S$ is completely determined by any one of the pairs $(\calC,S(\calC))$ where $S(\calC)\ne 0$.
\end{itemize}
\item For each category $\calC$, the number of isomorphism types of simple functors $S$ with $S(\calC)\ne 0$ is equal to the number of isomorphism types of simple modules for $\End_{ \BB_R^{\calX,\calY}}(\calC)$.
\end{enumerate}
\end{corollary}

These properties of the restriction of simple functors to full subcategories, and the extension of simple functors from full subcategories raise interesting questions about the relationship between the simple functors on different categories. For example, every simple biset functor defined on finite groups  extends uniquely to a simple biset functor defined on all finite categories, which then restricts to a simple biset functor or zero on (for example) posets. Equally, every simple biset functor defined on posets extends to a simple functor on all categories which then restricts to a simple functor or zero on groups. What are the simple functors that extend and restrict like this in a non-zero fashion?

\begin{example}
\label{discrete-categories-simple}
Consider the discrete categories $[m]$, for which we established that $\End_\BB([m]) \cong \Mat_{m,m}($R$)$ in Proposition~\ref{discrete-category-proposition}. When $R$ is a field these endomorphism rings each have only one simple module, and so, for each $[m]$, there is (up to isomorphism) only one simple biset functor that does not vanish on it, by Corollary~\ref{simple-evaluations}(4). This applies, in particular, to the category $\one=[1]$ and, in the context of groups, the simple functor that does not vanish on $[1]$ is denoted $S_ {1,R}$. We claim that this is also the simple functor that does not vanish on $[m]$, for every $m$. We argue that the identity endomorphism $1_{[1]}$ of $[1]$ in $\BB$ factors as 
$$1_{[1]}=\begin{bmatrix}1&0&\cdots&0\\
\end{bmatrix}
\begin{bmatrix}
1\\
0\\
\vdots\\
0\\
\end{bmatrix}$$
It acts in a non-zero fashion on $S_ {1,R}([1])$, so the morphism $S_ {1,R}([m])\to S_ {1,R}([1])$ induced by the biset $\begin{bmatrix}1&0&\cdots&0\\
\end{bmatrix}$ is also non-zero, and hence  $S_ {1,R}([m])\ne 0$. Thus $S_ {1,R}$ is the unique simple functor that does not vanish on $[m]$.
\end{example}

\subsection{A characterization of simple biset functors}
The characterization we describe is similar to that of \cite{TW}. It underlies the method we shall give in Subsection~\ref{calculation-subsection} for calculating the values of simple biset functors. We continue to write  $\BB$ for any of the categories $\BB_R^{\calX,\calY}$ with $\calX,\calY\in\{1,\all\}$.

\begin{definition}
Let $\calX$ be a collection of categories and $M$ a biset functor. 
We define $I_\calX M$ to be the subfunctor of $M$ \textit{generated} by its values at categories in $\calX$. This means that it is the smallest subfunctor of $M$ whose values on categories $\calC\in\calX$ are $M(\calC)$. It is thus the smallest subfunctor  $I_\calX M$ of $M$ for which the inclusion $I_\calX M\to M$ is an isomorphism on the categories $\calC\in\calX$. We say $M$ itself is \textit{generated} at $\calX$ if $M=I_\calX M$.

 Dually, we define $R_\calX M$ to be the largest subfunctor of $M$ so that the quotient map  $M\to M/R_\calX M$ is an isomorphism on the categories $\calC\in\calX$. We say that $M/R_\calX M$ is the quotient of $M$ \textit{cogenerated} by its values in $\calX$, and $M$ itself is \textit{cogenerated} by its values in $\calX$ if $M=M/R_\calX M$ (so $R_\calX M=0$).
\end{definition}

Analogues of $R_\calX M$ and $I_\calX M$ play an important role in the theory of Mackey functors and biset functors for groups, for instance in induction theorems (see \cite{Webf}), and in the stratification of Mackey functors \cite{Webg, Webc}.

\begin{proposition}
Let $\calX$ be a collection of categories, $\calD$ a category, and $M$ a biset functor. Then
$$
I_\calX M(\calD)=\sum_{\calC\in\calX,\; \Omega\in\Hom_\BB(\calC,\calD)} \Im M(\Omega), 
$$
and 
$$
R_\calX M(\calD)=\bigcap_{\calC\in\calX,\; \Omega\in\Hom_\BB(\calD,\calC)}\Ker M(\Omega).
$$
\end{proposition}

Another way to write the first equation is
$$
I_\calX M(\calD)=\Hom_\BB(\calC,\calD)\cdot  M(\calC).
$$

\begin{proof}
We check first that the two specifications on the right do give subfunctors of $M$. This is a question of checking, for an element $M(\Omega)(u)$ with $u\in M(\calC)$, that if $\Psi:\calD\to\calE$ is a biset then $M(\Psi)(M(\Omega)(u))$ has the form $M(\Omega')(u)$ for some $\Omega': \calC\to\calE$. It does on taking $\Omega'=\Psi\Omega$. This shows that $\sum_{\calC\in\calX,\; \Omega\in\Hom_\BB(\calC,\calD)} \Im M(\Omega)$ defines a subfunctor of $M$. 

Similarly, for $v\in M(\calD)$, if $M(\Omega)(v)=0$ for all $\Omega:\calD\to\calC$ with $\calC\in \calX$ then $M(\Omega')(M(\Psi)(v))=0$ for all $\Omega':\calE\to\calC$, because it equals $M(\Omega'\Psi)(v))$, which is zero on taking $\Omega=\Omega'\Psi$. This shows that $\bigcap_{\calC\in\calX,\; \Omega\in\Hom_\BB(\calD,\calC)}\Ker M(\Omega)$ defines a subfunctor of $M$.

Now $\sum_{\calC\in\calX,\; \Omega\in\Hom_\BB(\calC,\calD)} \Im M(\Omega)$ contains the values of $M$ at $\calC\in \calX$, so is contained in $I_\calX M(\calD)$, and it is necessarily contained in every subfunctor of $M$ containing these values, so we have equality.

The argument for $R_\calX M$ is similar. Write
$$
L(\calD):=\bigcap_{\calC\in\calX,\; \Omega\in\Hom_\BB(\calD,\calC)}\Ker M(\Omega).
$$
If $\calC\in\calX$ then $L(\calC)=0$ because the intersection includes $\Ker M({}_\calC\calC_\calC)=0$, so $M\to M/L$ is an isomorphism on $\calC$ in $\calX$. Thus $L\subseteq R_\calX M$. On the other hand, if $N$ is a subfunctor of $M$ so that $M\to M/N$ is an isomorphism on $\calC\in\calX$ then $N(\calC)=0$ for such $\calC$. Thus $N(\calD)\subseteq \Ker M(\Omega)$ for all $\Omega\in\Hom_\BB(\calD,\calC)$. This means $N\subseteq L$ and shows that $L = R_\calX M$.
\end{proof}

The following characterization of simple biset functors is similar to the characterization of simple Mackey functors that appeared in \cite{TW}, and the proof is the same.

\begin{theorem}
\label{simplicity-criterion}
Let $M$ be a biset functor. Then $M$ is a simple biset functor if and only if there exists a category $\calC$ so that
\begin{enumerate}
\item $M(\calC)$ is a simple $\End_\BB(\calC)$-module,
\item $M$ is generated by its value at $\calC$; that is, $I_{\{\calC\}} M = M$, and
\item $M$ is cogenerated by its value at $\calC$, that is, $R_{\{\calC\}} M = 0$. 
\end{enumerate}
When $M$ is simple, every category $\calC$ for which $M(\calC)\ne 0$ satisfies conditions (1), (2) and (3).
\end{theorem}

Conditions (2) and (3) in this theorem say that for all categories $\calD$,
$$
M(\calD)=\sum_{\Omega\in\Hom_\BB(\calC,\calD)} \Im M(\Omega), 
$$
and
$$
\bigcap_{\Omega\in\Hom_\BB(\calD,\calC)}\Ker M(\Omega) =0.
$$

\begin{proof}
Suppose that $M$ is simple. Let $\calC$ be any category for which $M(\calC)\ne 0$. We have seen in Proposition~\ref{simple-restriction-extension} that $M(\calC)$ is a simple $\End_\BB(\calC)$-module. Because $I_{\{\calC\}}M$ and $R_{\{\calC\}}M$ are subfunctors of $M$, they are either $M$ or $0$. We also know that $I_{\{\calC\}}M(\calC)=M(\calC)$ is non-zero, so we have $I_{\{\calC\}}M=M$. Because $R_{\{\calC\}}M(\calC)=0$ the subfunctor $R_{\{\calC\}}M$ is not the whole of M, so $R_{\{\calC\}}M=0$.

Conversely, suppose the three conditions are satisfied and let $M_1$ be a subfunctor of $M$ that is not zero. Let $\calD$ be a category for which $M_1(\calD)\ne 0$. Condition (3) implies that there is a homomorphism $\Omega:\calD\to \calC$ so that $M_1(\Omega)$ is not the zero homomorphism, so that $M_1(\calC)\ne0$. It follows that $M_1(\calC)=M(\calC)$ by simplicity of $M(\calC)$. Now $M_1=M$ because $M$ is generated by its value at $\calC$. Hence $M$ is simple.
\end{proof}

\subsection{Parametrizing simple biset functors}
Simple representations of categories have been parametrized in many contexts, the result often being that their isomorphism types biject with pairs consisting of an object of the category, and a simple representation of a group or algebra associated to that object. An early form of such a parametrization is the theory of Munn and Ponizovski\v{i} that parametrizes simple modules for finite monoids (see \cite{Ste} for an account). This was extended later to representations of finite categories by Linckelmann and Stolarz~\cite{LS}.  For biset functors the monoid approach does not immediately apply because the biset category is $R$-linear and we consider only $R$-linear-functors, but there are still similarities. A parametrization of the simple biset functors defined on groups was given by Bouc \cite{Bou1} (and in particular cases by Webb  in \cite{Weba}) in terms of  the unique smallest group $G$ on which the simple functor is non-zero, together with a simple representation of the outer automorphism group $\Out(G)$. For each pair $(G,V)$, where $V$ is a simple representation of $\Out(G)$, there is a simple biset functor $S_{G,V}$, and this accounts for all simple biset functors on groups.

A parametrization of simple functors in terms of pairs $(\calC,V)$ where $V$ is a simple representation of $\Out(\calC)$ does not work so well for categories (recall that $\Out(\calC)$ was defined in \ref{out-definition}). Firstly, if $S$ is a simple biset functor, it is not clear what we should mean by a smallest category $\calC$ with $S(\calC)\ne 0$. The issue is that inequivalent categories may be isomorphic in $\BB$ (Theorem~\ref{idempotent-completion-isomorphism-theorem}) and $S$ has isomorphic values on such categories, so which one of these categories should we prefer? 

It is also the case that not all simple representations of $\Out(\calC)$ determine simple biset functors. We see what goes wrong in the next result with the discrete categories $[m]$ of Proposition~\ref{discrete-category-proposition}, where it was shown that $\End_\BB([m])\cong \Mat_{m,m}(R)$. 

\begin{proposition}
\label{discrete-category-proposition2}
For the discrete categories $[m]$ where $m\ge 1$ we have
\begin{enumerate}
\item $\Out([m])\cong S_m$, the symmetric group of degree $m$.
\item $R\Out([m])$ is not a homomorphic image of $\End_\BB([m])$ when $m\ge 2$.
\item When $m\ge 2$ and $V$ is a simple $S_m$-module, there is no simple biset functor $S$ with $S(\calC)=V$.
\end{enumerate}
\end{proposition}

\begin{proof}
Any self equivalence of $[m]$ is a permutation of the objects, and the only natural isomorphisms between them are the identity, because $[m]$ has only identity morphisms. (2) is immediate in view of Proposition~\ref{discrete-category-proposition}. For (3), we have seen in Example~\ref{discrete-categories-simple} that the only simple functor $S$ with $S([m])\ne 0$ is $S_{1,R}$, and $S_{1,R}([m]) = R^m$, which is not a simple module for $S_m$ when $m\ge 2$.
\end{proof}

This shows that simple biset functors $S$ are not parametrized by pairs $(\calC, V)$ where $V$ is a simple module for $\Out(\calC)$ in such a way that $S(\calC)=V$. 
What we can say about parametrizing simple functors is that, according to  Corollary~\ref{simple-evaluations} part (3), they determine and are determined by their list of non-zero evaluations at categories $\calC$, which are simple representations of $\End_\BB(\calC)$. Finding a parametrization of $S$ involves finding a distinguished pair $(\calC,S(\calC))$ where $S(\calC)\ne 0$, and finding a good description of the simple $\End_\BB(\calC)$-modules $S(\calC)$ that can appear.

We resolve these issues by describing a way to parametrize simple functors in terms of \textit{essential algebras} relative to a well ordering of finite categories. This is similar to what has been done elsewhere, but has the difference that we have no compelling natural order on categories. 

\begin{definition}
\label{essential algebras}
Place the isomorphism types (in $\BB$) of finite categories in any well order $<$. For each category $\calC$,  we define the \textit{inessential ideal} $\Iness^<(\calC)$ to be 2-sided ideal in $\End_{\BB}(\calC)$ that is the $R$-span of the bisets that factor through categories $\calD$ with $\calD<\calC$. Thus  $\Iness^<(\calC)$ is the span of bisets of the form $\Omega=\Psi\circ\Theta$ where $\Psi$ is a $(\calC,\calD)$-biset and $\Theta$ is a $(\calD,\calC)$-biset, with $\calD<\calC$. We term the quotient $\Ess^<(\calC):=\End_{\BB}(\calC)/\Iness^<(\calC)$ the \textit{essential algebra} of $\calC$ with respect to $<$.
\end{definition} 

Notice that we could have defined the inessential ideal to be the span of linear combinations of bisets in $\End_\BB(\calC)$ that factor through earlier categories, but this would define the same ideal: if a linear combination of bisets factors through an earlier category, then so must each term in the linear combination.

The essential algebra we have just defined depends on the choice of an arbitrary well order on categories, but in most cases we would choose a well order with the property that if $\calC$ has fewer morphisms than $\calD$ then $\calC < \calD$. 

\begin{proposition}
\label{inessential-zero}
Let $S$ be a biset functor and $\calC$ the least category with respect to the well order $<$ such that $S(\calC)\ne 0$. Then $\Iness^<(\calC)$ acts as 0 on $S(\calC)$, so that the structure of $S(\calC)$ as an $\End_\BB(\calC)$-module is the same as its structure as an $\Ess^<(\calC)$-module. Thus if $S$ is a simple biset functor, then $S(\calC)$ is a simple $\Ess^<(\calC)$-module.
\end{proposition}

\begin{proof}
We know from Proposition~\ref{simple-restriction-extension}(3) that $S(\calC)$ is a simple module for $\End_\BB(\calC)$.
The image on $S(\calC)$ of a biset in $\Iness^<(\calC)$ must be 0 because it factors through a category $\calD$ with $\calD<\calC$, and $S(\calD)=0$. Thus $\Iness^<(\calC)$ acts as 0 on $S(\calC)$, and the rest is immediate.
\end{proof}

It is the next result that gives a parametrization of the simple biset functors.

\begin{theorem}
\label{simple-parametrization}
With respect to the well order $<$, the isomorphism types of simple biset functors are parametrized by the pairs $(\calC,V)$, where $\calC$ is a category and $V$ is a simple $\Ess^<(\calC)$-module, in such a way that a simple biset functor $S$ corresponds to the pair where $\calC$ is the least category with $S(\calC)\ne 0$ and $V=S(\calC)$.
\end{theorem}

\begin{proof}
For every simple functor there is a least category $\calC$ with $S(\calC)\ne 0$. The endomorphisms in $\Iness^<(\calC)$ are all zero on $S$ and so  $V$ is a simple  $\Ess^<(\calC)$-module by Proposition~\ref{inessential-zero}. Because $S$ is determined by each of its values on categories by Proposition~\ref{simple-restriction-extension}, the assignment from isomorphism classes of simple functors to pairs $(\calC,V)$ is injective.

We must also show that for each simple $\Ess^<(\calC)$-module $V$ there is a simple biset functor $S$ for which $\calC$ is the least category with $S(\calC)\ne 0$ and $S(\calC)=V$. Now, there does exist a simple biset functor $S$ with $S(\calC)=V$, by  Proposition~\ref{simple-restriction-extension}. If $\calD$ is another category with $S(\calD)\ne 0$ then both $S(\calC)$ and $S(\calD)$ generate $S$, and so
$$
\begin{aligned}
S(\calD) &= \Hom_\BB(\calC,\calD)\cdot S(\calC)\quad\hbox{and}\\
S(\calC) &= \Hom_\BB(\calD,\calC) \cdot S(\calD)\\
\end{aligned}
$$
Thus $S(\calC) = \Hom_\BB(\calD,\calC) \Hom_\BB(\calC,\calD)\cdot S(\calC)$. We cannot have $\calD<\calC$, because then the homomorphism composites in the last equation would lie in $\Iness^<(\calC)$, which acts as 0 on $\calC$. We deduce that $\calC$ is the least category with $S(\calC)\ne 0$ and $S(\calC)=V$.
\end{proof}

\begin{example}
Let $<$ be any well order on finite groups that refines the partial order given by the size of the group. Thus if $G$ has fewer elements than $H$ we require $G<H$. Then 
$$
\End_\BB(G)\cong R\Out G \oplus \Iness^<(G)
$$
so that $\Ess^<(G)\cong R\Out G$. This is Proposition 4.3.2 of \cite{Bou3}.  It shows that we recover the usual parametrization of simple biset functors for groups by means of pairs $(H,V)$ where $H$ is a group and $V$ is a simple module for $\Out H$.
This includes the case of the simple functor $S_{1,R}$ when $R$ is a field, where $1$ is the identity group.
\end{example}

\begin{example}
\label{discrete-essential-algebra-example}
We have already seen in Example~\ref{discrete-categories-simple} that, for each integer $m\ge 1$, $S_{1,R}$ is the unique simple functor that is non-zero on the discrete category $[m]$, and in fact $S_{1,R}([m])\cong R^m$. We will now see that this accords with the parametrization of simple functors of Theorem~\ref{simple-parametrization} in terms of essential algebras.

\begin{proposition}
\label{discrete-essential-algebra-proposition}
Put the discrete categories $[m]$ in the order $[m] < [n]$ if and only if $m<n$.
\begin{enumerate}
\item All $([m],[m])$-bisets factor through $\one=[1]$.
\item $\Iness^<([m]) = \End_\BB([m])$ and 
$$
\Ess^<([m])=\begin{cases}
0&\hbox{if }m>1\cr
R&\hbox{if }m=1.\cr
\end{cases}
$$
\item With this choice of well order on the discrete categories, the only simple functor parametrized on these categories is $S_{1,R}$
\end{enumerate}
\end{proposition}

\begin{proof}
(1) The indecomposable $([m],[m])$-bisets are the bisets $E_{ij}$ that consist of a single point set in position $(i,j)$ and are empty elsewhere. Now $E_{ij}=E_{i1}E_{1j}$ is a factorization through $[1]$. 

(2) is immediate and (3) follows because if $m>1$ there are no simple $\Ess^<([m])$-modules to form part of the parametrization.
\end{proof}

Compare this result with Proposition~\ref{discrete-category-proposition2} where it is observed that $\Ess^<([m])$ is not isomorphic to $R\Out ([m]) = RS_m$ when $m>1$ and that the pairs  $(\calC,V)$ where $V$ is a simple representation of $\Out(\calC)$ do not parametrize the simple functors $S$. 
\end{example}

\subsection{Calculating the values of simple biset functors}
\label{calculation-subsection}

We now provide a general way to compute the value of a simple biset functor $S$ at a category $\calD$, given the value of $S$ at some other category $\calC$. It is related to methods described in \cite{Bou1}, \cite{BST} and \cite{Webc}. For this approach we construct a matrix as follows. Let $\Omega_j$, $j\in J$ be the indecomposable  $(\calD,\calC)$-bisets, and let $\Psi_i$, $i\in I$ be the indecomposable $(\calC,\calD)$-bisets, all taken up to isomorphism. Here $I$ and $J$ are just some (possibly infinite) indexing sets for these bisets. We form a block matrix $A$ with blocks indexed by $I\times J$, and where the block in position $(i,j)$ is the matrix of the endomorphism $S(\Psi_i\Omega_j): S(\calC)\to S(\calC)$, taken with respect to some chosen basis of $S(\calC)$.

\begin{theorem}
\label{simple-dimension-method}
Let $R$ be a field, let $\calC$ and $\calD$ be categories, let $S$ be a simple biset functor, and suppose that $S(\calC)$ is non-zero. Then
 $\dim S(\calD)$ is the column rank of the matrix $A$ just constructed.
 \end{theorem}
 
 \begin{proof}
We continue with the notation established just before the theorem.
Consider the composite map
$$
\bigoplus_{j\in J} S(\calC)\xrightarrow{(S(\Omega_j))} S(\calD)\xrightarrow{(S(\Psi_i))}  \prod_{i\in I} S(\calC).
$$
By  Theorem~\ref{simplicity-criterion} the map on the left is surjective and the map on the right is injective. The matrix of this composite map is $A$. It follows that the dimension of the image of this map equals $\dim S(\calD)$.
\end{proof}

We exemplify Theorem~\ref{simple-dimension-method} by calculating the values of $S_{\one,R}^{1,1}(\calP)$ when $\calP$ is a poset. Here $R$ is a field, $\one$ is the category that has a single object $*$ and only the identity morphism, and  $S_{\one,R}^{1,1}$ is the simple functor defined on $\BB_R^{1,1}$ that is non-zero on $\one$.
It is already known that if $G$ is a finite group then $S_{\one,R}^{1,1}(G)$ is $R$ when $|G|$ is invertible in $R$ and is 0 otherwise, according to a formula in \cite{Weba}. 

We regard the poset $\calP$ as a category in which the poset elements are the objects and there is a unique morphism $x\to y$ if and only if $x\le y$. This includes not only the discrete categories we studied in Example~\ref{discrete-category-bisets} and the posets $2^X$ that we shall study in the next section, but also representations of quivers. Given a quiver $Q$ without oriented cycles, its representations are the same thing as representations of the free category $FQ$ constructed from $Q$ (see \cite{MacLane}), and $FQ$ is a poset. Thus representations of quivers are also representations of posets. We recall (see \cite{Webb}) that, over any field $k$, a poset $\calP$ with $n$ elements has $n$ simple representations, so that the rank of the Grothendieck group $G_0(k\calP)$ is $n$.

\begin{theorem}
Let $R$ and $k$ be fields. Then $S_{\one,R}^{1,1}\cong R\otimes_\ZZ G_0(k(-))$ as biset functors on the full subcategory of $\BB^{1,1}$ whose objects are posets.
\end{theorem}

\begin{proof}
 We will first compute the dimension of $\dim S_{\one,R}^{1,1}(\calP)$ using the method of Theorem~\ref{simple-dimension-method}.  Place the elements of $\calP$ in a total order that refines the partial order on $\calP$.
The  $(\calP,\one)$-bisets are really just left $\calP$-sets, and they are birepresentable if and only if they are representable as $\calP$-sets, and similarly for $(\one,\calP)$-bisets. For $x\in \calP$ let $\Omega_x$ be the $(\one,\calP)$-biset generated at $(*,x)$, and let $\Psi_y$ be the $(\calP,\one)$-biset generated at $(y,*)$. Then the composition $\Omega_x\Psi_y\ne\emptyset$ if and only if  $y\le x$, by examining the definition of the composition. This means the matrix $A$ constructed for Theorem~\ref{simple-dimension-method} is triangular. Its diagonal entries are 1 because the diagonal entries are the functorial effect of sets 
$$
\Omega_x\Psi_x=\bigsqcup_{y\in\calP}(\Omega_x(y)\times \Psi_x(y))/\sim
$$ 
and there is only one non-empty term in the disjoint union, namely $\Omega_x(x)\times \Psi_x(x)$, which is a single point.
Thus $A$ has rank $|\calP|=\dim S_{\one,R}^{1,1}(\calP)$.

The biset functor $R\otimes_\ZZ G_0(kQ)$ is non-zero on $\one$, so it has $S_{\one,R}^{1,1}$ as a composition factor, this being the only simple functor that is non-zero on $\one$. Because  $R\otimes_\ZZ G_0(k\calP)$ and $S_{\one,R}^{1,1}(\calP)$ have the same dimension they are isomorphic, and this is an isomorphism of biset functors.
\end{proof}

\begin{corollary}
Let $R$ and $k$ be fields and let $\calP$ be a finite poset. Then $R\otimes_\ZZ G_0(k(\calP))$ is a simple module for $\End_{\BB^{1,1}}(\calP)$. 
\end{corollary}

It is harder to determine the values of $S_{\one,R}^{\calX,\calY}(\calC)$ when $\calX\ne 1\ne \calY$. For example, we may calculate that $\dim_R S_{\one,R}^{\all,\all}(\calA_2)=3$ when $\calA_2$ is the  category 
$$\calA_2 =\mathop{\bullet}\limits_x\xrightarrow{\alpha}\mathop{\bullet}\limits_y$$
of Example~\ref{a2-sets}. The value $S_{\one,R}^{\all,\all}(\calA_2)$ identifies as the representation ring $K_0(R\calA_2,\oplus)$ for this particular category. We omit these calculations.

Finally, we note that calculations of the values of $S_{\one,R}^{\calX,\calY}(\calC)$ are given in \cite{Bou1}, \cite{Weba} and elsewhere when $\calC$ is a finite group and $\calX,\calY\in\{1,\all\}$.

\subsection{Finiteness and projective covers}

We continue with the understanding that  $\BB$ denotes any of the categories $\BB_R^{\calX,\calY}$ with $\calX,\calY\in\{1,\all\}$ and $\calF$ the category of biset functors defined on $\BB$. We will write $F_\calC:=\Hom_\BB(\calC,-)$ for the representable biset functor at $\calC$. 

\begin{proposition}
\begin{enumerate}
    \item $F_\calC$ is a projective biset functor. 
    \item The category $\calF$ has enough projectives.
    \item $\End_\calF(F_\calC)\cong \End_\BB(\calC)$.
\end{enumerate}
\end{proposition}

\begin{proof}
All this is a standard consequence of Yoneda's Lemma, see \cite{Webb}.
\end{proof}

 We now restrict $R$ to be a field or (more generally) a complete discrete valuation ring.
The endomorphism ring $\End_\BB(\calC)$ in the last proposition may or may not have finite rank. As an example of this we have seen in Example~\ref{Kronecker-biset-example} that the endomorphism ring of the Kronecker quiver in $\BB^{1,1}$ does not have finite rank.  When it has infinite rank this may affect questions such as the existence of projective covers and whether or not the Krull-Schmidt theorem holds. We avoid this problem and study only the case when it does have finite rank, which is equivalent to saying that there are only finitely many $(\calC,\calC)$-bisets up to isomorphism.

\begin{theorem}
\label{projective-cover-theorem}
Let $\calC$ be a finite category for which there are only finitely many isomorphism classes of $(\calC,\calC)$-bisets and let $S$ be a simple biset functor. Then:
\begin{enumerate}
\item  $\dim S(\calC)$ is finite.
\item There are only finitely many simple biset functors $S$ (up to isomorphism) for which $S(\calC)\ne 0$.
\item If $S$ is a simple biset functor with $S(\calC)\ne 0$ then $S$ has a projective cover $P_S$.
\item  The representable biset functor $F_\calC = \Hom_\BB(\calC,-)$ decomposes as a finite direct sum
$$
F_\calC \cong \bigoplus_{S} P_S^{\dim S(\calC)/\dim \End_\BB(S)}.
$$
\end{enumerate}
\end{theorem}

\begin{proof}
(1) We start by noting that $\End_{\BB\mod}(F_\calC)\cong \End_\BB(\calC)^\op$ has finite rank over $R$, and also that $\Hom(F_\calC,S)\cong S(\calC)$, both by Yoneda's lemma. In particular, $S(\calC)\ne 0$ if and only if there is a non-zero homomorphism $F_\calC\to S$, and such a homomorphism gives a map $\End(F_\calC)\to \Hom(F_\calC,S)$, which is surjective by projectivity of $F_\calC$. This implies that $\dim S(\calC) = \dim \Hom(F_\calC,S)\le \Rank\End(F_\calC)$, which is finite.

(2) This is immediate from part (4) of Corollary~\ref{simple-evaluations}.

(3) The argument here is the same as in \cite[Sect. 5]{Weba}. Let $\Rad F_\calC$ be the interesection of the maximal subfunctors of $F_\calC$, or in other words the intersection of the kernels of the homomorphisms from $F_\calC$ to simple functors. There are only finitely many simple functors that can appear, so $F_\calC/\Rad F_\calC$ is semisimple -- a finite direct sum of simple functors. The simple functor $S$ is one of these, so there is a primitive idempotent $e\in\End(F_\calC/\Rad F_\calC)$ so that $e(F_\calC/\Rad F_\calC) = S$. By projectivity, the map $\End(F_\calC)\to \End(F_\calC/\Rad F_\calC)$ is surjective, and these are finite rank $R$-orders, so $e$ lifts to a primitive idempotent $\hat e\in\End(F_\calC)$. It's image $P_S:= \hat e\cdot F_\calC$ is now an indecomposable projective summand with the property that $P_S/\Rad(P_S)\cong S$, and so this is a projective cover of $S$.

(4) Being semisimple, $F_\calC/\Rad F_\calC\cong \bigoplus_{S} S^{d_S}$ for certain integers $d_S$, corresponding to a decomposition $F_\calC\cong  \bigoplus_{S} P_S^{d_S}$. We have
$$
\dim\Hom(F_\calC,S) = d_S \dim\End_\BB(S),
$$
and from this the result follows.
\end{proof}

We conclude this section with a result that is an extension to categories of a well known result for groups.

\begin{corollary}
Let $R$ be a field or a complete discrete valuation ring. The Burnside ring functor $B=F_\one:=\Hom_\BB(\one,-)$ is indecomposable and projective. It is the projective cover of the simple functor $S_{1,\bar R}$, where $\bar R$ is the residue field of $R$.
\end{corollary}

\begin{proof}
We have seen in Example~\ref{Burnside-functor-example} that the Burnside ring functor identifies as the representable functor $F_\one$. The result follows from Theorem~\ref{projective-cover-theorem}(4), because, by Corollary~\ref{simple-evaluations}, $S_{1,\bar R}$ is the only simple functor that is non-zero on $\one$, and its value there is $\bar R$.
\end{proof}

\section{Correspondences and bisets}

In this section we describe a way that the theory of correspondences of Bouc and Th\'evenaz \cite{BT} fits in with the theory of biset functors on categories. The goal is to show in Corollary~\ref{simple-correspondence-corollary} that simple correspondence functors and simple biset functors on Boolean lattices, using only birepresentable bisets, are parametrized the same way.

 Given sets $X$ and $Y$, a \textit{correspondence} from $X$ to $Y$ is a subset of the cartesian product $Y\times X$, the order being reversed to facilitate the composition of correspondences, about to be defined. Correspondences are also called \textit{relations}. Given correspondences $U\subseteq Z\times Y$ and $V\subseteq Y\times X$ another correspondence $UV$ is defined by
$$
UV:=\{ (z,x)\in Z\times X\bigm| \exists\;  y\in Y\hbox{ so that }(z,y)\in U\hbox{ and } (y,x)\in V\}.
$$
This composition of correspondences is associative and for each set $X$ there is an identity correspondence $\Delta_X = \{ (x,x)\bigm| x\in X\}\subseteq X\times X$. We form a category $\Corresp$ whose objects are finite sets and where $\Hom_\Corresp(X,Y)$ is the set of correspondences from $X$ to $Y$. Given a commutative ring $R$, a \textit{correspondence functor} over $R$ is a representation over $R$ of $\Corresp$, namely, a functor $\Corresp\to R\hbox{-mod}$. Equivalently we may form the linearization $R\Corresp$ of $\Corresp$, which has the same objects and whose morphisms are formal $R$-linear combinations of the morphisms in $\Corresp$ between two objects. A correspondence functor is now an $R$-linear functor $R\Corresp\to R\hbox{-mod}$. See \cite{BT} for details.

For each set $X$ we let $2^X$ denote the poset of subsets of $X$. For each correspondence $U\subseteq Y\times X$ we get a mapping ${}^+U: 2^X\to 2^Y$, and also a mapping $U^+:2^Y\to 2^X$, defined as follows: for each subset $A\subseteq X$ and $B\subseteq Y$ we put
$$
{}^+U(A):= \{ y\bigm| \exists\; (y,x)\in U,\; x\in A\}.
$$
and
$$
U^+(B):= \{ x\bigm| \exists\; (y,x)\in U,\; y\in B\}.
$$

In the next result we identify posets as categories in which the poset elements are the objects and there is a unique morphism $x\to y$ if and only if $x\le y$.

\begin{proposition}
\label{faithful-functor}
There is a functor $\Corresp\to \Cat$ specified on objects by $X\to 2^X$ and on morphisms by $U\to {}^+U$. There is also a functor  $\Corresp^\op\to \Cat$ specified on objects by $X\to 2^X$ and on morphisms by $U\to U^+$.
Both these functors are faithful.
In particular, for each set $X$, the monoid of correspondences on $X$ embeds in the monoid of order-preserving maps $2^X\to 2^X$.
\end{proposition} 

\begin{proof}
Given a correspondence $U\subseteq Y\times X$ it is evident that both ${}^+U$ and $U^+$ are order preserving maps. We also check that if $V\subseteq X\times W$ then ${}^+(UV)={}^+U{}^+V$ and $(UV)^+=V^+ U^+$. For the first equation, for each subset $A$ of $W$, the check is
$$
\begin{aligned}
{}^+(UV)(A)&={}^+\{(y,w)\bigm| \exists\; x\in X,\; (y,x)\in U, \; (x,w)\in V\} (A)\cr
&=\{ y\bigm| \exists\; x\in X,\;\exists\; w\in A,\; (y,x)\in U,\; (x,w)\in V\}\cr
&={}^+U\{x\bigm|\exists \;w\in A, \; (x,w)\in V\}\cr
&={}^+U({}^+V(A)).\cr
\end{aligned}
$$
The second equation is checked similarly. Each identity correspondence is sent to an identity functor.

These functors are one-to-one on objects. We check that they are one-to-one on morphisms. Suppose that $U,U'\subseteq Y\times X$ and that ${}^+U={}^+U'$. Then for each singleton $\{x\}\subseteq X$ we have ${}^+U \{x\}= {}^+U'\{x\}$, which shows that for each $x$, the pairs $(y,x)$ in $U$ are the same as for $U'$. It follows that $U=U'$.
\end{proof}

Not all order preserving maps $2^X\to 2^X$ arise in the form ${}^+U$ or $U^+$. The ones that do are exactly those that preserve joins of elements.

We have already examined in Proposition~\ref{cat-embedding} two functors $\phi:\Cat\to\BB$ and $\hat\phi:\Cat^\op\to\BB$. Composing the two functors that are covariant, and also composing the two functors that are contravariant, we obtain two covariant functors $\Corresp\to\BB$, both of which take a set $X$ to the poset $2^X$. The first of these functors takes a correspondence $U\subseteq Y\times X$ to the biset ${}_{2^{Y}}2^Y_{{}^{({}^+U)}2^X}$, and the second takes $U$ to the biset ${}_{2^{Y(U^+)}}2^X_{2^X}$.  As a consequence of this we note that any biset functor $\BB\to R\mod$ defines by composition a correspondence functor, in two ways. We will examine the nature of this relationship further.

\begin{example}
\label{bisets-from-correspondences-example}
Let $X=Y=[2]:=\{1,2\}$ and let $U\subset Y\times X$ be the correspondence $U:=\{(1,1),(1,2),(2,2)\}$. We obtain two bisets that we describe by matrices in the manner of Example~\ref{discrete-category-bisets}:
$$
 {}_{2^{[2]}}2^{[2]}_{{}^{({}^+U)}2^{[2]}} =
 \begin{bmatrix}
*&\emptyset&\emptyset&\emptyset\cr
*&*&\emptyset&\emptyset\cr
*&\emptyset&\emptyset&\emptyset\cr
*&*&*&*\cr
\end{bmatrix},
\qquad
{}_{2^{{[2]}(U^+)}}2^{[2]}_{2^{[2]}} =
\begin{bmatrix}
*&\emptyset&\emptyset&\emptyset\cr
*&*&*&*\cr
*&\emptyset&*&\emptyset\cr
*&*&*&*\cr
\end{bmatrix}.
$$
The rows and columns of these matrices are both indexed by the subsets
$$(\emptyset,\{1\}, \{2\}, \{1,2\})$$
 of $[2]$, taken in that order. Thus, for example, the third column of the second matrix is indexed by $\{2\}$ and indicates a $2^Y$-set that has its support on $\{1\},\{2\}$ and $\{1,2\}$. It can be described diagrammatically as the set
 $$
 \diagram{\emptyset&\umapright{}&*\cr \rmapdown{}&&\rmapdown{}\cr *&\umapright{}&*\cr} \quad\hbox{for the category}\quad 2^Y=  \diagram{\emptyset&\umapright{}&\{1\}\cr \rmapdown{} &&\rmapdown{}\cr \{2\}&\umapright{}&\{1,2\}\cr}.
 $$
 It is calculated by putting a $*$ in the position $(B,\{2\})$ whenever $\{2\}\subseteq U^+(B)$, for each subset $B\subseteq Y$.
\end{example}

We have seen that the functors $\phi$ and $\hat\phi$ are not faithful, in general. They are, however, if we restrict to the posets $2^X$.

\begin{proposition}
\label{correspondence-embedding}
The functors $\phi$ and $\hat\phi$, when restricted to the full subcategory of $\BB$ whose objects are the posets $2^X$, are faithful. Consequently the two functors $\Corresp\to \BB$ just described are faithful.
\end{proposition}

\begin{proof}
The functors $\phi:\Cat\to\BB$ and $\hat\phi:\Cat^\op\to\BB$ are the identity on objects. According to Proposition~\ref{cat-embedding}, two functors (i.e. order-preserving maps) $F, F':2^X\to 2^Y$ are sent to isomorphic bisets if and only if they are naturally isomorphic. This means there is a family of isomorphisms $\eta_x:F(x)\to F'(x)$ in $2^Y$ so that certain diagrams commute. The only isomorphisms in $2^Y$ are identity maps, so we deduce that $F=F'$ in this case. Thus $\phi$ and $\hat\phi$ are faithful. Also the two functors $\Corresp\to\Cat$ are faithful by Proposition~\ref{faithful-functor}, so the composites are faithful as well.
\end{proof}

The next lemmas will be applied to bisets for the posets $2^X$.

\begin{lemma}
\label{terminal-objects}
Let $\calC$ be a category with a terminal object $t$. Then a $\calC$-set $\Omega$ is $\sqcup$-indecomposable if and only if $\Omega(t)$ is a single point.
\end{lemma}

\begin{proof}
(1) Every non-empty $\calC$-set $\Omega$ has $\Omega(t)\ne\emptyset$, so if $\Omega=\Omega_1\sqcup\Omega_2$ is a disjoint union of non-empty $\calC$-sets  then $\Omega(t)$ has size at least 2. Thus if $\Omega(t)$ is a single point, it is indecomposable. Conversely, if $\Omega(t)=\{u_1,\ldots,u_r\}$ consists of $r$ points we may define $\Omega_i(x):= f_x^{-1}(u_i)$ where $f_x:x\to t$ is the unique morphism to the terminal object, and now $\Omega=\Omega_1\sqcup\cdots\sqcup\Omega_r$ decomposes.

\end{proof}

\begin{definition}
For any $(\calC,\calD)$-biset $\Omega$ we will say that a \textit{row} of $\Omega$ is the set of values $\Omega(u,v)$ where $u$ is fixed and $v$ is allowed to vary; and a \textit{column}  of $\Omega$ is the set of values $\Omega(u,v)$ where $v$ is fixed and $u$ is allowed to vary. A row has the structure of a $\calD$-set and a column has the structure of a $\calC$-set. This is consistent with the matrix notation of Example~\ref{bisets-from-correspondences-example}.
\end{definition}

\begin{lemma}
An indecomposable $(2^Y,2^X)$-biset $\Omega$ is birepresentable if and only if 
\begin{enumerate}
\item for each row of $\Omega$, the support of that row consists of the subsets of some particular subset of $X$;
\item for each column of $\Omega$, the support of that column  consists of the supersets of some particular subset of $Y$; and 
\item for each pair of subsets $A\subseteq X$ and  $B\subseteq Y$ the value $\Omega(B,A)$ has size 0 or 1.
\end{enumerate}
\end{lemma}

\begin{proof}
If each row and column have the given form, it is representable, represented by the maximal subset of the support in the case of rows, and the minimal subset of the support in the case of columns. 

Conversely, suppose that $\Omega$ is birepresentable. Then each row and column is representable, and so is a disjoint union of representable $(2^X)^\op$-sets (for rows) or  $2^Y$-sets (for columns). We claim that each row and column is indecomposable. To see this we observe that  $2^Y\times (2^X)^\op$ has a terminal element $(Y,\emptyset)$, and $\Omega$ is an indecomposable set for this category, so $\Omega(Y,\emptyset)$ is  a single element by Lemma~\ref{terminal-objects}.  It follows that the row of $\Omega$ determined by $Y$ is indecomposable (again  by Lemma~\ref{terminal-objects}), as is the column determined by $\emptyset$, and hence this row and column only have entries of size $0$ or $1$, because hom sets of  $2^X$ and $2^Y$ have size $0$ or $1$. From this it follows that every row and column of $\Omega$ is indecomposable, because they are sets for $2^X$ or $2^Y$, which have terminal objects, and the values at these objects have size $0$ or $1$. This means each row and column is representably generated by a single element, and each evaluation of $\Omega$ consists of a single element or is empty.
\end{proof}

\begin{lemma}
\label{birepresentable-product-indecomposable}
Let $\Omega$ be an indecomposable birepresentable $(2^Y,2^X)$-biset and $\Psi$  an indecomposable birepresentable $(2^Z,2^Y)$-biset. Then $\Psi\circ\Omega$ is either an indecomposable $(2^Z,2^X)$-biset or is empty.
\end{lemma}

\begin{proof}
By Lemma~\ref{terminal-objects} it suffices to show that $(\Psi\circ\Omega)(Z,\emptyset)$ consists of a single point, or is empty. From the definition of biset composition this set is
$$
\bigsqcup_{U\subseteq Y}\Psi(Z,U)\times\Omega(U,\emptyset)/\sim
$$
Because $\Psi$ and $\Omega$ are birepresentable there are subsets $A,B$ of $Y$ so that $\Psi(Z,U)\ne\emptyset$ if and only if $U\subseteq B$ and $\Omega(U,\emptyset)\ne\emptyset$ if and only if $U\supseteq A$, so that the disjoint union is non-empty if and only if $A\subseteq U\subseteq B$, in which case $\Psi(Z,U)\times\Omega(U,\emptyset)= \{(b(i_U^{B})^*,(i_A^U)_*a)\}$, where $a\in\Omega(A,\emptyset)$ and $b\in\Psi(Z,B)$ are generators of the first column of $\Omega$ and the last row of $\Psi$ respectively, and $i_A^U, i_U^B$ are the inclusion morphisms in $2^Y$. All of these pairs are equivalent under $\sim$ because $(b(i_U^{B})^*,(i_A^U)_*a)\sim (b,(i_A^B)_*a)$, for instance, independently of $U$. This shows that  $(\Psi\circ\Omega)(Z,\emptyset)$ consists of a single point.
\end{proof}

\begin{proposition}
\label{correspondence-bisets}
Let $U\subseteq Y\times X$ be a correspondence and let $\Omega = {}_{2^{Y}}2^Y_{{}^{({}^+U)}2^X}$ as a $(2^Y,2^X)$-biset. Then
\begin{enumerate} 
\item $\Omega$ is birepresentable,
\item the row of entries $\Omega(Y,A)$ as $A$ ranges through subsets of $X$, and the column of entries $\Omega(B,\emptyset)$ as $B$ ranges through subsets of $Y$, are all non-empty with a single point $*$ in every entry.
\end{enumerate}
\end{proposition}

\begin{proof}
1. As a right $2^X$-set $\Omega$ is 
$$
\bigsqcup_{B\in 2^Y} {}_{2^{Y}}2^Y_{{}^{({}^+U)}2^X}(B,-) = \bigsqcup_{B\in 2^Y} \Hom_{2^Y}({}^+U(-),B)
\cong  \bigsqcup_{B\in 2^Y} \Hom_{2^X}(-,\hat B)
$$
where $\hat B = \{ x\in X\bigm| {}^+U(\{x\})\subseteq B\}$, and this is representable by the subsets $\hat B$ of $X$. This is because there is a homomorphism ${}^+U(A)\to B$ in $2^Y$ if and only if ${}^+U(A)\subseteq B$, if and only if
$$
A\subseteq \{ x\in X\bigm| {}^+U(x)\subseteq B\} = \hat B.
$$

As a left $2^Y$-set $\Omega$ is $\bigsqcup_{A\in 2^X}\Hom_{2^Y}({}^+U(A),-)$, which is also representable, by the subsets ${}^+U(A)$ of $Y$.

2. This is a matter of computing these bisets. We have $\Omega(Y,A) = \Hom({}^+U(A), Y)=\{*\}$ because ${}^+U(A)\subseteq Y$ regardless of what $A$ is, and $\Omega(B,\emptyset)=\Hom({}^+U(\emptyset),B) = \Hom(\emptyset,B) = \{*\}$ regardless of what $B$ is, because $\emptyset\subseteq B$ always.
\end{proof}

\begin{definition}
Let us call a birepresentable $(2^Y,2^X)$-biset $\Omega$ a \textit{correspondence biset} if it satisfies the condition that the values $\Omega(Y,A)$ as $A$ ranges through subsets of $X$, and $\Omega(B,\emptyset)$ as $B$ ranges through subsets of $Y$, are all a single point $*$.
\end{definition}

We have seen in Proposition~\ref{correspondence-bisets} that the bisets of the form $ {}_{2^{Y}}2^Y_{{}^{({}^+U)}2^X}$ are correspondence bisets, and we will see in Corollary~\ref{biset-bijection} that all correspondence bisets have this form, giving a bijection between correspondences and correspondence bisets. Note that correspondence bisets are necessarily indecomposable by Lemma~\ref{terminal-objects}, because on the terminal object $(Y,\emptyset)$ their value has size 1.

\begin{example}
We consider again the bisets in Example~\ref{bisets-from-correspondences-example} and we refer to the matrices shown there. We have
 $X=Y=[2]:=\{1,2\}$ and $U\subset Y\times X$ is the correspondence $U:=\{(1,1),(1,2),(2,2)\}$. We see directly from the matrix that describes the biset $ {}_{2^{[2]}}2^{[2]}_{{}^{({}^+U)}2^{[2]}}$ that it is birepresentable: each row is a $(2^X)^\op$-set that has value a single point $*$ on all the subsets of a particular subset of $X$, and each column is a $2^Y$-set that has value a single point $*$ on all the supersets of a particular subset of $Y$. 
 
 By contrast to Proposition~\ref{correspondence-bisets}, the biset ${}_{2^{Y(U^+)}}2^X_{2^X}$ is not birepresentable. 
 This biset is not representable as a left $2^Y$-set because the column indexed by $\{2\}$ (the third column) is an indecomposable $2^Y$-set that requires two generators but takes the value a single point on the terminal object.
 
 Because it is not birepresentable ${}_{2^{Y(U^+)}}2^X_{2^X}$ is not a correspondence biset. On the other hand, the extra condition beyond birepresentability to be a correspondence biset is that the left column and bottom row of the matrix of the biset matrix have entries $*$, so we confirm that $ {}_{2^{[2]}}2^{[2]}_{{}^{({}^+U)}2^{[2]}}$ is a correspondence biset.
 \end{example}

In the next proposition the initial sentence says that $\Omega$ is not a correspondence biset.

\begin{proposition}
\label{factoring}
Let $\Omega$ be  an indecomposable birepresentable $(2^Y,2^X)$-biset for which either there is a subset $A\subseteq X$ with $\Omega(Y,A)=\emptyset$ or there is a subset $B\subseteq Y$ with $\Omega(B,\emptyset)=\emptyset$. Then $\Omega$ factors through $2^{X'}$ in the first case for some proper subset $X'\subseteq X$, and $\Omega$ factors through $2^{Y'}$ in the second case, for some proper subset $Y'\subseteq Y$.
\end{proposition}

\begin{proof}
Suppose there is a subset $A\subseteq X$ with $\Omega(Y,A)=\emptyset$. The row $\Omega(Y,-)$ has a unique largest set $X'\subseteq X$ for which $\Omega(Y,X')\ne \emptyset$ and $X'\ne X$ because $X'\not\supseteq A$. We define the $(2^Y,2^{X'})$-biset $\Omega'$ by restriction: $\Omega'(V,U) = \Omega(V,U)$ for $U\subseteq X'$ and $V\subseteq Y$. Now $\Omega = \Omega'\circ {}_{2^{X'\phi}}2^X_{2^X}$ where $\phi:2^{X'}\to 2^X$ is inclusion, showing that $\Omega$ factors through $2^{X'}$. The second case is similar: the column $\Omega(-,X)$ has a unique minimal set $Y'\subseteq Y$ so that $\Omega(Y',X)\ne\emptyset$. Now $\Omega$ factors through $2^{Y-Y'}$ in the following way. We define a $(2^{Y-Y'},2^{X})$-biset $\Omega'$ as follows: $\Omega'(V,U):=\Omega(V\sqcup Y',U)$ for $U\subseteq X$ and $V\subseteq Y-Y'$. Now
$\Omega ={}_{2^{Y}}2^Y_{{}^\phi2^{Y-Y'}}\circ  \Omega'$ where $\phi:2^{Y-Y'}\to 2^Y$ is inclusion.
\end{proof} 


\begin{proposition}
\label{determined-bisets}
Let $\Omega$ be a correspondence biset. Then
\begin{enumerate}
\item $\Omega$ is completely determined by the values $\Omega(B,A)$ where $A\subseteq X$ is a subset of size 1, and $B\subseteq Y$ is a subset of size $|Y|-1$, and
\item $\Omega= {}_{2^{Y}}2^Y_{{}^{({}^+U)}2^X}$ where $U\subseteq Y\times X$ is
$$
U=\{(y,x)\bigm| \Omega(Y-\{y\},\{x\})=\emptyset \}.
$$
\end{enumerate}
\end{proposition}

\begin{proof}
(1) Each row $\Omega(B,-)$ of $\Omega$ has a unique maximal set $A\subseteq X$ on which it is non-empty, and it has size at least 1, so it is the union of the one-points sets in the same row on which $\Omega$ is non-empty. Hence it is determined by the value on one-point sets in that row. Similarly, each column $\Omega(-,A)$ of $\Omega$ has a unique minimal set $B\subseteq Y$ on which it is non-empty, and it has size at most $|Y|-1$, so it is the intersection of the subsets of size $|Y|-1$ in the same column on which $\Omega $ is non-empty. Hence it is determined by the value on sets of size $|Y|-1$ in that column. Putting this together, the $\Omega(B,A)$ are determined in the first instance when $|A|=1$ and $B$ is arbitrary, and when $A$ is arbitrary and $|B|=|Y|-1$. From this we see that $\Omega$ is determined on all pairs of subsets.

(2) Write $\Psi={}_{2^{Y}}2^Y_{{}^{({}^+U)}2^X}$. We show that $\Psi = \Omega$. Now 
$$
\begin{aligned}
\Psi(Y-\{y\},\{x\})=\emptyset &\Leftrightarrow {}^+U(\{x\})\not\subseteq Y-\{y\}\\ 
&\Leftrightarrow (y,x)\in U\\ 
&\Leftrightarrow \Omega(Y-\{y\},\{x\})=\emptyset.\\
\end{aligned}
$$ 
Thus $\Psi$ and $\Omega$ agree on pairs $(B,A)$ where $A\subseteq X$ is a subset of size 1, and $B\subseteq Y$ is a subset of size $|Y|-1$. By part (1), $\Psi=\Omega$.
\end{proof}

\begin{corollary}
\label{biset-bijection} 
The map from correspondences $U$ on $Y\times X$ to (isomorphism classes of) correspondence $(2^Y,2^X)$-bisets, given by $U\mapsto  {}_{2^{Y}}2^Y_{{}^{({}^+U)}2^X}$, is a bijection. Consequently the number of such correspondence bisets is $2^{|X|\cdot |Y|}$.
\end{corollary}

\begin{proof}
We already know from Proposition~\ref{faithful-functor} that this map is one-to-one, and we see this again from Proposition~\ref{determined-bisets} (2). We see from Proposition~\ref{determined-bisets} (1) that the map is surjective.
\end{proof}

\begin{lemma}
\label{correspondence-biset-factorization}
Let $\Omega$ be a correspondence $(2^Y,2^X)$-biset. Suppose that $\Omega=\Psi\circ\Theta$ where $\Psi$ is a birepresentable $(2^Y,2^Z)$-biset and $\Theta$ is a birepresentable $(2^Z,2^X)$-biset. Then there is a factorization $\Omega=\Psi'\circ\Theta'$ where $\Psi'$ is a correspondence $(2^Y,2^Z)$-biset  and $\Theta'$ is a correspondence $(2^Z,2^X)$-biset.
\end{lemma}

\begin{proof}
We may assume $\Omega$ and $\Psi$ are indecomposable.
Whenever $\Psi(B,\emptyset)=\emptyset$ we have $\Psi(B,A)=\emptyset$ for all $A\subseteq X$ by birepresentability, and this forces $\Omega(B,\emptyset)=\emptyset$, which is not the case. Thus $\Psi(B,\emptyset)$ is a single point for all $B\subseteq Y$, and by a similar argument $\Theta(Y,A)$  is a single point for all $A\subseteq X$.

The row $\Psi(Y,A)$ as $A\subseteq X$ appears only in the calculation of the row $\Omega(Y,A)$, whose entries are each a single point. It follows that the biset $\Psi'$, obtained by changing the entries $\Psi(Y,A)$ to a single point if necessary, can also appear in the product $\Omega=\Psi'\circ\Theta$ without changing the answer. Making a similar construction of $\Theta'$ that is the same as $\Theta$, except perhaps at entries $\Theta'(B,\emptyset)$ which may be changed to be a single point, we obtain the desired factorization $\Omega=\Psi'\circ\Theta'$.
\end{proof}

We now compare the essential algebra of correspondences on a set $X$ with the essential algebra of bisets for $2^X$. Throughout, we consider algebras over a commutative ring $R$ that does not figure in our notation, and our result does not depend on the choice of $R$. We gave in \ref{essential algebras} a definition of the essential algebra of a category $\calC$ in the biset category, and it depends on placing a well order on categories. We are now considering only categories $2^X$. We put $2^X < 2^Y$ if and only if $|X|<|Y|$, and we consider essential algebras in $\BB^{1,1}$. Thus the inessential ideal $\Iness_{2^X}^<$ is the 2-sided ideal in $\End_{\BB^{1,1}}(2^X)$ that is the $R$-span of the birepresentable bisets that factor birepresentably through categories $2^Z$ with $|Z|<|X|$. We put $\Ess_{\BB^{1,1}}^<(2^X):=\End_{\BB^{1,1}}(2^X)/\Iness_{2^X}^<$. This follows the similar definition by Bouc and Th\'evenaz \cite{BT} of the essential algebra $\Ess_\Corresp(X)$ of $X$ in the correspondence category, which  is the quotient of $\End_\Corresp(X)$ by the ideal of correspondences that factor through smaller sets than $X$.

We come to our main theorem in this section.

\begin{theorem}
\label{essential-algebra-isomorphism}
The functor $\Corresp\to\BB^{1,1}$ that sends a set $X$ to the poset $2^X$, and sends a correspondence $U\subseteq Y\times X$ to the biset ${}_{2^{Y}}2^Y_{{}^{({}^+U)}2^X}$, induces an isomorphism of essential algebras
$$
\Ess_{R\Corresp}(X)\to \Ess_{\BB^{1,1}}^<(2^X).
$$
\end{theorem}

\begin{proof}
By functoriality this functor gives a ring homomorphism
$$
\End_{R\Corresp}(X)\to \End_{\BB^{1,1}}(2^X).
$$
If a correspondence factors through a smaller set $Z$ than $X$, then the biset that is its image factors through the smaller category $2^Z$ and so lies in $\Iness_{2^X}^<$. For this formal reason we obtain a homomorphism $\Ess_{R\Corresp}(X)\to \Ess_{\BB^{1,1}}^<(2^X)$. We will show that it is an isomorphism.

By Proposition~\ref{factoring},  $\Ess_{\BB^{1,1}}^<(2^X)$ is spanned by the images of bisets $\Omega$ for which  $\Omega(X,A)\ne\emptyset$  and   $\Omega(A,\emptyset)\ne\emptyset$ for all $A\subseteq X$. Corollary~\ref{biset-bijection} shows that these bisets are all images of correspondences. This shows that the homomorphism is a surjection.

Finally, we show that the homomorphism is an injection. Suppose we have an element $\sum_i\lambda_i U_i\in \End_{R\Corresp}(X)$ that maps to an element of $\Iness_{2^X}^<$. We will show that it lies in the inessential ideal $\Iness_X$ of the correspondence algebra at $X$, consisting of elements that factor through smaller sets. Now  $\Iness_{2^X}^<$ is spanned by bisets $\Psi\circ\Omega$ where $\Psi$ is a $(2^X,2^Z)$-biset and $\Omega$ is a $(2^Z,2^X)$-biset, for some set $Z$ smaller than $X$, and these bisets are indecomposable and birepresentable. By Lemma~\ref{birepresentable-product-indecomposable} this product is either empty, or indecomposable, so that if a linear combination of bisets factors through a smaller set than $X$, so does each indecomposable biset in that linear combination. This implies that if a linear combination of correspondences maps to an element of $\Iness_{2^X}^<$, so does each correspondence in that linear combination.

We are now reduced to the situation of a correspondence $U\subseteq X\times X$ for which the associated biset ${}_{2^{X}}2^X_{{}^{({}^+U)}2^X}$ factors through some $2^Z$, where $Z$ is smaller than $X$, as a product of birepresentable bisets. By Lemma~\ref{correspondence-biset-factorization} we may assume that they are correspondence bisets ${}_{2^{X}}2^X_{{}^{({}^+U)}2^X}=\Psi\circ\Theta$. Each of $\Psi$ and $\Theta$ is the image of a correspondence $V$ and $W$. Because $U$ and $VW$ have the same bisets as their image we have $U=VW$ by Proposition~\ref{correspondence-embedding}, showing that  $U$ lies in the inessential ideal $\Iness_X$. This completes the proof that the homomorphism is an injection.
\end{proof}

\begin{corollary}
\label{simple-correspondence-corollary}
The simple biset functors defined on the full subcategory of $\BB^{1,1}$ whose objects are the posets $2^X$ are parametrized the same way as the simple functors on the correspondence category $\Corresp$. Specifically, if $S$ is a simple biset functor on $\BB^{1,1}$ with $S(2^X)\ne 0$ for some set $X$ of smallest size, then the correspondence functor $\hat S$ obtained from $S$ by composition with the inclusion functor $\Corresp \to \BB^{1,1}$  has a unique simple composition factor $T$ with $X$ minimal such that $T(X)\ne 0$. This mapping $S\mapsto T$ is a bijection between simple biset functors that do not vanish on the $2^X$ and simple correspondence functors.
\end{corollary}

\begin{proof}
We take the well order given by $2^X < 2^Y$ if and only if $|X| < |Y|$ and now Theorem~\ref{simple-parametrization} implies that the simple biset functors $S$ that do not vanish on posets $2^X$ are parametrized by pairs $(2^X,V)$ where $V$ is a simple module for $\Ess_{\BB^{1,1}}^<(2^X)$, in such a way that $S(2^X) = V$ and $X$ is minimal such that $S(2^X)\ne 0$. Restricting to correspondence functors, $V$ remains simple for $\Ess_{R\Corresp}(X)$ because of the isomorphism in Theorem~\ref{essential-algebra-isomorphism}, and so the pair $(X,V)$ parametrizes a simple correspondence functor $T$ with $T(X)=V$ and $T(Z)=0$ if $|Z|<|X|$; and $T$ is characterized by this evaluation, by the theory of \cite{BT}. As a correspondence functor, $S$ must have a composition factor whose evaluation at $X$ is $V$, and so $T$ is this composition factor. All simple correspondence functors arise in this way, and we obtain a bijection between these simple biset functors that do not vanish on the $2^X$ and simple correspondence functors.
\end{proof}

We refer to \cite{BT} for the specifics of the parametrization of simple correspondence functors.

\section{Further structures on biset functors}

There are tensor product and internal Hom structures on the biset category, and also on the category of biset functors, extending the corresponding structures for biset functors on groups. On groups, the definitions and properties are described  in Chapter 8 of \cite{Bou3}, and our development here proceeds in parallel with that. Many arguments from that exposition are similar to those that appear here, and sometimes they apply verbatim. Sometimes inverses of group elements arise, but in most, if not all, cases this is because an isomorphism to the opposite group is being  invoked implicitly, and a valid argument can be found without inverses, using opposite categories. In group theory it is unusual to consider the opposite of a group, it being isomorphic to the group itself, but in our situation the arguments at times gain transparency because we are forced to use opposite categories.

In this section we will only consider the biset category $\BB=\BB_R$ where \textit{all} bisets are used. In some places the arguments we will present do not work for bisets that are representable on the left or on the right, and so we exclude these from consideration.
We write $\calF=\calF_R$ for the category of biset functors $\BB\to R\mod$, usually suppressing the notation for the commutative ring $R$.

We start by observing that $\BB$ is a rigid symmetric monoidal additve $R$-linear category, with $R$-bilinear product. It also possesses a contravariant self-equivalence and an internal Hom. From this we go on to deduce that $\calF$ acquires the structure of a symmetric monoidal category and an internal Hom. From that, we define Green functors and Green modules on categories. We conclude by showing how fibered biset functors may be made to work on categories. For this we consider a theory using bi-objects in categories other than $\Set$. All these constructions extend what is done for biset functors on groups.

\subsection{Structures on the biset category}

In the following we note that the term \textit{tensor category} is used in more than one way in the literature. 

\begin{proposition}
The biset category $\BB$ is a tensor category over $R$, in the sense that it is an $R$-linear additive category with a symmetric monoidal structure $\otimes$ that is $R$-bilinear.
\end{proposition}

\begin{proof}
The tensor product operation $\otimes$ on objects of $\BB$ is the product of categories, and the unit is the category $\one$. We note that there are canonical isomorphisms in $\BB$: 
$$
\begin{aligned}
(\calC \times\calD)\times\calE &\cong \calC \times (\calD\times\calE)\\
\calC \times\calD&\cong \calD \times\calC\\
\one\times\calC\cong&\calC\cong\calC\times\one\\
\end{aligned}
$$
showing that the identities for a symmetric monoidal structure are satisfied. We observe that in $\BB$ the coproduct and product of objects $\calC,\calD$ is their disjoint union $\calC\sqcup\calD$, because every $(\calC\sqcup\calD,\calE)$-biset is uniquely the disjoint union of a $(\calC,\calE)$-biset and a $(\calD,\calE)$-biset, and similarly with $(\calE,\calC\sqcup\calD)$-bisets. Thus $\BB$ is additive.
\end{proof}

For any object $X$ in a tensor category, a \textit{dual} of $X$ is an object $X^\vee$ equipped with \textit{coevaluation} and \textit{evaluation} maps
$$
\alpha: \one\to X\otimes X^\vee,\quad \beta: X^\vee\otimes X\to\one,
$$
so that the compositions
$$
\begin{aligned}
X\xrightarrow{\alpha\otimes\mathrm{id}} &X\otimes X^\vee\otimes X \xrightarrow{\mathrm{id}\otimes\beta} X\\
X^\vee\xrightarrow{\mathrm{id}\otimes\alpha} &X^\vee\otimes X\otimes X^\vee \xrightarrow{\beta\otimes\mathrm{id}} X^\vee\\
\end{aligned}
$$
are the identity. We say the tensor category is \textit{rigid} if every object has a dual.

We are about to show that objects in $\BB$ have duals.
The coevaluation and evaluation maps in $\BB$ are constructed in the following way.
For each category $X=\calC$ the identity $(\calC,\calC)$-biset ${}_\calC\calC_\calC$ may also be regarded as a $(\calC\times\calC^\op,\one)$-biset
$$
\alpha={}_{\calC\times\calC^\op}\calC_\one : \one\to \calC\times\calC^\op
$$
and also as a $(\one,\calC^\op\times\calC)$-biset
$$
\beta={}_\one\calC_{\calC^\op\times\calC}:\calC^\op\times\calC\to\one.
$$
All three are $\calC$ regarded as a set for $\calC\times\calC^\op$, after identifying $\calC\times\calC^\op$ with $\one\times(\calC^\op\times\calC)^\op$ and with $(\calC\times\calC^\op)\times\one^\op$. For example, if $x,y$ are objects of $\calC$ the value of ${}_{\calC\times\calC^\op}\calC_\one$ at $((y,x),1)$ is $\Hom_\calC(x,y)$ and this is also the value of ${}_\calC\calC_\calC$ at $(y,x)$ and ${}_\one\calC_{\calC^\op\times\calC}$ at $(1,(y,x))$. 

\begin{theorem}
\label{biset-dual-theorem}
Each object $\calC$ of $\BB$ has a dual in the sense of tensor categories, namely the opposite category $\calC^\op$. Thus the tensor category $\BB$ is rigid.
\end{theorem}

\begin{proof}
We show that the bisets $\alpha$ and $\beta$ defined above satisfy the identities to show that $\calC^\op$ is the dual of $\calC$. To do this we define bisets
$$
\begin{aligned}
{}_\calC\Omega_{\calC\times\calC^\op\times\calC} &= 
{}_\calC\calC_\calC\times {}_\one\calC_{\calC^\op\times\calC}\\
{}_{\calC\times\calC^\op\times\calC}\Psi_\calC
&={}_{\calC\times\calC^\op}\calC_\one\times {}_\calC\calC_\calC,\\
\end{aligned}
$$
so that if $r,s,x,z$ are objects of $\calC$ and $y$ is an object of $\calC^\op$ we have 
$$
\Omega(r,(x,y,z)) = \{(\alpha,\beta)\bigm| \alpha:x\to r,\;\beta:z\to y\}
$$
and
$$
\Psi((x,y,z),s)=\{(\gamma,\delta)\bigm| \gamma:y\to x,\;\delta:s\to z\}.
$$

We verify that 
    $$
    \Omega\circ\Psi \cong {}_\calC\calC_\calC
    $$
this being one of the identities we need to check.    The other identity is obtained by changing the roles of $\calC$ and $\calC^\op$ to obtain similar bisets $\Omega'$ and $\Psi'$ with $\Omega'\circ\Psi' \cong {}_{\calC^\op}\calC^\op_{\calC^\op}$ and it proceeds similarly.

    The definition of the biset product $\Omega\circ\Psi$ at a pair of objects $(r,s)$ has the form
    $$
    (\Omega\circ\Psi)(r,s) =\bigsqcup_{(x,y,z) \in \calC\times\calC^\op\times\calC} \Omega(r,(x,y,z))\times\Psi((x,y,z),s) \big/ \sim
    $$
    We define a morphism of bisets
    $$\Omega\circ\Psi(r,s)\to {}_\calC\calC_\calC(r,s)
    $$ by $((\alpha,\beta),(\gamma,\delta))\mapsto \alpha\gamma\beta\delta$. This specification is constant on equivalence classes of the relation $\sim$ because if $(\zeta,\theta,\psi)$ is a triple of morphisms in $\calC\times\calC^\op\times\calC$, with suitable domains and codomains, then the pairs
    $$
    ((\alpha,\beta)(\zeta,\theta,\psi),(\gamma,\delta)) = ((\alpha\zeta,\theta\beta\psi),(\gamma,\delta))
    $$ and
    $$
    ((\alpha,\beta),(\zeta,\theta,\psi)(\gamma,\delta))= ((\alpha,\beta),(\zeta\gamma\theta,\psi\delta))
    $$
    both map to the same morphism in ${}_\calC\calC_\calC$, namely $\alpha\zeta\gamma\theta\beta\psi\delta$. This means we have defined a morphism of bisets. In the opposite direction we have an inverse map $\Omega\circ\Psi \gets{}_\calC\calC_\calC$ given at $(r,s)$ by $\nu\mapsto ((1_r,1_r),(1_r,\nu))$. These two maps compose in both directions to give the identity, so we have verified the desired identity.

\end{proof}

A $(\calC,\calD)$-biset ${}_\calC\Omega_\calD$ is  the same thing as a $\calC\times\calD^\op$-set, which we can also regard as a $\calD^\op\times (\calC^{\op})^{\op}$-set in a canonical way. This gives ${}_\calC\Omega_\calD$ the structure of a $(\calD^\op,\calC^\op)$-biset ${}_{\calD^\op}\Omega_{\calC^\op}$. 

\begin{definition}
\label{duality-definition}
There is a contravariant functor $\tau:\BB^\op\to\BB$ given by $\tau(\calC) = \calC^\op$ on objects, and on morphisms $\tau({}_\calC\Omega_\calD) = {}_{\calD^\op}\Omega_{\calC^\op}$. 
\end{definition}

Evidently $\tau$ is a contravariant equivalence that is own inverse: a duality on $\BB$. It will be an important ingredient in defining the internal tensor product of biset functors.

\begin{proposition}
There is an internal Hom construction on $\BB$, given by $\calH(\calC,\calE)=\calE\times\tau(\calC)$. It satisfies $\Hom_\BB(\calC\times\calD,\calE)\cong \Hom_\BB(\calD,\calH(\calC,\calE))$, and also $\calH(\calC\times\calD,\calE)\cong \calH(\calD,\calH(\calC,\calE))$.
\end{proposition}

\begin{proof}
Both sides of the Hom isomorphism yield the Grothendieck group of sets for $\calE\times\calC^\op\times\calD^\op$. The two sides of the $\calH$ isomorphism are $\calE\times\tau(\calC\times\calD)$ and $(\calE\times\tau\calC)\times\tau\calD$, and these are isomorphic.
\end{proof}

Thus, for example, $\tau \cong \calH(-,\one)$.

\subsection{Structures on biset functors: the Yoneda-Dress construction}

We define two of these constructions, depending on whether we do the construction with respect to a category, or to its opposite.

\begin{definition}
For any category $\calK$ let $\ydp_\calK,\ydp_\calK^\op : \BB\to\BB$ be the functors defined as $\ydp_\calK(\calC) :=\calC\times\calK$ and $\ydp_\calK^\op(\calC) :=\calC\times\calK^\op$ on categories $\calC$. On morphisms of $\BB$ that are bisets ${}_\calC\Omega_\calD$ the definitions are 
$$
\begin{aligned}
\ydp_\calK({}_\calC\Omega_\calD) = {}_\calC\Omega_\calD\times {}_\calK\calK_\calK &: \calD\times\calK \to \calC\times\calK\cr
\ydp_\calK^\op({}_\calC\Omega_\calD) = {}_\calC\Omega_\calD\times {}_{\calK^\op}\calK^\op_{\calK^\op} &: \calD\times\calK^\op \to \calC\times\calK^\op\cr
\end{aligned}
$$
Here ${}_\calK\calK_\calK$ and ${}_{\calK^\op}\calK^\op_{\calK^\op}$ are the identity bisets at $\calK$ and $\calK^\op$. Note that $\ydp_\calK^\op$ is the same as $\ydp_{\tau(\calK)}$ and could also be written $\ydp_{\calK^\op}$, but this makes it harder to notate the functorial dependence on $\calK$.

We define $\ydP_\calK, \ydP_\calK^\op :\calF\to\calF$ to be the functors that are precomposition with $\ydp_\calK$ and $\ydp_\calK^\op$. 
Thus if $M$ is a biset functor and $\calC$ a category then 
$$\begin{aligned}
\ydP_\calK (M)(\calC) &= (M\circ \ydp_\calK) (\calC) = M(\calC\times\calK),\cr
\ydP_\calK^\op (M)(\calC) &= (M\circ \ydp_\calK^\op) (\calC) = M(\calC\times\calK^\op).\cr
\end{aligned}
$$
These are the \textit{Yoneda-Dress construction} of $M$ at $\calK$ and the \textit{opposite Yoneda-Dress construction} of $M$ at $\calK$.  In other contexts the Yoneda-Dress construction is also written $M_\calK$ instead of $\ydP_\calK (M)$. 
\end{definition}

\begin{lemma}
Let $\calJ$ and $\calK$ be categories.
\begin{enumerate}
\item $\ydp_\calJ\circ\ydp_\calK = \ydp_{\calJ\times\calK} = \ydp_\calK\circ\ydp_\calJ $ and
$\ydp_\calJ^\op\circ\ydp_\calK^\op = \ydp_{\calJ\times\calK}^\op = \ydp_\calK^\op\circ\ydp_\calJ^\op $ as functors $\BB\to\BB$.
\item $\ydP_\calJ\circ\ydP_\calK = \ydP_{\calJ\times\calK} = \ydP_\calK\circ\ydP_\calJ $ and
$\ydP_\calJ^\op\circ\ydP_\calK^\op = \ydP_{\calJ\times\calK}^\op = \ydP_\calK^\op\circ\ydP_\calJ^\op $ as functors $\calF\to\calF$.
\end{enumerate}
\end{lemma}

\begin{proof}
This is immediate from the definitions, using  the canonical isomorphisms $\calJ\times\calK\cong \calK\times\calJ$ and $(\calJ\times\calK)^\op\cong \calJ^\op\times\calK^\op$.
\end{proof}

Given a $(\calK,\calJ)$-biset $\Omega$ we obtain a natural transformation $\ydP_\Omega: \ydP_\calJ\to \ydP_\calK$ given at $M$ and $\calC$ by 
$$
\ydP_\calJ M(\calC) = M(\calC\times\calJ) \xrightarrow{M({}_\calC\calC_\calC \times \Omega)}
\ydP_\calK M(\calC) = M(\calC\times\calK)
$$
and we also obtain a natural transformation in the opposite direction $\ydP_\Omega^\op: \ydP_\calK^\op\to \ydP_\calJ^\op$ given at $M$ and $\calC$ by 
$$
\ydP_\calK^\op M(\calC) = M(\calC\times\calK^\op) \xrightarrow{M({}_\calC\calC_\calC \times \tau(\Omega))}
\ydP_\calJ^\op M(\calC) = M(\calC\times\calJ^\op)
$$
using the duality $\tau$. Extending this linearly, if $f$ is an $R$-linear combination of $(\calK,\calJ)$-bisets we obtain natural transformations $\ydP_f: \ydP_\calJ\to \ydP_\calK$ and $\ydP_f^\op: \ydP_\calK^\op\to \ydP_\calJ^\op$.

\begin{corollary}
$\ydP_\bullet : \BB\to \mathrm{Fun}(\calF,\calF)$ and $\ydP_\bullet^\op : \BB^\op\to \mathrm{Fun}(\calF,\calF)$ are functors, where $\mathrm{Fun}(\calF,\calF)$ denotes the category of $R$-linear functors from $\calF$ to $\calF$. 
\end{corollary}

\begin{proof}
We must check that $\ydP_g\circ\ydP_f = \ydP_{gf}$, $\ydP_f^\op\circ\ydP_g^\op = \ydP_{gf}^\op$ and that both $\ydP_\bullet$ and $\ydP_\bullet^\op$ send the identity to the identity. These are immediate.
\end{proof}

\begin{definition}
For each category $\calD$ let $F_\calD$ denote the \textit{representable biset functor} at $\calD$, so $F_\calD(\calC)= \Hom_\BB(\calD,\calC)= A(\calC,\calD)$ is the $R$-linear span of the indecomposable $(\calC,\calD)$-bisets. In particular $F_\one = B$ is the Burnside ring functor.
\end{definition}

\begin{proposition}
\label{Yoneda-Dress-representable-proposition}
We have isomorphisms
$$
\ydP_\calK F_\calD \cong F_{\calK^\op\times\calD}\quad \hbox{and} \quad\ydP_\calK^\op F_\calD \cong F_{\calK\times\calD}
$$
as biset functors. As a particular case, $F_\calD\cong \ydP_\calD^\op B$. 
\end{proposition}

\begin{proof}
    For any category $\calC$ we have
    $$
    (\ydP_\calK F_\calD)(\calC) \cong F_\calD(\calC\times\calK) \cong \Hom_\BB(\calD,\calC\times\calK),
    $$
    using Yoneda's lemma. This is the Grothendieck group of $\calC\times\calK\times \calD^\op$-sets, which is also a description of $F_{\calK^\op\times\calD}(\calC)$. The remaining statements follow.
\end{proof}

\begin{proposition}
\label{left-right-adjoint-proposition}
    For each category $\calK$, the functor $P_\calK^\op:\calF\to\calF$ is both left and right adjoint to $P_\calK$.
\end{proposition}

\begin{proof}
    We show that $P_\calK^\op$ is right adjoint to $P_\calK$. The other adjoint property follows by exchanging the roles of $\calK$ and $\calK^\op$. We construct the unit and counit of the adjunction using the coevaluation and evaluation bisets introduced before Theorem~\ref{biset-dual-theorem}, and verify the triangle identities. We define the unit $\eta$ of the adjunction as the natural transformation $\eta:1\to \ydP_\calK^\op \ydP_\calK$ that, at the biset functor $M$, takes the values
    $$
    \eta_M=M({}_\calC\calC_\calC\times {}_{\calK\times\calK^\op}\calK_\one): M(\calC)\to \ydP_\calK^\op \ydP_\calK M(\calC).
    $$
    The counit $\epsilon:  \ydP_\calK \ydP_\calK^\op\to 1$ at $M$ takes the values
    $$
    \epsilon_M=M({}_\calC\calC_\calC\times {}_\one\calK_{\calK^\op\times\calK}):\ydP_\calK \ydP_\calK^\op M(\calC)\to M(\calC).
    $$
    We claim that the two composites
    $$
    \ydP_\calK\xrightarrow{\ydP_\calK\eta}\ydP_\calK \ydP_\calK^\op\ydP_\calK
    \xrightarrow{\epsilon_{\ydP_\calK}}\ydP_\calK
    $$
    and
    $$
    \ydP_\calK^\op\xrightarrow{\ydP_\calK\eta}\ydP_\calK^\op\ydP_\calK \ydP_\calK^\op
    \xrightarrow{\ydP_\calK^\op\epsilon}\ydP_\calK^\op
    $$
    are both the identity, and this will verify the adjunction. In the case of the first composite this is because, on evaluation at $M$, it gives a natural transformation which, at a category $\calC$, is the morphism $M({}_\calC\calC_\calC\times\Psi)$ followed by $M({}_\calC\calC_\calC\times\Omega)$, which compose to give the identity on $(\ydP_\calK M)(\calC)$, by  the calculation of Theorem~\ref{biset-dual-theorem}. The second composite is the identity by a similar argument.
\end{proof}

\subsection{Structures on biset functors: internal Hom}
The definition of internal Hom in \cite[Sec. 8.3]{Bou3} works verbatim with groups replaced by categories.

\begin{definition}
Let $M$ and $N$ be biset functors. We define the internal Hom of $M$ and $N$ to be
$$
\calH(M,N):=\Hom_\calF (M,\ydP_\bullet(N)),
$$
using the functoriality of the Yoneda-Dress construction in its subscript to make $\calH$ a biset functor.
\end{definition}

Evidently $\calH(-,-)$ is a bilinear functor $\calF\times \calF\to\calF$. We explore some of its properties.

\begin{proposition}
    $\calH (M,N)$ is left exact in each of its variables $M$ and $N$.
\end{proposition}

\begin{proof}
    Functors are (left) exact if and only if so are their evaluations.  In the first variable $\calH$ is left exact because, on evaluation at $\calK$, it is the Hom functor $\Hom_\calF(-,\ydP_\calK(N))$ and $\Hom$ functors are left exact. In the second variable $\calH$ is left exact because, after evaluation at $\calK$, the functor $\ydP_\calK$ is exact, and so $\Hom_\calF(M,\ydP_\calK(-))$  is left exact because Hom is left exact 
\end{proof}

There is an asymmetry in our definition of $\calH$ in that we used functoriality in the second variable of $\Hom$. It is also possible to use functoriality in the first variable of $\Hom$ to get an isomorphic result. We see a difference with the usual version for groups, in that a superscript ${}^\op$ is inserted in the first variable.

\begin{proposition} 
\label{H-first-variable-proposition}
We have
    $\calH(M,N)\cong \Hom_\calF(\ydP_\bullet^\op(M),N)$ as biset functors.
\end{proposition}

\begin{proof}
The proof is simply that, on evaluation at $\calK$, we have
$$
\Hom_\calF(\ydP_\calK^\op(M),N) \cong \Hom_\calF (M,\ydP_\calK(N))
$$
because $\ydP_\calK^\op$ is left adjoint to $\ydP_\calK$ by Proposition~\ref{left-right-adjoint-proposition}.
\end{proof}

\begin{proposition}
\label{yd-to-H-proposition}
    Let $\calK$ be a category. Then
    $$
    \ydP_\calK\circ \calH(-,-) \cong \calH(-,\ydP_\calK(-))
    \cong \calH(\ydP_\calK^\op(-),-).
    $$
\end{proposition}

\begin{proof}
    Evaluating at biset functors $M$ and $N$ and at a category $\calC$ we have
    $$
    \begin{aligned}
    (\ydP_\calK \circ \calH(M,N))(\calC)
    &\cong \calH(M,N)(\calC\times\calK)\cr
    &\cong \Hom_\calF (M,\ydP_{\calC\times\calK}(N))\cr
    &\cong\Hom_\calF (M,\ydP_\calC(\ydP_\calK (N)))\cr
    &\cong\calH(M,\ydP_\calK(N))(\calC)\cr
    \end{aligned}
    $$
    and this establishes the first isomorphism. For the second isomorphism we do a similar thing with the first variable. Starting with the second term of the above isomorphisms, we have (using Proposition~\ref{H-first-variable-proposition})
    $$
    \begin{aligned}
        \calH(M,N)(\calC\times\calK)
        &\cong \Hom_\calF(\ydP_{\calC\times\calK}^\op(M),N)\cr
        &\cong \Hom_\calF(\ydP_\calC^\op(\ydP_\calK^\op(M)), N)\cr
        &\cong\calH(\ydP_\calK^\op(M),N)(\calC),\cr
    \end{aligned}
    $$
    which gives the second isomorphism.
\end{proof}

\begin{corollary}
\label{hom-representable-corollary}
    Let $M$ be a biset functor, let $F_\calJ$ and $F_\calK$ be the representable biset functors at categories $\calJ,\calK$, and let $B$ be the Burnside ring functor. Then
    $$
    \calH(F_\calK,M)\cong \calH(\ydP_\calK^\op(B)M) \cong \ydP_\calK(M)
    $$ 
    In particular, $\calH (B,M)\cong M$ and $\calH(F_\calK,F_\calJ)\cong F_{\calJ\times\calK^\op}$. It follows that if $P$ and $Q$ are projective biset functors, then so is $\calH(P,Q)$.
\end{corollary}

\begin{proof}
    We already know from Proposition~\ref{Yoneda-Dress-representable-proposition} that $F_\calK \cong \ydP_\calK^\op F_\one \cong \ydP_\calK^\op B $ so that by Proposition~\ref{yd-to-H-proposition} it suffices to show that $\calH (F_\one,M)\cong M$. Evaluating the left side of this at $\calC$ we get
    $$
    \calH (F_\one,M)(\calC) \cong \Hom_\calF(F_\one,\ydP_\calC(M))
    \cong \ydP_\calC(M)(\one) \cong M(\one\times\calC) \cong M(\calC),
    $$
    which is what is needed.
    For the penultimate statement we apply Proposition~\ref{Yoneda-Dress-representable-proposition}.
\end{proof}

\subsection{Structures on biset functors: tensor product}
\label{tensor-product-section}

\begin{definition}
    If $L$ is a right biset functor over $R$ (that is, an $R$-linear functor $\BB^\op\to R\mod$) and $M$ is a (left) biset functor over $R$, we will write $L\otimes_\BB M$ for the $R$-module $\bigoplus_{\calD\in\BB} L(\calD)\otimes_R M(\calD)$ modulo relations $(a\cdot {}_\calE\Omega_\calD) \otimes_R b = a\otimes_R ({}_\calE\Omega_\calD \cdot b)$ when $a\in L(\calE)$ and $b\in M(\calD)$ and ${}_\calE\Omega_\calD$ is a $(\calE,\calD)$-biset.
\end{definition}

This construction produces the coend $\int^\calD L(\calD)\otimes_R M(\calD)$ described in \cite{MacLane}.
The notation $\otimes_\BB$ that we propose to use is consistent with the usual notation $L\otimes_\Lambda M$ for the tensor product of a right $\Lambda$-module $L$ and a left $\Lambda$-module $M$, when $\Lambda$ is an $R$-algebra. In our situation $L$ and $M$ may be termed right and left $\BB$-modules, that is, functors $L:\BB^\op\to R\mod$ and $M:\BB\to R\mod$.

In the next definition we use the duality $\tau$ on $\BB$ defined in Definition~\ref{duality-definition} to turn a left biset functor into a right biset functor.

\begin{definition}
    We define the \textit{internal tensor product} of (left) biset functors $L$ and $M$ to be the biset functor
    $$
    L\otimes M:= (\ydP_\bullet(L)\circ\tau) \otimes_\BB M
    $$
    where the functoriality on $\BB$ is obtained via the subscript $\bullet$. Thus the value of $L\otimes M$ at a category $\calC$ is $$
    (L\otimes M)(\calC) = (\ydP_\calC(L)\circ\tau) \otimes_\BB M = L(\tau(-)\times \calC)\otimes_\BB M(-).
    $$
\end{definition}

\begin{theorem}
\label{tensor-adjoint-property-theorem}
    For all biset functors $L,M$ and $N$ we have
    \begin{enumerate}
        \item $\Hom_\calF (L\otimes M,N) \cong \Hom_\calF(M,\calH(L,N))$, and 
        \item $\calH(L\otimes M,N)\cong \calH(M,\calH(L,N))$.
    \end{enumerate}
\end{theorem}

\begin{proof}
    1. We will show that
    $$
    \Hom_\calF ((\ydP_\bullet(L)\circ\tau)\otimes_\BB M,N) \cong \Hom_\calF(M,\Hom_\calF(\ydP_\bullet^\op L,N)).
    $$
    The term on the left consists of natural transformations $\phi$ that are families of $R$-module homomorphisms
$$
\phi_\calC:(\ydP_\calC (L)\circ\tau)\otimes_\BB M \to N(\calC).
$$
The tensor product is spanned by basic tensors $a\otimes_\BB b$ where $a\in L(\calD^\op\times\calC)$ and $b\in M(\calD)$ for some category $\calD$, and each such tensor determines an element $\phi_\calC(a\otimes b)\in N(\calC)$.

From this we may construct a natural transformation
$$
\psi: M\to\Hom_\calF(\ydP_\bullet^\op L,N)
$$
as a family of maps $\psi_\calD : M(\calD)\to \Hom_\calF(\ydP_\calD^\op L,N)$ as follows. If $b\in M(\calD)$ and $a\in \ydP_\calD^\op L(\calC)= L(\calC\times\calD^\op)\cong  L(\calD^\op\times\calC)$, let $\hat a\in L(\calD^\op\times\calC)$ be the element corresponding to $a$ under this isomorphism, determined by the canonical isomorphism $\calC\times\calD^\op\cong  \calD^\op\times\calC$.  We define $\psi_\calD(b)$ to be the natural transformation $\ydP_\calD^\op L\to N$ that is the family of maps
$\psi_\calD(b)_\calC :(\ydP_\calD^\op L)(\calC)\to N(\calC)$ where $\psi_\calD(b)_\calC(a) = \phi_\calC(\hat a\otimes_\BB b)\in N(\calC)$. 

We similarly construct an inverse mapping in the opposite direction from the right to the left of the second line of this proof, by reversing all the constructions just indicated. This establishes 1.

2. We deduce this from part 1. We have 
$$
\begin{aligned}
    \calH(L\otimes M,N) &=\Hom_\calF(L\otimes M, \ydP_\bullet(N))\cr
    &\cong \Hom_\calF (M,\calH(L,\ydP_\bullet(N))\cr
    &=\Hom_\calF(M,\ydP_\bullet\calH(L,M))\cr
    &=\calH(M,\calH(L,N))\cr
\end{aligned}
$$
using Proposition~\ref{yd-to-H-proposition}.
\end{proof}

\begin{corollary}
    For all biset functors $L,M$ the functors $L\otimes -$ and $-\otimes M$ are right exact.
\end{corollary}

\begin{proof}
It is a corollary of Theorem~\ref{tensor-adjoint-property-theorem} that $L\otimes -$ is right exact, because it is the left adjoint of the functor $\calH(L,-)$.

To see that $-\otimes M$ is right exact we employ a separate argument. We show that if $0\to U\to V\to W\to 0$ is an exact sequence of biset functors, then $0\to U\otimes M\to V\otimes M\to W\otimes M\to 0$ is exact, and to test this it suffices to check on each evaluation. The evaluation $U\otimes M(\calC)$ is $(\ydP_\calC U\circ\tau)\otimes_\BB M$ and this is the composite of the functors $\ydP_\calC$ (right adjoint $\ydP_{\calC^\op}$), $\tau$ (right adjoint $\tau$), and $-\otimes_\BB M$. This latter has a right adjoint by the standard argument for tensor products: if $K:\BB^\op\to R\mod$ is a right biset functor, we claim
$$
\Hom_R(K\otimes_\BB M, V) \cong \Hom_\calF(K,\Hom_R(M,R))
$$
via an isomorphism
$$
(\phi: K\otimes M\to V) \mapsto \left(\psi_\calC : K(\calC)\to \Hom(M(\calC),V)\right)
$$
where the natural transformation $\psi = (\psi_\calC)$ is specified by
$$
\psi_\calC(a) = (b\mapsto \phi(a\otimes b)).
$$
There is an inversely defined morphism in the opposite direct that shows we have an isomorphism. It follows from all this that $-\otimes M$ has a right adjoint, and so is right exact.
\end{proof}

\begin{theorem}
\label{tensor-properties-theorem}
For all biset functors $L,M$ and $N$ and categories $\calJ, \calK$ we have
    \begin{enumerate}
        \item $F_\calJ\otimes F_\calK \cong F_{\calJ\times\calK}$,
        \item $L\otimes M\cong M\otimes L$,
        \item $B\otimes M\cong M\cong M\otimes B$ where $B=F_\one$ is the Burnside ring functor, 
        \item $(L\otimes M)\otimes N\cong L\otimes (M\otimes N)$, and
        \item $\ydP_\calJ(L\otimes M)\cong (\ydP_\calJ L)\otimes M\cong L\otimes (\ydP_\calJ M)$.
    \end{enumerate}
\end{theorem}

\begin{proof}
    1. For any biset functor $N$ we have
    $$
    \begin{aligned}
    \calH(F_\calJ\otimes F_\calK, N) &\cong \calH(F_\calK,\calH(F_\calJ,N))\cr 
    &\cong \calH(F_\calK,\ydP_\calJ(N)) \cr 
    &\cong \ydP_\calJ\calH(F_\calK,N) \cr 
    &\cong \ydP_\calJ \ydP_\calK (N) \cr 
    &\cong \ydP_\calK \ydP_\calJ (N) \cr 
    &\cong \ydP_{\calJ\times\calK} (N) \cr 
    &\cong  \calH(F_{\calJ\times\calK},N). \cr 
    \end{aligned}
    $$
    Noting that $\calH(M,N)(\one)=\Hom_\calF(M,N)$ we deduce that 
    $$
    \Hom_\calF(F_\calJ\otimes F_\calK, N)\cong\Hom_\calF(F_{\calJ\times\calK}, N)\cong N(\calJ\times\calK),
    $$
    so that $\Hom_\calF(F_\calJ\otimes F_\calK, -)$ is isomorphic to the representable functor $\Hom_\calF(F_{\calJ\times\calK}, -)$. Because the Yoneda embedding is faithful we have $F_\calJ\otimes F_\calK\cong F_{\calJ\times\calK}$.

    2. When $L$ and $M$ happen to be representable functors $L=F_\calJ$ and $M=F_\calK$ the result holds by part 1. From this the result holds when $L$ and $M$ are projective functors.  For general $L$ and $M$ we employ an argument with projective presentations. Let $P_1\to P_0 \to L\to 0$ and $Q_1\to Q_0\to M\to 0$ be exact sequences, where the $P_i$ and $Q_i$ are projective. Now $L$ and $M$ are the zero homology of the complexes $P_1\to P_0$ and $Q_1\to Q_0$ and, using right exactness on both sides of $\otimes$, $L\otimes M$ is the zero homology of $(P_1\to P_0)\otimes (Q_1\to Q_0)$. This complex is isomorphic to $(Q_1\to Q_0)\otimes (P_1\to P_0)$ because $\otimes$ is commutative on projective functors, and this has zero homology $M\otimes L$, which is thus isomorphic to $L\otimes M$.

    3. When $M$ is a representable functor $M=F_\calK$ then $F_\one\otimes F_\calK\cong F_{\one\times\calK}\cong M\cong F_\calK \otimes F_\one$ by parts 1 and 2. Thus $B\otimes -$ is the identity functor on projective functors.  If $M$ is arbitrary we resolve it by projectives $P_1\to P_0\to M\to 0$ and apply the right exact functor $B\otimes -$ to get an exact sequence $ P_1\to  P_0\to B\otimes M\to 0$ showing that $B\otimes M\cong M$, and similarly with $M\otimes B\cong M$.

4. For any biset functor $Q$ we have 
$$
\begin{aligned}
    \Hom_\calF((L\otimes M)\otimes N,Q)&\cong \Hom_\calF(N,\calH(L\otimes M, Q))\cr
    &\cong \Hom_\calF(N,\calH(M, \calH(L,Q)))\cr
    &\cong \Hom_\calF(M\otimes N,\calH(L,Q))\cr
    &\cong \Hom_\calF(L\otimes(M\otimes N),Q)\cr
\end{aligned}
$$
so that the two representable functors $\Hom_\calF((L\otimes M)\otimes N,-)$ and $\Hom_\calF(L\otimes(M\otimes N),-)$ are isomorphic. Because the Yoneda embedding is faithful the result follows.

5. Using Propositions~\ref{H-first-variable-proposition}, \ref{yd-to-H-proposition} and Theorem~\ref{tensor-adjoint-property-theorem} we have
$$
\begin{aligned}
    \Hom_\calF(\ydP_\calJ(L\otimes M),N)
    &\cong \Hom_\calF(L\otimes M,\ydP_\calJ^\op N)\cr
    &\cong \Hom_\calF(M,\calH(L,\ydP_\calJ^\op N))\cr
    &\cong \Hom_\calF(M,\calH(\ydP_\calJ L, N))\cr
    &\cong \Hom_\calF((\ydP_\calJ L)\otimes M,N)\cr
    &\cong \Hom_\calF(L\otimes (\ydP_\calJ M),N) \quad\hbox{similarly.}\cr
\end{aligned}
$$
The representable functors $\Hom_\calF(\ydP_\calJ(L\otimes M),-)$, $\Hom_\calF((\ydP_\calJ L)\otimes M,-)$ and $\Hom_\calF(L\otimes (\ydP_\calJ M),-)$ are all isomorphic, so the result follows because the Yoneda embedding is faithful.
\end{proof}

To summarize:

\begin{corollary}
\label{symmetric-monoidal-corollary}
With the product operation $\otimes$ and unit $B$ the category of biset functors $\calF$ is a symmetric monoidal tensor category. Furthermore, if $P$ and $Q$ are projective biset functors then $P\otimes Q$ and $\calH(P,Q)$ are also projective.   
\end{corollary}

\begin{proof}
    The last sentence is a consequence of Theorem~\ref{tensor-properties-theorem}(1) and Corollary~\ref{hom-representable-corollary}.
\end{proof}

\subsection{Structures on biset functors: a contravariant equivalence}

For each biset functor $M$ we define another biset functor $M^*$ as $M^*(\calC):= \Hom_R(M(\tau(\calC)),R)$, by pre-composing $M$ with the contravariant functor $\tau$ and post-composing with the contravariant functor $\Hom_R(-,R)$. Then $M^*$ is a \textit{dual biset functor} of $M$ in a certain sense, not to be confused with the notation of a dual object in the context of tensor categories.  When $M(\calC)$ and $M(\calC^\op)$ are finitely generated projective $R$-modules, we have  $M^{**}\cong M$. 

\subsection{Green biset functors and modules}
\label{Green-section}

Green biset functors are biset functors with a product and unit that make them behave like rings.
The development of their theory given in section 8.5 of \cite{Bou3}, where they are defined on groups, works also for biset functors defined on categories. We sketch briefly how this goes, referring to \cite{Bou3} for greater detail.

\begin{definition}
    A Green biset functor is a monoid $A$ in $\calF$.
\end{definition}

This means there are maps of biset functors
$$
\mu: A\otimes A\to A, \quad e:B\to A 
$$
where the Burnside ring functor $B$ is the unit of $\otimes$.
There is a requirement that the following two diagrams commute, one ensuring that $\mu$ is associative, and the other that $B$ multiplies as the identity:
\begin{center}
\scalebox{0.9}{
$
\diagram{
{A\otimes(A\otimes A)} &\umapright{\mathrm{Id}\otimes\mu}& {A\otimes A}&&\cr
&&&\mu\atop\searrow&\cr
\rmapdown{\alpha}&&&&A\cr
&&&\mu\atop\nearrow&\cr
{(A\otimes A)\otimes A} &\umapright{\mu\otimes \mathrm{Id}}& {A\otimes A}&&\cr
}
\quad
\diagram{
{B\otimes A} &\umapright{e\otimes \mathrm{Id}}& {A\otimes A} &\umapleft{\mathrm{Id}\otimes e}& {A\otimes B}\cr
&\lambda\atop\searrow&\rmapdown{\mu}&\rho\atop\swarrow&\cr
&&A&&\cr
}
$
}
\end{center}
where $\alpha, \lambda$ and $\rho$ are the isomorphisms that are part of the symmetric monoidal structure.

By the adjoint property of $\otimes$ and $\calH$ (Theorem~\ref{tensor-adjoint-property-theorem}), the map $\mu$ corresponds to a natural transformation of biset functors $\nu:A\to\calH(A,A)$. It means that for each category $\calD$ there is a map of $R$-modules
$$
\nu_\calD: A(\calD)\to\calH(A,A)(\calD) = \Hom_\calF(A,\ydP_\calD A)
$$
so that, for each element $b\in A(\calD)$ and for each category $\calC$, there is a map of $R$-modules
$$
\nu_\calD(b)_\calC: A(\calC)\to A(\calC\times\calD).
$$
Putting this together, $\mu$ determines for each pair of categories $\calC,\calD$ a bilinear map of $R$-modules
$$
A(\calC)\times A(\calD) \to A(\calC\times\calD)
$$
that sends each pair of elements $(a,b)$ to $\nu_\calD(b)_\calC(a)$.

\begin{definition}
    Let $a\in A(\calC)$ and $b\in A(\calD)$. We define
    $$
    a\times b:=\nu_\calD(b)_\calC(a)\in A(\calC\times\calD)
    $$
We also define $\varepsilon_A\in A(\one)$ to be the image of the single point in $B(\one)$ under the map $e$ (at $\one$).
\end{definition}

This specification satisfies the following conditions:

\begin{itemize}
\item (Associativity) Let $\calC,\calD$ and $\calE$ be categories and $$
\alpha: \calC\times(\calD\times\calE) \to (\calC\times\calD)\times\calE
$$
the canonical isomorphism functor. Let $\Omega_\alpha$ be the biset for these two product categories determined by $\alpha$ (as in Section~\ref{functors-to-bisets-section}). Then if $a\in A(\calC)$, $b\in A(\calD)$ and $c\in A(\calE)$
$$
(a\times b)\times c = A(\Omega_\alpha)(a\times (b\times c)).
$$
\item (Identity element) For each category $\calC$ Let $\lambda_\calC: \one\times\calC\to\calC$ and $\rho_\calC:\calC\times\one\to\calC$ denote the canonical isomorphism functors, and $\Omega_{\lambda_\calC},\Omega_{\rho_\calC}$ the corresponding $(\calC,\one\times\calC)$- and $(\one\times\calC,\calC)$-bisets. Then for every $a\in A(\calC)$
$$
a= A(\Omega_{\lambda_\calC})(\varepsilon_A\times a) = A(\Omega_{\rho_\calC})(a\times\varepsilon_A).
$$
\item (Functoriality) If ${}_{\calC'}\Omega_\calC$ and ${}_{\calD'}\Psi_\calD$ are bisets, then for any $a\in A(\calC)$ and $b\in A(\calD)$
$$
A(\Omega\times\Psi)(a\times b) = A(\Omega)(a) \times A(\Psi)(b).
$$
\end{itemize}

\begin{proposition}
    A biset functor $A$ is a Green functor if and only if there are bilinear maps of $R$-modules $A(\calC)\times A(\calD)\to A(\calC\times\calD)$ and an element $\varepsilon_A\in A(\one)$ satisfying the conditions just listed.
\end{proposition}

\begin{proof}
    We leave the verification of this discussion to the reader.
\end{proof}

\begin{definition}
    Let $A$ be a Green biset functor. A (left) $A$-module $M$ is a biset functor equipped with a morphism of biset functors $\mu_M:A\otimes M\to M$, with the requirement that a certain two diagrams commute imposing associativity and the unital property. See  \cite{Bou3} for details.
\end{definition}

Again, it is equivalent to require that there are bilinear maps 
$$
A(\calC)\times M(\calD)\to M(\calC\times\calD)
$$
satisfying conditions described in \cite{Bou3}.

\begin{example}
As with Green biset functors defined on groups, the Burnside functor $B$ is a Green biset functor. The map $B(\calC)\times B(\calD)\to B(\calC\times\calD)$ sends a pair $(\Omega,\Psi)$ consisting of a $\calC$-set and a $\calD$-set to the $(\calC\times\calD)$-set $\Omega\times\Psi$. Not only that, but every biset functor has the structure of a $B$-module, where the maps $B(\calC)\times M(\calD)\to M(\calC\times\calD)$ are specified as follow: if ${}_\calC\Omega_\one$ is a $\calC$-set in $B(\calC)$, regarded as a $(\calC,\one)$-biset, and $u\in M(\calD)\cong M(\one\times\calD)$, then the pair $(\Omega,u)$ is sent to $M(\Omega\times {}_\calD\calD_\calD)(u)\in M(\calC\times\calD)$. Thus biset functors are the same thing as modules for the Green biset functor $B$.
\end{example}

There is a very useful category $\calP_A$ associated to a Green biset functor $A$ with the property that $A$-modules may be identified as $R$-linear functors $\calP_A\to R\mod$. In the context of groups it appears as part (5) of \cite[Prop. 8.6.1]{Bou3} and is described more fully in \cite[Sec. 2.1]{Rom}, where it plays an important role in describing simple $A$-modules. The description given there also works for categories in general, not just groups, provided we distinguish appropriately between a category and its opposite. 

\begin{definition}
    Given a Green biset functor $A$, we define a category $\calP_A$ whose objects are finite categories, and if $\calC,\calD$ are finite categories then
    $$
    \Hom_{\calP_A}(\calC,\calD) = A(\calD\times\calC^\op).
    $$
    The composition of $\beta\in A(\calE\times\calD^\op)$  and $\alpha\in A(\calD\times\calC^\op)$ is the following:
    $$
    \beta\circ\alpha = A({}_\calE\calE_\calE\times{}_\one\calD_{\calD^\op\times\calD}\times{}_{\calC^\op}\calC^\op_{\calC^\op})(\beta\times\alpha).
    $$
    Here $\beta\times\alpha$ denotes the image of $(\beta,\alpha)$ under the map 
    $$
    A(\calE\times\calD^\op)\times A(\calD\times\calC^\op) \to A(\calE\times\calD^\op\times \calD\times\calC^\op)
    $$
    that is part of the Green biset functor.
\end{definition}

\begin{proposition}
    Given a Green biset functor $A$, the specification of $\calP_A$ defines an $R$-linear category. At each category $\calC$ the identity morphism of $\calC$ in $\calP_A$ is equal to $A({}_{\calC\times\calC^\op}\calC_\one)(\varepsilon_A)\in A(\calC\times\calC^\op)$. It has the property that the category of $A$-modules is equivalent to the category of $R$-linear functors $\calP_A\to R\mod$. That is, $A\mod\simeq\Fun(\calP_A,R\mod)$. 
\end{proposition}

\begin{proof}
    The proof is the same as that of Proposition 2.11 of \cite{Rom}, writing categories instead of groups, and making use of opposite categories as appropriate. The dictionary that translates \cite{Rom} to the situation of categories is as follows. For each biset ${}_\calD\Omega_\calC$ the notation ${}_\calD\overrightarrow \Omega_\calC$ is replaced by ${}_{\calD\times\calC^\op} \Omega_\one$, and ${}_\calD\overleftarrow \Omega_\calC$ is replaced by ${}_\one \Omega_{\calD^\op\times\calC}$. Groups $G,H$ should be replaced by categories $\calC,\calD$ and a product $G\times H$ is replaced by $\calC\times\calD^\op$. There is one place in \cite{Rom} where a biset $X_H=(H\times H)/\Delta(H)$ appears, and this should be translated as ${}_{\calD\times\calD^\op}\calD_\one$.
    
    To give the essence of the proof, we describe the functors in the equivalence. We define a functor  $A\mod\to\Fun(\calP_A,R\mod)$ as follows. Given an $A$-module $M$ we construct an $R$-linear functor $F_M:\calP_A\to R\mod$ by $F_M(\calC)=M(\calC)$ on objects; for morphisms $\alpha\in A(\calD\times\calC^\op)$ we define $F_M(\alpha):M(\calC)\to M(\calD)$ by
    $$
    F_M(\alpha)(m) = M({}_\calD\calD_\calD\times {}_\one\calC_{\calC^\op\times\calC})(\alpha\times m).
    $$

    In the opposite direction, if $F:\calP_A\to R\mod$ is an $R$-linear functor we define an $A$-module $M$ as follows. On objects, $M(\calC)=F(\calC)$. On a biset ${}_\calD\Omega_\calC$ we use the unit $e:B\to A$ to obtain an element $e({}_{\calD\times\calC^\op}\Omega_\one) \in A(\calD\times\calC^\op)$ and now $M({}_\calD\Omega_\calC): M(\calC)\to M(\calD)$ is defined to be $F(e({}_{\calD\times\calC^\op}\Omega_\one))$. This establishes $M$ as a biset functor, and to make $M$ an $A$-module we define $A(\calC)\times M(\calD)\to M(\calC\times\calD)$
    by $(a,m)\mapsto F(A({}_\calC\calC_\calC\times {}_{\calD\times\calD^\op}\calD_\one)(a))(m)$. The verification that these are functors establishing an equivalence follows the lines of \cite[Prop. 2.11]{Rom}.
\end{proof}

\subsection{Functors to categories other than $\Set$}
We have developed a theory of functors with values in $\Set$ because this is the classical situation and the questions that arise occur at their most basic level. We may, more generally, take a symmetric monoidal category $\SS$ with product $\diamond$, and with the property that finite colimits exist and commute with $-\diamond X$ for all objects $X$ in $\SS$. We will write $X\sqcup Y$ for the coproduct of $X$ and $Y$ in $\SS$. There is a monoidal functor $\sigma: \Set\to\SS$ that sends a set $X$ to $\sigma(X):=\bigsqcup_{x\in X} \one_\SS $, the coproduct of copies of the unit $\one_\SS$ indexed by $X$, and sends a morphism of sets to the induced map on the coproducts. 

\begin{definition}
    For each category $\calC$ we define a \textit{$\calC$-object} in $\SS$ to be a functor $\calC\to\SS$, and a \textit{$(\calC,\calD)$-bi-object in $\SS$} to be a functor $\calC\times\calD^\op\to\SS$.  
\end{definition} 

We can now form the Grothendieck group $A^\SS(\calD,\calC)$ of $(\calD,\calC)$-bi-objects with respect to $\sqcup$, extend the scalars to $R$ giving $A_R^\SS(\calD,\calC)$, and construct the \textit{bi-object category $\BB^\SS$} whose objects are finite categories, and where $\Hom(\calC,\calD) = A_R^\SS(\calD,\calC)$. In this category, the composition of two bi-objects is defined by the same formula as in Definition~\ref{biset-composition-definition}, with the direct product of sets replaced by the $\diamond$ product of objects of $\SS$. It is associative. The unit bi-object at a category $\calC$ is the composite $\sigma\circ {}_\calC\calC_\calC :\calC\times\calC^\op \to \SS$. 

\begin{definition}
\textit{Bi-object functors with respect to $\SS$} are defined to be $R$-linear functors $\BB^\SS\to R\mod$. These are the objects of a category $\calF^\SS$ in which the morphisms are natural transformations.
\end{definition}

The functor $\sigma:\Set\to\SS$ gives rise to an $R$-linear functor $\Sigma:\BB\to\BB^\SS$ that is the identity on objects (which are categories in both cases) and which sends a morphism $\Omega:\calD\times\calC^\op\to\Set$ to the composite $\Sigma(\Omega):=\sigma\circ\Omega$. Composition with $\Sigma$ now provides a functor $\calF^\SS\to\calF$, so that bi-object functors with respect to $\SS$ may be regarded, by restriction along $\Sigma$ as biset functors.

All of the categories we consider have monoidal structures. On categories $\calC,\calD$ the tensor product in $\BB$ and in $\BB^\SS$ is given by direct product $\calC\times\calD$. Using this structure we may define Yoneda-Dress functors $\ydP_\calC$ on $\calF^\SS$ and hence a symmetric monoidal structure on $\calF^\SS$ in the same way as in Section~\ref{tensor-product-section}. From this we may construct \textit{Green bi-object functors with respect to $\SS$} in the same way as in Section~\ref{Green-section}. We leave the details to the interested reader, and merely summarize some of the key points.

\begin{proposition}
Let $\SS$ be a symmetric monoidal category with
product $\diamond$, and with the property that finite colimits exist and commute with $\diamond$.
\begin{enumerate}
    \item The bi-object category $\BB^\SS$ is a category. 
    \item The specification of $\sigma$ defines a functor.
    \item The functors $\sigma$ and $\Sigma$, as well as the  restriction $\calF^\SS\to\calF$ along $\Sigma$, are all monoidal. 
    \item Every bi-object functor with respect to $\SS$ can be regarded also as a biset functor and every Green bi-object functor with respect to $\SS$ can be regarded as a Green functor.
\end{enumerate}
    \end{proposition}

\subsection{Fibered biset functors}

As an example of the use of bi-object functors with respect to a category $\SS$ other than $\Set$ we show how a (modified) theory of fibered biset functors may be extended to categories.
The theory of $A$-fibered biset-functors was defined in \cite{BC}, where $A$ is an abelian group. For groups $G,H$, an $A$-fibered $(H,G)$-biset is an $(H,G)$-biset with a commuting action of $A$, such that the orbits under the action of $A$ are free. The Grothendieck over $R$ of $(H,G)$-bisets with respect to disjoint union is denoted $B_R^A(H,G)$ in \cite{BC}, and the fibered biset category, which is denoted $\calC_R^A$ there, has finite groups as its objects and $\Hom(G,H) = B_R^A(H,G)$. Now $A$-fibered biset functors are $R$-linear functors $\calC_R^A\to R\mod$.

In order to realize this for categories, a first thought is to take $\SS$ to be category whose objects are free $A$-sets, with $A$-equivariant maps as morphisms, and a certain product $\diamond$. The trouble with this is that finite colimits do not exist in this category: the colimit of a diagram of $A$-sets with free orbits need not have free orbits, as remarked in \cite{BC}. The solution there is simply to remove any orbits that are not free after taking the colimit, using the fact that commuting group actions of $G$ and $H$ preserve free orbits. In the context of categories, where the action of a morphism in a category on orbits under a commuting action need not be invertible, it is no longer always the case that free $A$-orbits are sent to free $A$-orbits, and we have to be more careful.
The solution to this issue that we provide here is to start by allowing all $A$-orbits, not just the free ones. 

\begin{definition}
    Let $A$ be an abelian group and let $\SS$ be the category of $A$-sets with finitely many $A$-orbits. 
If $\calC, \calD$ are categories then we will call a $\calC$-object in $\SS$, namely a functor $\Omega:\calC\to\SS$, a \textit{weakly $A$-fibered $\calC$-set}. A $(\calC,\calD)$-bi-object in $\SS$, namely a functor $\Omega:\calC\times\calD^\op\to\SS$, is a \textit{weakly $A$-fibered $(\calC,\calD)$-biset}. 
\end{definition}


The category $\SS$ just defined does possess finite colimits (take the colimit of the diagram of sets, with the induced $A$-action). We define a product on $A$-sets $X$ and $Y$ by $X\diamond Y:= (X\times Y)/\sim$ where the equivalence relation $\sim$ is specified by $(ax,y)\sim (x,ay)$ for all $x\in X$, $y\in Y$ and $a\in A$. Thus $X\diamond Y$  is the set of orbits of the subgroup $\{(a,a^{-1}) \bigm| a\in A\}\le A\times A$ on $X\times Y$, including the orbits that are not free. If $[x,y]$ denotes the equivalence class of $(x,y)$ then $a[x,y] = [ax,y]=[x,ay]$ defines an action of $A$ on $X\diamond Y$. With this definition, $\SS$ is a symmetric monoidal category, where the unit is the set $A$ itself, permuted regularly.

\begin{definition}
    We call the bi-object category $\BB^\SS$ in this situation the \textit{weakly $A$-fibered biset category}, and $R$-linear functors $\BB^\SS\to R\mod$ \textit{weakly $A$-fibered biset functors}.
\end{definition}

 \begin{example}
      If the categories $\calC, \calD, \calE$ are groups and ${}_\calE\Psi_\calD$, ${}_\calD\Omega_\calC$ are  bi-objects in $\SS$ with free $A$-orbits, then forming the product $\Psi\circ\Omega$ and subsequently removing the non-free $A$-orbits we obtain the product of $A$-fibered bisets defined in \cite{BC}.
 \end{example}

 \begin{proposition}
Let $A$ be an abelian group, let $\BB^\SS(\hbox{Groups})$ be the full subcategory of the weakly $A$-fibered biset category whose objects are groups, and let the categories $\calC,\calD$ be groups. Let $I(\calC,\calD)$ be the span in $\Hom_{\BB^\SS}(\calC,\calD)$ of the indecomposable weakly $A$-fibered $(\calD,\calC)$-bisets $\Omega$ that have a non-free $A$-orbit. Then 
\begin{enumerate}
    \item $I$ is an ideal in $\BB^\SS(\hbox{Groups})$, and
    \item $A$-fibered biset functors on groups may be identified as the weakly $A$-fibered functors on groups that vanish on the ideal $I$.
\end{enumerate}
\end{proposition}

 To say that $I$ is an ideal means that if the $A$-fibered biset $\Omega$ has a non-free $A$-orbit, then so do $\Psi\circ\Omega$ and $\Omega\circ\Theta$, whenever $\Psi,\Theta$ are weakly $A$-fibered bisets for which the composite is defined.

\begin{proof}
The proof in \cite{BC} that $A$-fibered bisets with free $A$-orbits span the space of morphisms in a category is equivalent to showing that $I$ is an ideal. We see that the quotient category $\BB^\SS(\hbox{Groups})/I$ is equivalent to the $A$-fibered biset category of \cite{BC}, and from this statement (2) follows.
\end{proof}

The quotient category in the last proof has the same objects as $\BB^\SS(\hbox{Groups})$ and
$$
\Hom_{\BB^\SS(\mathrm{Groups})/I}(\calC,\calD):=\Hom_{\BB^\SS(\mathrm{Groups})}(\calC,\calD)/I(\calC,\calD),
$$
which identifies as the span of indecomposable $(\calD,\calC)$-bisets with free $A$-orbits.
The simple $A$-fibered biset functors have been studied extensively in \cite{BC}, \cite{Rom} and elsewhere. The next corollary shows that this is part of the study of weakly $A$-fibered biset functors defined on all categories.

\begin{corollary}
    Simple $A$-fibered biset functors defined on groups are also simple weakly $A$-fibered biset functors defined on groups. They extend uniquely to simple weakly $A$-fibered biset functors defined on all categories.
\end{corollary}

\begin{proof}
    The extension to all categories is a consequence of Proposition~\ref{simple-restriction-extension}.
\end{proof}

\end{document}